\documentclass{amsart}

\newtheorem{theorem}{Theorem}[section]
\newtheorem{lemma}[theorem]{Lemma}
\newtheorem{corollary}[theorem]{Corollary}
\newtheorem{proposition}[theorem]{Proposition}

\theoremstyle{definition}
\newtheorem{definition}[theorem]{Definition}
\newtheorem{example}[theorem]{Example}

\theoremstyle{remark}

\numberwithin{equation}{section}

\begin{document}

\title{Right coideal subalgebras in $U_q(\mathfrak{sl}_{n+1})$}
\author{V.K. Kharchenko, A.V. Lara Sagahon}
\address{FES-Cuautitl\'an, Universidad Nacional Aut\'onoma de M\'exico, 
Centro de Investigaciones Te\'oricas, 
Primero de Mayo s/n, Campo 1, Cuautitl\'an Izcalli, 
Edstado de M\'exico, 54768, M\'EXICO}
\email{vlad@servidor.unam.mx}


\date{}


Journal of Algebra, 319(2008), 2571--2625

\begin{abstract}
\small
We offer a complete classification of right coideal subalgebras which contain 
all group-like elements for the multiparameter version of the quantum group 
$U_q(\mathfrak{sl}_{n+1})$ provided that the main parameter $q$ is not a root of 1. 
As a consequence, we determine that for each subgroup $\Sigma $ of the group $G$
of all group-like elements  
the quantum Borel  subalgebra $U_q^+(\mathfrak{sl}_{n+1})$ containes  $(n+1)!$ different 
homogeneous right coideal subalgebras $U$ such that $U\cap G=\Sigma .$
If $q$ has a finite multiplicative order $t>2,$ the classification remains valid for homogeneous right 
coideal subalgebras of the multiparameter version of the Lusztig quantum group 
$u_q(\frak{ sl}_{n+1}).$ In the paper we consider the quantifications of Kac-Moody algebras as 
character Hopf algebras  [V.K. Kharchenko, A combinatorial approach to the quantifications of Lie algebras, Pacific J. Math., 203(1)(2002), 191-233].
\end{abstract}

\maketitle

\section{Introduction}
One of the reasons one-sided coideal subalgebras became more and more important 
is that Hopf algebras do not have ``enough" Hopf subalgebras.
The one-sided coideal condition instead plays
prominent roles in constructions and developing theory.
The very one-sided comodule subalgebras, but not the Hopf subalgebras, turn out to be the Galois 
objects in the Galois theory for Hopf algebra actions (A. Milinski \cite{Mil1, Mil2}, see also a detailed  
survey by T.Yanai \cite{Y}). In particular, the Galois correspondence theorem for the actions on free 
algebra set up a one to one correspondence between right coideal subalgebras and
intermediate free subalgebras (see,  V.O. Ferreira, L.S.I. Murakami, and A. Paques \cite{FMP}).
A recent survey by G. Letzter \cite{Let} provides a panorama of the use
of one-sided coideal subalgebras in constructing quantum symmetric pairs, in forming quantum
Harish-Chandra modules, and in producing quantum symmetric spaces
(T. Koornwinder \cite{Koo},  A. Joseph and G. Letzter \cite{JL},
M. Noumi and T. Sugitani \cite{NS}, M. Noumi \cite{Nou},
M. Dijkhuizen \cite{Dij}, M.S. K\'eb\'e \cite{Keb, Keb1},
G. Letzter \cite{Let1, Let2, Let3}, S. Sinelshchikov and L. Vaksman \cite{SV}, 
M. Dijkhuizen and M. Noumi \cite{DN}).

In the present paper we offer a complete classification of right coideal subalgebras which contain the 
coradical for the multiparameter version of the quantum group $U_q(\frak{sl}_{n+1}),$
see \cite{Res,CV,CM}, provided that the main parameter $q$ is not a root of 1. If $q$ has 
a finite multiplicative order $t>2,$ this classification remains valid for homogeneous right coideal 
subalgebras of the multiparameter version, see \cite{Luz1, Tow}, of the Lusztig quantum group 
$u_q(\frak{ sl}_{n+1}).$ We are reminded that any Hopf algebra generated by group-like and skew-
primitive elements is pointed, while in a pointed Hopf algebra the group-like elements span
the coradical, see \cite[Definition 5.1.5]{Mon}.

In the second section we introduce main concepts and provide the general results 
on the structure of the character Hopf algebras that are of use for classification. 
In Lemma \ref{qsim} we note that 
if the given character Hopf algebra $H=A\# {\bf k}[G]$ is a bosonisation
of a quantum symmetric algebra $A,$ then each invariant 
differential subspace $U$ of $A$ defines a right coideal $U\# {\bf k}[G].$
This statement allows one to use noncommutative differential calculus, 
\cite[p.6]{Luz2}, \cite{MS}, \cite{Kh1},
due to P. Schauenburg's characterization of quantum Borel subalgebras \cite{Scha}.
The key point of the section is the construction of a PBW-basis
over the coradical for a right coideal 
subalgebra by means of  \cite{KhT, KhQ}. This basis, in particular, provides some 
invariants for right coideal subalgebras (Definition \ref{root}). 

In the third section we define the multiparameter quantification of a Kac-Moody algebra as a character 
Hopf algebra. This approach \cite{Kh4} combines and generalizes all known quantifications. We do 
not put unnecessary restrictions on the characteristic and on the quantification parameters. This allows 
one, for example, to define a new class of finite Frobenius algebras as the Lusztig quantum groups 
over a finite field. All their right coideal subalgebras are also Frobenius \cite{Skr} (finite Frobenius 
algebras, in turn, have a significant r\^{o}le in the coding theory \cite{GMO}). 
In Proposition \ref{rtri} we provide 
a short proof of the so called  ``triangular decomposition" in a quite general form.

In the fourth section (Proposition \ref{phi}) we show that each homogeneous right coideal subalgebra
in the quantum Borel algebra $U_q^+(\frak{ sl}_{n+1})$ has PBW-generators over {\bf k}$\, [G]$
of the following form
\begin{equation} \Psi ^{\hbox{\bf s}}(k,m)=
 \hbox{\Large [[}\ldots \hbox{\Large [}u[1+s_r,m], u[1+s_{r-1},s_r]\hbox{\Large ]}, \ldots \ ,
u[1+s_1,s_2]\hbox{\Large ]},\, u[k,s_1]\hbox{\Large ]},
\label{cbr1int}
\end{equation}
where the brackets are defined by the structure 
of a character Hopf algebra, $[u,v]=uv-\chi ^u(g_v)vu;$ 
$u[i,j]=[\ldots [x_i,x_{i+1}],\ldots ,x_j];$ $k\leq s_1<s_2<\ldots <s_r<m,$
 ${\bf S}\cap [k,m-1]=\{ s_1,s_2,\ldots , s_r\},$
while $x_i,$ $1\leq i \leq n$ are the main skew-primitive generators of $U_q^+(\frak{ sl}_{n+1}).$
This certainly implies that the set of all right coideal subalgebras which contain
the coradical is finite (Corollary \ref{fin1}). To precisely describe the right coideal
subalgebras  we associate a right coideal subalgebra ${\bf U}_{\theta }$ to 
each sequence of integer  numbers $\theta =(\theta _1, \theta _2, \ldots , \theta _n),$
$0 \leq \theta _i\leq n-i+1,$ $1\leq i\leq n$ in the following way.

We define  subsets $R_k,$ $T_k,$ $1\leq k\leq n$
of the interval $[k, n]$ thus: if $\theta _n=0,$ we put $R_n=T_n=\emptyset $ and if 
$\theta _n=1,$ we put $R_n=T_n=\{ n \}.$ Suppose that 
$R_i,$ $T_i,$ $k<i\leq n$ are already defined. If $\theta _k=0,$ then we set 
$R_k=T_k=\emptyset .$ If $\theta _k\neq 0,$ then 
by definition $R_k$ contains 
$\tilde{\theta }_k=k+\theta _k-1$
and all $m$ that satisfy the following three properties
$$
\begin{matrix} \smallskip
a)\ k\leq m<\tilde{\theta }_k; \hfill \cr \smallskip 
b)\ \tilde{\theta }_k \notin T_{m+1}; \hfill \cr
c) \  \forall r (k\leq r<m)\ \ m\in T_{r+1}\Longleftrightarrow \tilde{\theta }_k\in T_{r+1}. \hfill
\end{matrix}
$$
Respectively, $T_k\stackrel{df}{=}R_k\cup \bigcup_{s\in R_k\setminus \{ n\} } T_{s+1}.$
The algebra ${\bf U}_{\theta }$ by definition is
generated over the coradical  by all $\Psi ^{T_k}(k,m),$ $1\leq k\leq m\leq n$ with 
$m\in R_k.$

Theorem \ref{teor} shows that if $q$ is not a root of 1 then
all right coideal subalgebras over the coradical
have the form ${\bf U}_{\theta }.$ In particular, the exact number of
right coideal subalgebras, which include the coradical,
 in the quantum Borel algebra $U_q^{\pm }(\frak{ sl}_{n+1})$ equals $(n+1)!.$ 
If $q$ has a finite multiplicative order $t>2,$ then this is the case for homogeneous
right coideal subalgebras of $u_q^{\pm }(\frak{ sl}_{n+1}).$ (If $q$ is not a root of 1 then
all right coideal subalgebras that contain $G$ are homogeneous, Corollary \ref{odn1}).  

In Section 6 we consider right coideal subalgebras
in the quantum Borel algebra that do not contain the coradical. 
Note that for every submonoid 
$\Omega \subseteq G$ the set of all linear combinations {\bf k}$\, [\Omega]$
is a right coideal subalgebra. We show that if the intersection $\Omega $ of a homogeneous 
right coideal subalgebra $U$ with $G$ is a subgroup, then
$U=\, ${\bf U}$_{\theta }^{1}\, ${\bf k }$[\Omega ].$ Here {\bf U}$_{\theta }^{1}$
is a subalgebra generated by $g^{-1}_aa$ when $a$ runs through the described above 
generators of {\bf U}$_{\theta }.$

In Section 7 we characterize ad$_r$-invariant right coideal subalgebras that have 
trivial intersection with the coradical in terms of   K\'eb\'e's construction \cite{Keb, Keb1}. 

We see that the construction of {\bf U}$_{\theta }$ is completely
constructive, although it is not straightforward. 
Hence by means of computer calculations one may find 
all necessary invariants of the coideal subalgebras and relations between them.
In the eighth section  we provide a tableaux of the coideal subalgebras and their main
characteristics for $n=3$ that was found by means of computer calculations.

In Sections 9-11 we consider the whole of $U_q(\frak{ sl}_{n+1}).$
The triangular decomposition, 
\begin{equation}
U_q(\frak{ sl}_{n+1})= U_q^-(\frak{ sl}_{n+1})\otimes _{{\bf k}[F]} {\bf k}[H]
\otimes _{{\bf k}[G]}U_q^+(\frak{ sl}_{n+1}),
\label{trint}
\end{equation}
provides a hope that any (homogeneous) right coideal subalgebra that contains the coradical
has the triangular decomposition as well, and for any two right coideal subalgebras 
$U_{\theta }\subseteq U_q^+ (\frak{ sl}_{n+1}),$ $U_{\theta ^{\prime }}\subseteq U_q^-(\frak{ sl}_{n+1})$ 
the tensor product 
\begin{equation}
 \hbox{\bf U}\, =U_{\theta ^{\prime }}\otimes _{{\bf k}[F]} {\bf k}[H]
\otimes _{{\bf k}[G]}U_{\theta }
\label{truint}
\end{equation}
is a right coideal subalgebra. In this hypothesis just one statement fails, 
the tensor product indeed is a right coideal but not always a subalgebra.

To describe conditions when (\ref{truint}) is a subalgebra
we display the element  $\Psi ^{\hbox{\bf s}}(k,m)$
schematically as a sequence of black and white points labeled by the numbers
$k-1,$ $k,$ $k+1, \ldots $ $m-1,$ $m,$ where the first point is always white, and
the last one is always black, while an intermediate point labeled by $i$ is black if and only if 
$i\in {\bf S}:$  
$$
 \stackrel{k-1}{\circ } \ \ \stackrel{k}{\circ } \ \ \stackrel{k+1}{\circ } 
\ \ \stackrel{k+2}{\bullet }\ \ \ \stackrel{k+3}{\circ }\ \cdots
\ \ \stackrel{m-2}{\bullet } \ \ \stackrel{m-1}{\circ }\ \ \stackrel{m}{\bullet }
$$
Consider two elements 
$\Psi ^{T_k}(k,\tilde{\theta }_k)$ and $\Psi ^{T_i^{\prime }}(i,\tilde{\theta }_i^{\prime }),$
where $T_k,$ $T_i^{\prime }$ are defined as above by $\theta $ and $\theta ^{\prime},$
respectively. Let us display these elements graphically
\begin{equation}
\begin{matrix}
\stackrel{k-1}{\circ } \ & \cdots \ & \stackrel{i-1}{\bullet } 
\ & \stackrel{i}{\bullet }\ \ & \stackrel{i+1}{\circ }\ & \cdots &
\ & \stackrel{\tilde{\theta }_k}{\bullet } \ & \ & \stackrel{\tilde{\theta }_i^{\prime }}{\cdot } \cr
\ \ & \  \ & \circ  
\ & \circ \ \ & \bullet \ & \cdots &
\ & \bullet  \ & \cdots  \ & \bullet
\end{matrix}
\label{gra1int}
\end{equation}
In Theorem \ref{osn5} we prove that (\ref{truint}) is a subalgebra if and only if for each pair
$(k, i),$ $1\leq k,i\leq n$ one of the following two options is fulfilled:

a) Representation (\ref{gra1int}) has no fragments of the form 
$$
\begin{matrix}
\stackrel{t}{\circ } \ & \cdots & \stackrel{l}{\bullet } \cr
\circ  
\ & \cdots  & \bullet 
\end{matrix}
$$

b) Representation (\ref{gra1int}) has the form 
$$
\begin{matrix}
\stackrel{k-1}{\circ } \ & \cdots & \circ & \cdots & \bullet & \cdots & \stackrel{m}{\bullet } \cr
\circ  
\ & \cdots  &  \bullet & \cdots & \circ & \cdots &  \bullet 
\end{matrix}
$$
where no one of the intermediate columns has points of the same color.

The obtained criterion allows use of the computer in order to find the total number $C_n$ of right coideal subalgebras which contain the coradical: 
$$
C_2=26; \ C_3=252; \ C_4=3,368; \ C_5=58,810; \ C_6=1,290,930; \ C_7=34,604,844.
$$

\noindent
{\bf Remark}. If a Hopf algebra $H$ has a Hopf algebra pairing 
$\langle , \rangle : M\times H\rightarrow {\bf k}$ with a Hopf algebra $M,$
then $M$ acts on $H$ via $m\rightharpoonup h=\sum h^{(1)}\langle m, h^{(2)}\rangle .$
Certainly, in this case each right coideal is $M$-invariant. Conversely, if the pairing is left faithful
(that is, $\langle M, h\rangle =0$ implies $h=0$)  then each $M$-invariant subspace is a right coideal.
For $H=U_q(\frak{sl}_{n+1})$ (or for $H=u_q(\frak{sl}_{n+1})$ if $q^t=1$) there exists a Hopf algebra pairing with $M=GL_q(n),$ see \cite{RTF, CM, Tow}.
Hence, alternatively, our main result provides a classification of $GL_q(n)$-invariant subalgebras that contain the coradical.

\smallskip

The computer part of this work has been done by the second author, 
while the proofs are due to the first one.

\section{Preliminaries}
\noindent
{\bf PBW-generators}.
Let $S$ be an algebra over a field {\bf k} and $K$ its subalgebra
with a fixed basis $\{ g_j\, |\,  j\in J\} .$ A linearly ordered subset $W\subseteq S$ is said to be
a {\it set of PBW-generators of $S$ over $K$} if there exists 
a function $h:W\rightarrow {\bf Z}^+\cup {\infty },$
called the {\it height function}, such that the set of all products
\begin{equation}
g_jw_1^{n_1}w_2^{n_2}\cdots w_k^{n_k}, 
\label{pbge}
\end{equation}
where $j\in J,\ \ w_1<w_2<\ldots <w_k\in W,\ \ n_i<h(w_i), 1\leq i\leq k$
is a basis of $S.$ The value $h(w)$ is referred to as the {\it height} of $w$ in $W.$
If $K={\bf k}$ is the ground field, then we shall call 
$W$  simply as a {\it set of PBW-generators of} $S.$ 

\noindent 
{\bf Character Hopf algebras}.
Recall that a Hopf algebra $H$ is referred to as a {\it character} Hopf 
algebra if the group $G$ of all group like elements is commutative
and $H$ is generated over {\bf k}$[G]$ by skew primitive 
semi-invariants $a_i,\ i\in I:$
\begin{equation}
\Delta (a_i)=a_i\otimes 1+g_i\otimes a_i,\ \ \ \
g^{-1}a_ig=\chi ^i(g)a_i, \ \ g, g_i\in G,
\label{AIc}
\end{equation}
where $\chi ^i,\, i\in I$ are characters of the group $G.$
By means of the Dedekind Lemma it is easy to see that 
every character Hopf algebra is graded by the monoid $G^*$ of characters
generated by $\chi ^i:$
\begin{equation}
H=\bigoplus_{\chi \in G^*} H^{\chi }, \ \ H^{\chi }=\{ a\in H\ |\ g^{-1}ag=\chi (g)a, \  g\in G\}.
\label{grad}
\end{equation}
Let us associate a ``quantum" variable $x_i$ to $a_i.$
For each word $u$
in $X=\{ x_i\, |\, i\in I\}$
we denote by $g_u$ or gr$(u)$
an element of $G$
that appears from $u$
by replacing each $x_i$ with $g_i.$
In the same way we denote by $\chi ^u$
a character that appears from $u$
by replacing each $x_i$ with $\chi ^i.$
We define a bilinear skew commutator on homogeneous 
linear combinations of words by the formula
\begin{equation}
[u,v]=uv-\chi ^u(g_v) vu,
\label{sqo}
\end{equation}
where sometimes for short we use the notation  $\chi ^u(g_v)=p_{uv}=p(u,v).$
These brackets satisfy the following Jacobi identities, see \cite[(8)]{Kh3}:
\begin{equation}
[[u, v],w]=[u,[v,w]]+p_{wv}^{-1}[[u,w],v]+(p_{vw}-p_{wv}^{-1})[u,w]\cdot v.
\label{jak1}
\end{equation}
\begin{equation}
[[u, v],w]=[u,[v,w]]-p_{vu}^{-1}[v,[u,w]]+(p_{vu}^{-1}-p_{uv})v\cdot [u,w].
\label{ja}
\end{equation}
In particular the following conditional  identities are valid
\begin{equation}
[[u, v],w]=[u,[v,w]],\hbox{ provided that } [u,w]=0.
\label{jak3}
\end{equation}
\begin{equation}
[u,[v,w]]=p_{uv}[v,[u,w]],\hbox{ provided that } [u,v]=0 \hbox{ and }p_{uv}p_{vu}=1.
\label{jak4}
\end{equation}
The brackets are related to the product by the following ad-identities
\begin{equation}
[u\cdot v,w]=p_{vw}[u,w]\cdot v+u\cdot [v,w], 
\label{br1f}
\end{equation}
\begin{equation}
[u,v\cdot w]=[u,v]\cdot w+p_{uv}v\cdot [u,w].
\label{br1}
\end{equation}
In particular, if $[u,w]=0,$ we have
\begin{equation}
[u\cdot v,w]=u\cdot [v,w].
\label{br2}
\end{equation}
The antisymmetry identity takes the form
\begin{equation}
[u,v]=-p_{uv}[v,u] \ \ \hbox{ provided that } \ \ p_{uv}p_{vu}=1.
\label{bri}
\end{equation}
The group  $G$ acts on the free algebra ${\bf k}\langle X\rangle $
by $ g^{-1}ug=\chi ^u(g)u,$ where $u$ is an arbitrary monomial 
in $X.$
The skew group algebra $G\langle X\rangle $
has the natural Hopf algebra structure 
$$
\Delta (x_i)=x_i\otimes 1+g_i\otimes x_i, 
\ \ \ i\in I, \ \ \Delta (g)=g\otimes g, \ g\in G.
$$
We fix a Hopf algebra homomorphism
\begin{equation}
\xi :G\langle X\rangle \rightarrow H, \ \ \xi (x_i)=a_i, \ \ \xi (g)=g, \ \ i\in I, \ \ g\in G.
\label{gom}
\end{equation}

\noindent 
{\bf Algebra ${{\mathfrak A}}_n$}. Suppose that the quantification parameters
$p_{ij}\stackrel{df}{=}p(x_i,x_j)=\chi ^i(g_j)$ satisfy
$p_{ij}p_{ji}=1$ provided $|i-j|>1.$
In this case all elements $[x_i,x_j]=x_ix_j-p_{ij}x_jx_i,$ $|i-j|>1$ are skew primitive.
Therefore the ideal of $G\langle X\rangle $ generated by these elements is a Hopf ideal.
We denote by ${{\mathfrak A}}_n$ the quotient character 
Hopf algebra, 
$$
{{\mathfrak A}}_n=G\langle X\, ||\, [x_i,x_j]=0,\  j-i>1 \rangle \stackrel{df}{=} 
G\langle X\rangle / {\rm Id}\, \langle [x_i,x_j],\ j-i>1\rangle .
$$

\begin{definition} \rm
The elements $u,v$ are said to be 
{\it separated} if there exists an index $j,$ $1\leq j\leq n,$
such that either $u\in {\bf k}\langle x_i\ |\ i<j\rangle ,$
$v\in {\bf k}\langle x_i\ |\ i>j\rangle $ or vise versa 
$u\in {\bf k}\langle x_i\ |\ i>j\rangle ,$
$v\in {\bf k}\langle x_i\ |\ i<j\rangle .$
\label{sep}
\end{definition}
In the algebra ${{\mathfrak A}}_n$ every two separated
homogeneous elements $u,v$  (skew)commute, $[u,v]=0,$ due to (\ref{br1f}), (\ref{br1}).

Let us consider a word $u(k,m)=x_kx_{k+1}\ldots x_m.$ There are a lot of options to 
rearrange the brackets in this word. For example 
$[x_k,[x_{k+1},\ldots [x_{m-1},x_m]\ldots]],$ or  $[[\ldots [x_k,x_{k+1}],\ldots x_{m-1}],x_m].$
\begin{lemma}
The value in ${{\mathfrak A}}_n$ of the bracketed continuous word $u(k,m)$
is independent of the alignment of brackets. In particular for the element 
\begin{equation}
u[k,m]\stackrel{df}{=}[x_k,[x_{k+1},\ldots [x_{m-1},x_m]\ldots]]=[x_k,u[k+1, m]]
\label{bcaw}
\end{equation}
we have the following equalities in 
${{\mathfrak A}}_n:$
\begin{equation}
u[k,m]=\hbox{\rm \Large [}u[k,s], u[s+1,m]\hbox{\rm \Large ]}, \ \ k\leq s<m.
\label{alin}
\end{equation}
\label{ind}
\end{lemma}
\begin{proof} We use induction on $m-k.$
If $[a]=[[v],[w]],$ $vw=u(k,m)$ then by the inductive supposition it suffices to consider
the case when $[v],$ $[w]$ are bracketed like $u[k,m]$ in (\ref{bcaw}); that is, $[v]=[x_k,[v_1]],$
$[w]=[x_m,[w_1]].$ 
If $v_1$ is empty then we 
come to $u[k,m].$ If $v_1$ is not empty, then $x_k$ and $w$ are separated, hence 
$[x_k,[w]]=0$ in ${{\mathfrak A}}_n.$
According to (\ref{jak3}) this implies 
{\Large [}$[x_k,[v_1]],[w]${\Large ]}$=$ {\Large [}$x_k,[[v_1],[w]]${\Large ]}, and the inductive supposition
applied to $[[v_1],[w]]$ again brings the value to that of $u[k,m].$ \end{proof}

\smallskip
\noindent 
{\bf PBW-basis of a character Hopf algebra}.
A {\it constitution} of a word $u$ in $G \cup X$
is a family of non-negative integers $\{ m_x, x\in X\} $
such that $u$ has $m_x$ occurrences of $x.$
Certainly almost all $m_x$ in the constitution are zero.
We fix an arbitrary complete order, $<,$ on the set $X.$

Let $\Gamma ^+$ be the free additive (commutative) monoid generated by $X.$
The monoid $\Gamma ^+$ \label{Gamma} is a completely ordered monoid with respect to 
the following order:
\begin{equation}
m_1x_{i_1}+m_2x_{i_2}+\ldots +m_kx_{i_k}>
m_1^{\prime }x_{i_1}+m_2^{\prime }x_{i_2}+\ldots +m_k^{\prime }x_{i_k}
\label{ord}
\end{equation}
if the first from the left nonzero number in
$(m_1-m_1^{\prime}, m_2-m_2^{\prime}, \ldots , m_k-m_k^{\prime})$
is positive, where $x_{i_1}>x_{i_2}>\ldots >x_{i_k}$ in $X.$
We associate a formal degree $D(u)=\sum _{x\in X}m_xx\in \Gamma ^+$
to a word $u$ in $G\cup X,$ where $\{ m_x\, | \, x\in X\}$ is the constitution of $u$
(in \cite[\S 2.1]{FG} the formal sum $D(u)$ is called the {\it weight} of $u$).
Respectively, if $f=\sum \alpha _i u_i\in G\langle X\rangle ,$ $0\neq \alpha _i\in {\bf k}$
then 
\begin{equation}
D(f)={\rm max}_i\{ D(u_i)\} .
\label{degr}
\end{equation}

On the set of all words in $X$ we fix the lexicographical order
with the priority from the left to the right,
where a proper beginning of a word is considered to 
be greater than the word itself.

A non-empty word $u$
is called a {\it standard} word (or {\it Lyndon} word, or 
{\it Lyndon-Shirshov} word) if $vw>wv$
for each decomposition $u=vw$ with non-empty $v,w.$
A {\it nonassociative} word is a word where brackets 
$[, ]$ somehow arranged to show how multiplication applies.
If $[u]$ denotes a nonassociative word then by $u$ we denote 
an associative word obtained from $[u]$ by removing the brackets
(of course, $[u]$ is not uniquely defined by $u$ in general, however Lemma \ref{ind}
says that the value of $[u]$ in ${{\mathfrak A}}_n$ is uniquely defined provided that $u=u(k,m)$).
The set of {\it standard nonassociative} words is the biggest set $SL$
that contains all variables $x_i$
and satisfies the following properties.

1)\ If $[u]=[[v][w]]\in SL$
then $[v],[w]\in SL,$
and $v>w$
are standard.

2)\  If $[u]=[\, [[v_1][v_2]]\, [w]\, ]\in SL$ then $v_2\leq w.$

\noindent
Every standard word has
only one alignment of brackets such that the appeared
nonassociative word is standard (Shirshov theorem \cite{pSh1}).
In order to find this alignment one may use the following
procedure: The factors $v, w$
of the nonassociative decomposition $[u]=[[v][w]]$
are the standard words such that $u=vw$
and  $v$ has the minimal length (\cite{pSh2}, see also \cite{Lot}).
\begin{definition}  \rm A {\it super-letter}
is a polynomial that equals a nonassociative standard word
where the brackets mean (\ref{sqo}).
A {\it super-word} is a word in super-letters. 
A $G$-{\it super-word} is a super-word multiplied from the left by a group-like element. 
\label{sup1}
\end{definition}

By Shirshov's theorem every standard word $u$
defines only one super-letter, in what follows we shall denote it by $[u].$
The order on the super-letters is defined in the natural way:
$[u]>[v]\iff u>v.$
\begin{definition} \rm
A super-letter $[u]$
is called {\it hard in }$H$
provided that its value in $H$
is not a linear combination
of values of super-words of the same degree (\ref{degr})
in smaller than $[u]$ super-letters, 
\underline{and $G$-super-words of smaller degrees}.
\label{tv1}
\end{definition}
\begin{definition} \rm
We say that a {\it height} of a hard in $H$ super-letter $[u]$
equals $h=h([u])$ if  $h$
is the smallest number such that: first, $p_{uu}$
is a primitive $t$-th root of 1 and either $h=t$
or $h=tl^r,$ where $l=$char({\bf k}); and  then the value in $H$
of $[u]^h$
is a linear combination of super-words of the same degree (\ref{degr})
in  less than $[u]$ super-letters,
\underline{and $G$-super-words of smaller degrees}.
If there exists no such number  then the height equals infinity.
\label{h1}
\end{definition}
Certainly, if the algebra $H$ is homogeneous in each $a_i$ then one may omit
the underlined parts of the definitions.
\begin{theorem} $(${\rm \cite[{Theorem 2}]{Kh3}}$).$
The values of all hard in $H$ super-letters with the above defined height function
form a set of PBW-generators for $H$ over {\bf k}$[G].$
\label{BW}
\end{theorem}
 
\noindent 
{\bf PBW-basis of a coideal subalgebra}. According to \cite[Theorem 1.1]{KhT, KhQ}
every right coideal subalgebra {\bf U} that contains all group-like elements has a PBW-basis
over {\bf k}$[G]$ which can be extended up to a PBW-basis of $H.$

The PBW-generators $T$ for {\bf U} 
can be obtained from the PBW-basis of $H$ given in Theorem \ref{BW}
in the following way. 

Suppose that for a given hard super-letter $[u]$ there exists  an element $c\in ${\bf U}
with the leading term $[u]^s$ in the PBW-decomposition given in Theorem \ref{BW}:    
\begin{equation}
c=[u]^s+\sum \alpha _iW_i+\ldots \in \hbox{\bf U},
\label{vad1}
\end{equation}
where $W_i$ are the basis super-words starting with less than $[u]$ super-letters, 
$D(W_i)=sD(u),$ and by the dots we denote a linear combination of $G$-super-words 
of $D$-degree less than $sD(u).$ We fix one of the elements with the minimal $s,$ and denote it by $c_u.$ Thus, for every hard in $H$ super-letter $[u]$ we have at most one element $c_u.$
We define the height function by means of the following lemma.

\begin{lemma}$\!\!\! ($\cite[{\rm Lemma 4.3}]{KhT, KhQ}$).$ 
In the representation $(\ref{vad1})$ of the chosen element $c_u$
either $s=1,$ or $p(u,u)$ is a primitive $t$-th root of $1$ and $s=t$ or 
$($in the case of positive characteristic$)$ $s=t({\rm char}\, {\bf k})^r.$
\label{nco1}
\end{lemma}
If the height of $[u]$ in $H$ is infinite, then the height of $c_u$ in {\bf U}
is defined to be infinite as well. If the height of $[u]$ in $H$ equals $t,$ then, due to 
the above lemma, $s=1$ (in the PBW-decomposition (\ref{vad1}) the exponent 
$s$ must be less than
the height of $[u]$). In this case the height of $c_u$ in {\bf U} is supposed to be $t$ as well.
If the characteristic $l$ is positive, and the height of $[u]$ in $H$ equals
$tl^r,$ then we define the height of $c_u$ in {\bf U} to be equal to $tl^r/s$
(thus, in characteristic zero the height of $c_u$ in {\bf U} always
equals the height of $[u]$ in $H$).

\begin{proposition}
The set of all chosen $c_u$ with the above defined height function
forms a set of PBW-generators for {\bf U} over {\bf k}$[G].$  
\label{pro}
\end{proposition}
\begin{proof} See, \cite[Proposition 4.4]{KhT, KhQ}. \end{proof}

We note that there is an essential freedom in construction of the PBW-generators 
for a right coideal subalgebra.  In particular the PBW-basis is not uniquely defined in the above process. Nevertheless the set of leading terms of the PBW-generators indeed is uniquely defined.

\begin{definition} \rm A hard super-letter $[u]$ is called {\bf U}-{\it effective} if there exists
$c\in \, ${\bf U} of the form (\ref{vad1}). The degree $sD(u)\in \Gamma ^+ $ of $c$ with minimal $s$
is said to be an {\bf U}-{\it root}. An {\bf U}-root $\gamma \in \Gamma ^+ $
 is called a {\it simple} {\bf U}-{\it root} if it is not a sum of two or more other {\bf U}-roots.
\label{root}
\end{definition}
Thus, the set of {\bf U}-effective super-letters, the set of {\bf U}-roots, and the set of 
simple  {\bf U}-roots are invariants of any right coideal subalgebra {\bf U}.

{\bf Remark}. There is already a fundamental for Lie theory notion of roots associated
to semisimple Lie algebras. Certainly, the set of PBW-generators for 
the universal enveloping algebra $U({\frak g})$ 
coincides with a basis of the Lie algebra ${\frak g}.$ If we apply our definition to  $U({\frak g})$ then $U({\frak g})$-roots are the formal degrees of basis elements related to a fixed set of generators
$x_i, i\in I.$ At the same time the formal degrees of basis elements for the Borel subalgebra 
are in one-to-one correspondence with positive roots: to each root 
$\alpha _{i_1}+\alpha _{i_2}+\cdots +\alpha _{i_k}$
corresponds a basis element $[\ldots [x_{i_1},x_{i_2}],\ldots ,x_{i_k}],$ 
see \cite[Chapter IV, \S 3, Statement XVII]{Jac}. Therefore our definition of a root
is a natural generalization of the classical notion. Probably the analogy would be 
more clear if in our definition of the formal degree 
we  will replace the symbols $x_i$  with the characters $\chi ^i$
and identify the generators $g_i$ of the group $G$ with (exponents of the) basis elements $h_i$
of the Cartan subalgebra since the classical roots are elements of the dual space
$(\sum_i {\bf k}h_i)^*.$ Lemma \ref{odn} below shows that this replacement is admissible.
We belive that by this very reason in \cite{FG} the formal degree is referred to as {\it weight},
the notion already well defined in the Lie theory. 

\smallskip
\noindent 
{\bf Differential calculi}. 
The free algebra ${\bf k}\langle X\rangle $ has a coordinate differential calculus
\begin{equation}
\partial_i(x_j)=\delta _i^j,\ \ \partial _i (uv)=\partial _i(u)\cdot v+\chi ^u(g_i)u\cdot \partial _i(v).
\label{defdif}
\end{equation}
The partial derivatives connect the calculus with the coproduct on $G\langle X\rangle$ via
\begin{equation}
\Delta (u)\equiv u\otimes 1+\sum _ig_i\partial_i(u)\otimes x_i\ \ \ 
(\hbox{mod }G\langle X\rangle \otimes \hbox{\bf k}\langle X\rangle ^{(2)}),
\label{calc}
\end{equation}
where ${\bf k}\langle X\rangle ^{(2)}$ is the set (an ideal)
of noncommutative polynomials without free and linear terms.
Symetrically the equation
\begin{equation}
\Delta (u)\equiv g_u\otimes u+\sum _ig_ug_i^{-1}x_i\otimes\partial _i^*(u)\ \ \ 
(\hbox{mod }G\langle X\rangle ^{(2)}\otimes \hbox{\bf k}\langle X\rangle )
\label{calcdu}
\end{equation}
defines a dual differential calculus on ${\bf k}\langle X\rangle $
where the partial derivatives satisfy
\begin{equation}
\partial _j^*( x_i)=\delta _i^j,\ \ \partial _i^*(uv)=
\chi ^{i}(g_v) \partial _i^*(u)\cdot v+u\cdot \partial _i^*(v).
\label{difdu}
\end{equation}
Here $G\langle X\rangle ^{(2)}{\bf k}$ is an ideal of $G\langle X\rangle $
of elements without free and linear terms.
If the kernel of $\xi $ defined in (\ref{gom}) is contained in $G\langle X\rangle ^{(2)}$ then formule
(\ref{calc}), (\ref{calcdu}) with the $a$'s in place of the $x$'s define  coordinate differential calculi
on the subalgebra $A$ of $H$ generated by the $a$'s. In this case restriction of $\xi $
on $\hbox{\bf k}\langle X\rangle$ is a differential homomorphism, while
 (\ref{calc}) and (\ref{calcdu}) imply that each skew-primitive element $u$
from $A^{(2)}=\xi ({\bf k}\langle X\rangle ^{(2)})$ is a constant with respect to both calculi,
$\partial _i(u)=\partial _i^*(u)=0,$ $1\leq i\leq n.$
More details one can find in \cite{MS, Kh1, Kh5}.

\smallskip
\noindent 
{\bf Shuffle representation}. Let Ker$\, \xi \subseteq G\langle X\rangle ^{(2)}.$ 
In this case there exists a Hopf algebra projection $\pi : H\rightarrow {\bf k}[G],$
$a_i\rightarrow 0, $ $g_i\rightarrow g_i .$ Hence by the Radford theorem
\cite{Rad}  we have a decomposition in a biproduct,
$H=A\# {\bf k}[G],$ by means of the isomorphism
$u\rightarrow \vartheta (u^{(1)})\# \pi(u^{(2)})$ with $\vartheta (u)=\sum _{(u)}u^{(1)}\pi (S(u^{(2)})),$
see details in  \cite[\S 1.5, \S 1.7]{AS}.

If Ker$\, \xi $ is the biggest Hopf ideal
in $G\langle X\rangle ^{(2)},$ or, equivalently, if $H$ is a Hopf algebra 
of type one in the sense of Nichols \cite{Nic}, or, equivalently, 
if $A$ is a {\it quantum symmetric algebra} (a Nichols algebra \cite[\S 1.3, Section 2]{AS}),
then $A$ has a shuffle representation as follows.

The algebra $A$ has a structure of a {\it braided Hopf algebra}, \cite{Tak1}, 
with a braiding $\tau (u\otimes v)=p(v,u)^{-1}v\otimes u.$
The braided coproduct 
$\Delta ^b$ is connected with the coproduct on $H$ in the following way,
see \cite[p. 93, (3.18)]{Kh5},
\begin{equation}
\Delta ^b(u)=\sum _{(u)}u^{(1)}\hbox{gr}(u^{(2)})^{-1}\underline{\otimes} u^{(2)}.
\label{copro}
\end{equation}
At the same time the tensor space $T(V),$ $V=\sum_i {\bf k}x_i$ also has a structure of a braided Hopf algebra.
This is the {\it quantum shuffle algebra} $Sh_{\tau }(V)$ with the coproduct 
\begin{equation}
\Delta ^b(u)=\sum _{i=0}^m(z_1\ldots z_i)\underline{\otimes} (z_{i+1}\ldots z_m),
\label{bcopro}
\end{equation}
where $z_i\in X,$ and $u=(z_1z_2\ldots z_{m-1}z_m)$ is the tensor
$z_1\otimes z_2\otimes \ldots \otimes z_{m-1}\otimes z_m$ 
considered as an element of $Sh_{\tau }(V).$
The map $a_i\rightarrow (x_i)$ defines an embedding of the braided Hopf algebra
$A$ into the the braided Hopf algebra $Sh_{\tau }(V).$
More details can be find in \cite{Nic, Wor, Scha, Ros, AG, FG1,Tak1, Flo, Kh5, Kh1}.

\smallskip
\noindent 
{\bf Differential subalgebras}. 
If {\bf U} is a right coideal subalgebra of $H$ and {\bf k}$[G]\subseteq \,${\bf U}, then 
$\vartheta (\hbox{\bf U})\subseteq \, ${\bf U}, hence 
{\bf U}$=U_A\# {\bf k}[G]$
with $U_A=\vartheta (\hbox{\bf U})=\hbox{\bf U}\cap A.$ Formula
(\ref{calc}) implies that $U_A$ is a differential subalgebra of $A,$
that certainly satisfies $gU_Ag^{-1}\subseteq U_A,$ $g\in G.$
The converse statement is valid if Ker$\, \xi $ is the biggest Hopf ideal.
\begin{lemma} Suppose that Ker$\, \xi $ is the biggest Hopf ideal
 in $G\langle X\rangle ^{(2)}.$
If $U$ is a differential subspace of $A={\bf k}\langle a_i\rangle $ $=\vartheta (H),$
and $gUg^{-1}\subseteq U,$ $g\in G,$
then $U\# {\bf k}[G]$ is a right coideal of $H.$
\label{qsim}
\end{lemma}
\begin{proof}
The braided coproduct (\ref{copro}) also defines a differential calculus 
\begin{equation}
\Delta ^b(u)\equiv u\underline{\otimes} 1+
\sum _i\frac{\partial ^bu}{\partial x_i}\underline{\otimes} x_i\ \ \ 
(\hbox{mod }A \underline{\otimes} A^{(2)}).
\label{calc1}
\end{equation}
In  \cite[Theorem 4.8]{Kh1} this calculus is denoted by $d^*.$
Formulae (\ref{calc}), (\ref{copro}), and (\ref{calc1}) imply
$\partial ^bu/\partial x_i=g_i\partial _i(u)g_i^{-1}.$

Since $A$ has a representation as a subalgebra of the 
quantum shuffle algebra $Sh_{\tau }(V),$  
by \cite[Theorem 5.1]{Kh1} applied to the calculus $d^*$
the restriction $\Omega =\xi |_{\hbox{\bf k}\langle X\rangle}$ has the following
differential form
\begin{equation}
\Omega (u)=\sum _{i_1,i_2,\cdots ,i_n}\frac{(\partial^b) ^nu}{\partial x_{i_1}\partial  x_{i_2}\cdots 
\partial x_{i_n}}(x_{i_n}x_{i_{n-1}}\cdots x_{i_1}),\ \ \ u\in V^{\otimes n},
\label{oderc3}
\end{equation} 
where as above $(x_{i_n}x_{i_{n-1}}\cdots x_{i_1})$ is the tensor
$x_{i_n}{\otimes } x_{i_{n-1}}{\otimes } \cdots {\otimes } x_{i_1}$
considered as an element  of $Sh_{\tau }(V).$ By means of (\ref{bcopro}) we have 
$$
\Delta ^b(\Omega (u))=\sum _{i_1,i_2,\cdots ,i_n}
\frac{(\partial ^b)^nu}{\partial x_{i_1}\partial x_{i_2}\cdots 
\partial x_{i_n}} \sum _{k=1}^{n+1}(x_{i_n}\cdots x_{i_k})
\underline{\otimes }(x_{i_{k-1}}\cdots x_{i_1})
$$

$$
=\sum _{k=1}^{n+1}\sum _{i_1,i_2,\cdots ,i_{k-1}}\left(\sum _{i_k,i_{k+1},\cdots ,i_n}
\frac{(\partial ^b)^{n-k+1}\left[ \frac{(\partial ^b)^{k-1}u}{\partial x_{i_1} \cdots 
\partial x_{i_{k-1}}}\right]}{\partial x_{i_k} \cdots 
\partial x_{i_n}} (x_{i_n}\cdots x_{i_k})\right)
\underline {\otimes }
 (x_{i_{k-1}}\cdots x_{i_1})
$$
\begin{equation}
=\sum _{k=1}^{n+1}\sum _{i_1,i_2,\cdots ,i_{k-1}}
\Omega \left( \frac{(\partial ^b)^{k-1} u}{\partial x_{i_1} \cdots 
\partial x_{i_{k-1}}}\right) 
\underline {\otimes }
(x_{i_{k-1}}\cdots x_{i_1}).
\label{oderc1}
\end{equation}
Since $\Omega $ is a $d^*$-differential map, we have got 
\begin{equation}
\Delta ^b(w)=\sum _{k=1}^{n+1}
\sum _{i_1,i_2,\cdots ,i_{k-1}}
\frac{(\partial ^b)^{k-1} w}{\partial x_{i_1}\cdots \partial x_{i_k}}
\underline {\otimes }
(x_{i_{k-1}}\cdots x_{i_1}), \ \ w=\Omega (u).
\label{oderc4}
\end{equation} 
This formula implies that each differential subspace $W\subseteq A$ with respect to $d^*$
 is a right coideal with respect to $\Delta ^b.$
Indeed, $\Delta ^b(W)\subseteq (A\underline{\otimes } A)\cap 
(W\underline{\otimes } Sh_{\tau})$ $=W\underline {\otimes } A.$
Since
$\partial ^bu/\partial x_i=g_i\partial_i(u)g_i^{-1},$
the space $U$ given in the lemma is a right coideal
with respect to the coproduct $\Delta ^b.$ Now (\ref{copro})
shows that $U{\bf k}[G]$ is a right coideal of $H.$  
The lemma is proved.\end{proof}

\section{Multiparameter quantification of Kac-Moody algebras\\ as character Hopf algebras}

\noindent
{\bf Quantification of Borel subalgebras}.
Let $C=||a_{ij}||$ be a symmetrizable by $D={\rm diag }(d_1, \ldots d_n)$ generalized Cartan matrix, $d_ia_{ij}=d_ja_{ji}.$
Denote by $\mathfrak g$  a Kac-Moody algebra defined by $C,$ see \cite{Kac}.
Suppose that the quantification parameters $p_{ij}=p(x_i,x_j)=\chi ^i(g_j)$ are related by
\begin{equation}
p_{ii}=q^{d_i}, \ \ p_{ij}p_{ji}=q^{d_ia_{ij}},\ \ \ 1\leq i,j\leq n. 
\label{KM1}
\end{equation}
In this case the multiparameter quantization $U^+_q ({\mathfrak g})$ of the 
Borel subalgebra ${\mathfrak g}^+$
is a character Hopf algebra defined by Serre relations 
with the skew brackets in place of the Lie operation:
\begin{equation}
[\ldots [[x_i,\underbrace{x_j],x_j], \ldots ,x_j]}_{1-a_{ji} \hbox{ times}}=0, \ \ 1\leq i\neq j\leq n.
\label{KM2}
\end{equation}
By \cite[Theorem 6.1]{Khar} the left hand sides of these relations are skew-primitive elements 
in $G\langle X\rangle .$ Therefore the ideal generated by  these elements is a Hopf ideal,
while $U^+_q ({\mathfrak g})$ indeed has a natural character Hopf algebra structure.

\begin{lemma} If $C$ is a Cartan matrix of finite type $($in particular the symmetric matrix $DC$ is positively defined $),$ and $q$ is not a root of $1$ then the grading of $U^+_q ({\mathfrak g})$  by 
characters $(\ref{grad})$ coincides with the grading by $\Gamma ^+.$ 
\label{odn}
\end{lemma}
\begin{proof} Since every homogeneous in $\Gamma ^+$ element is homogeneous with respect to 
(\ref{grad}), it suffices to show that the characters $\chi ^i=\chi ^{x_i},$ $1\leq i\leq n$
generate a free Abelian group. Suppose in contrary that
\begin{equation}
\chi_1 \stackrel{df}{=}(\chi ^1)^{k_1}\cdots (\chi ^n)^{k_n}=(\chi ^1)^{m_1}\cdots (\chi ^n)^{m_n}
\stackrel{df}{=} \chi_2,
\label{fch}
\end{equation}
where $k_i, m_i\geq 0, \ k_im_i=0,$ $1\leq i\leq n,$ and one of the $k_i$'s is nonzero.
Let $g=g_1^{k_1}\cdots g_n^{k_n},$ $h=g_1^{m_1}\cdots g_n^{m_n}.$
By means of (\ref{KM1}) we have
$$
\chi _1(g)=\prod _{1\leq i,j\leq n} p_{ij}^{k_jk_i}=  \prod _{i<j}(p_{ij}p_{ji})^{k_ik_j}
\cdot \prod _ip_{ii}^{k_i^2}=q^N,
$$
where due to (\ref{KM1}) for $\vec{k}=(k_1, \ldots , k_n)$ we have 
$$
N=\sum_{i<j}d_ia_{ij}k_ik_j+\sum_{i}d_ik_i^2
$$
$$
=\frac{1}{2}
(\sum_{i<j}d_ia_{ij}k_ik_j+ \sum_{i<j}d_ja_{ji}k_jk_i+\sum_{i}a_{ii}d_ik_i^2)
=(\vec{k}DC,\vec{k})>0.
$$
In the same way $\chi _2(h)=q^M,$ $M\geq 0.$
Relations (\ref{KM1}) imply 
$$
\chi _2(g)\chi _1(h)=\prod _{1\leq i,j\leq n}p_{ij}^{m_ik_j}\cdot \prod _{1\leq i,j\leq n}p_{ij}^{m_jk_i}
=\prod _{1\leq i, j\leq n}(p_{ij}p_{ji})^{m_ik_j}=q^{L},
$$
with $L=\sum _{i,j}d_ia_{ij}m_ik_j \leq 0$ since in the Catran matrix $a_{ij}\leq 0$ for $i\neq j,$
while $k_im_i=0.$

We have $\chi _1(g)=\chi _2(g),$ and $\chi_1(h)=\chi _2(h).$ Therefore 
$q^{M+N}$ $=\chi _1(g)\chi _2(h)$ $=\chi _2(g)\chi _1(h)$ $=q^{L}.$ A contradiction.
\end{proof}

\smallskip
\noindent
{\bf Remark}. Of course, if the characters $\chi_i,$ $1\leq i\leq n$ generate a 
free Abelian group then $g_i,$ $1\leq i\leq n$ generate a free Abelian group as well.
In particular relations (\ref{KM1}) imply that $G$ is a free Abelian group with the free generators
$g_i,$  $1\leq i\leq n$ provided that $q$ is not a root of 1 and $C$ is of finite type.
\begin{corollary} 
If $q$ is not a root of 1 and $C$ is of fifnite type, 
then every subalgebra $U$ of $U_q^+({\frak g})$ containing $G$
is homogeneous with respect to each of the variables $x_i$.
\label{odn1}
\end{corollary}
\begin{proof} By the above lemma it suffices to note that $U$ is homogeneous 
with respect to the grading by characters (\ref{grad}). If $c=\sum_ic_i\in U$
with $c_i\in H^{\chi _i}$ and different $\chi_i\in G^*,$ then 
\begin{equation}
g^{-1}cg=\sum _i \chi_i(g)c_i\in U, \ \ g\in G.
\label{eqw}
\end{equation}
According to the Dedekind Lemma there exist
elements $h_i\in G,$ such that the matrix $M=||\chi_i(h_j)||$ is invertible.
Hence we may solve the system of equations (\ref{eqw}) considering $c_i$
as variables. In particular $c_i\in U.$
\end{proof}

\smallskip
If the multiplicative order $t$ of $q$ is finite, 
then we define $u^+_q ({\mathfrak g})$ as $G\langle X\rangle/{\bf \Lambda },$
where ${\bf \Lambda }$ is the biggest Hopf ideal in $G\langle X\rangle ^{(2)}.$
This is a ${\Gamma }^+$-homogeneous ideal, see \cite[Lemma 2.2]{KA}. 
Certainly ${\bf \Lambda }$ contains all skew-primitive elements of
$G\langle X\rangle ^{(2)}$ (each of them generates a Hopf ideal). Hence
by  \cite[Theorem 6.1]{Khar} relations (\ref{KM2}) are still valid in $u^+_q ({\mathfrak g}).$

\smallskip
\noindent 
{\bf Quantification of Kac-Moody algebras}.
Consider a new set of variables $X^-=$ $\{ x^-_1, x^-_2,\ldots, x^-_n\} .$ Suppose that an Abelian group $F$ generated by the elements $f_1, f_2, \ldots ,f_n$ acts on the linear space spanned by 
$X^-$ so that $(x_i^-)^{f_j}=p_{ji}^{-1}x_i^-,$ where $p_{ij}$ are the same parameters, 
see (\ref{KM1}), that define $U_q^+(\frak{g}).$ Relations (\ref{KM1}) are invariant
under the substitutions $p_{ij}\leftarrow p_{ji}^{-1},$ $q\leftarrow q^{-1}.$ This allows us to define the
character Hopf algebra $U_q^-(\frak{g})$ as $U_{q^{-1}}^+(\frak{g})$  with the characters
$\chi ^i_-,$ $1\leq i\leq n$ such that $\chi ^i_-(f_j)=p_{ji}^{-1}.$

We may extend the characters $\chi ^i $ on $G\times F$ in the following way
\begin{equation}
\chi ^i(f_j)\stackrel{df}{=}p_{ji}=\chi ^j(g_i).
\label{shar1}
\end{equation}
Indeed, if $\prod_k f_k^{m_k}=1$ in $F,$ then application to $x_i^-$
implies $\prod_k p_{ki}^{-m_k}=1,$ hence $\chi ^i(\prod _k f_k^{m_k})=\prod p_{ki}^{m_k}$
equals 1 as well. In the same way we may extend the characters $\chi ^i_-$ on $G\times F$
so that 
\begin{equation}
\chi ^i_-=(\chi ^i)^{-1} \ \ \hbox{as characters of } G\times F.
\label{shar2}
\end{equation}

In what follows we denote by $H$ a quotient group $(G\times F)/N,$ where 
$N$ is an arbitrary subgroup with $\chi ^{i}(N)=1,$ $1\leq i\leq n.$ For example, if the quantification parameters satisfy additional symmetry conditions $p_{ij}=p_{ji},$  $1\leq i,j\leq n,$
as this is a case for the original Drinfeld-Jimbo and Lusztig quantifications, then 
$\chi ^i(g_k^{-1}f_k)=p_{ik}^{-1}p_{ki}=1,$ and we may take $N$ to be the subgroup generated by 
$g_k^{-1}f_k,$  $1\leq k\leq n. $ In this particular case the groups $H,$ $G,$ $F$ may be identified. 

In the general case without loss of generality we may suppose that $G,F\subseteq H.$
Certainly $\chi ^i, 1\leq i\leq n$ are characters of $H$ and $H$ still 
acts on the  space spanned by $X\cup X^-$ by means of these characters and their inverses. 
 Consider the skew  group algebra $H\langle X\cup X^-\rangle $ as a character Hopf algebra:
\begin{equation}
\Delta (x_i)=x_i\otimes 1+g_i\otimes x_i,\ \ \ \Delta (x_i^-)=x_i^-\otimes 1+f_i\otimes x_i^-,
\label{AIcm}
\end{equation}
\begin{equation}
g^{-1}x_ig=\chi ^i(g)\cdot x_i, \ \ g^{-1}x_i^-g=(\chi ^i)^{-1}(g)\cdot x_i^-, \ \ g\in H.
\label{AIcm2}
\end{equation}
We define  the algebra $U_q(\frak{g})$ as a quotient
of $H\langle X\cup X^-\rangle $ by the following relations:
\begin{equation}
[\ldots [[x_i,\underbrace{x_j],x_j], \ldots ,x_j]}_{1-a_{ji} \hbox{ times}}=0, \ \ 1\leq i\neq j\leq n;
\label{rela1}
\end{equation}
\begin{equation}
[\ldots [[x_i^-,\underbrace{x_j^-],x_j^-], \ldots ,x_j^-]}_{1-a_{ji} \hbox{ times}}=0, \ \ 1\leq i\neq j\leq n;
\label{rela2}
\end{equation}
\begin{equation}
[x_i, x_j^-]=\delta_i^j(1-g_if_i), \ \ \ \ \ 1\leq i,j\leq n
\label{rela3}
\end{equation}
where the brackets are defined on $H\langle X\cup X^-\rangle $ by the structure of character Hopf algebra, see (\ref{sqo}). Since due to (\ref{KM1}) and \cite[Theorem 6.1]{Khar} all polynomials in the above relations are skew primitive in  $H\langle X\cup X^-\rangle ,$ they define  a Hopf ideal of 
$H\langle X\cup X^-\rangle $; that is, the natural homomorphism 
\begin{equation}
H\langle X\cup X^-\rangle \rightarrow U_q(\frak{g})
\label{gom2}
\end{equation}
defines a Hopf algebra structure on $U_q(\frak{g}).$

If $q$ has a finite multiplicative order then  $u_q(\frak{g})$ is defined by relations (\ref{rela3})
and $u=0,$ $u\in {\bf \Lambda },$ $u^-=0,$ $u^-\in {\bf \Lambda }^-,$ where ${\bf \Lambda },$
${\bf \Lambda }^-$ are the biggest Hopf ideals respectively in $G\langle X\rangle ^{(2)}$
and $F\langle X^-\rangle ^{(2)}.$

Both algebras $U_q(\mathfrak{g}),$ and $u_q(\mathfrak{g})$ have a grading by the additive 
group $\Gamma $ generated by $\Gamma ^+,$
see p.\pageref{Gamma}, provided that we put $D(x_i^-)=-D(x_i)=-x_i,$ $D(H)=0$ 
since in this way relations (\ref{rela3}) became homogeneous. 

\begin{corollary} 
If $q$ is not a root of 1 and the Cartan matrix $C=||a_{ij}||$ is of finite type then every subalgebra 
$U$ of $U_q(\frak{g})$ containing $H$ is $\Gamma $-homogeneous.
\label{odn11}
\end{corollary}
\begin{proof}
By Lemma \ref{odn} 
and definition (\ref{shar2}) grading by $\Gamma $ coincides with the grading
(\ref{grad}) by the group of characters (freely) generated by $\chi ^i,$ $1\leq i\leq n.$
Hence every subspace invariant under the conjugations by $H$ is 
$\Gamma $-homogeneous. 
\end{proof}

\smallskip
The defined quantification reduces to known ones under a suitable choice of 
$x_i, x_i^-$ depending up the particular definition of  $U_q (\frak{g}).$ For example 
for classical case of one parameter quantification we have $G=F=H,$ and in the 
notations of \cite{Luz2} we may identify
$$
x_i= E_i,\ g_i= K_i,\ x_i^-= 
F_iK_i(v^{-d_i}-v^{d_i})^{-1}, \ p_{ij}= v^{-d_ia_{ij}}, 
$$
while in the notations of \cite{Luz1, Mul} we may take
$$
x_i= E_i,\ g_i=\tilde{K}_i, \ x_i^-= F_i\tilde{K}_i(v_i^{-1}-v_i)^{-1}, \
p_{i\mu }= v^{-\langle \mu ,i^{\prime }\rangle }.
$$
For two-parameter quantizations, say in the notations of \cite{BGH}, we may put
$$
x_i\leftarrow e_i, \ g_i\leftarrow \omega _i, \ x_i^-\leftarrow f_i(\omega _i^{\prime })^{-1}(r_i-s_i)^{-1},\ 
f_i\leftarrow (\omega _i^{\prime })^{-1},
$$
and find values of parameters $p_{ij}$ by means of \cite[(B2), (C2), (D2)]{BGH}.
For the multiparameter case of Reshetikhin or DeConcini-Kac-Procesi 
in the notations of \cite{CV}, we may take 
$$
x_i\leftarrow E_iL_{\beta _i}, \ g_i\leftarrow L_{\beta _i-\alpha _i+\gamma _i}, \ x_i^-\leftarrow F_iL_{\alpha _i+\beta _i}^{-1}(q_i-q_i^{-1})^{-1},\ f_i\leftarrow L_{\gamma _i+\alpha _i+\beta _i}^{-1}.
$$

\smallskip
\noindent
{\bf Triangular decomposition}.
One may prove that the subalgebra of 
$U_q (\frak{g})$ generated by $G$ and values of $x_i,$ $1\leq i\leq n$ is isomorphic to
 $U_q^+ (\frak{g})$ 
while the subalgebra generated by $F$ and values of $x_i^-,$ 
$1\leq i\leq n$ is isomorphic to $U_q^- (\frak{g}).$ 
Moreover, one has the following so called ``triangular decomposition'' for both algebras:
\begin{equation}
U_q(\frak{g})= U_q^-(\frak{g})\otimes _{{\bf k}[F]} {\bf k}[H]
\otimes _{{\bf k}[G]}U_q^+(\frak{g}),
\label{tr}
\end{equation}
\begin{equation}
u_q(\frak{g})= u_q^-(\frak{g})\otimes _{{\bf k}[F]} {\bf k}[H]
\otimes _{{\bf k}[G]}u_q^+(\frak{g}).
\label{tr1}
\end{equation}
Actually this is not so evident (see \cite{Luz1, Luz2} for standard one parameter version, 
\cite{CV} for the multiparameter version with Cartan matrix of finite type,
 \cite{BGH} for two-parameter version with particular Cartan matrices only). 
We shall provide here a relatively short proof in the general setting that uses
a lemma on tensor decomposition
for character Hopf algebras, \cite[Lemma 6.2]{Kh4}, and (in case $q^t=1$) the 
Heyneman--Radford theorem.
\begin{proposition}
Let $J\subseteq G\langle X\rangle ^{(2)} ,$ $J^-\subseteq F\langle X^-\rangle ^{(2)}$
be constitution homogeneous Hopf ideals of 
$G\langle X\rangle $ and $F\langle X^-\rangle $ respectively.
Denote by ${\mathfrak A}$ the algebra generated over $H$ by $X\cup X^-$ and defined by the relations $(\ref{rela3})$ and
$u_s=0,$ $s\in S,$ $u^-_t=0,$ $t\in T,$ where $\{ u_s,$ $s\in S\} $ 
$($respectively $\{ u_t^-,$ $t\in T\} )$ 
is a set of homogeneous generators of the ideal $J$ $($respectively $J^-).$  We have
\begin{equation}
{\mathfrak A}= (F\langle X^-\rangle/J^-)\otimes _{{\bf k}[F]} {\bf k}[H]
\otimes _{{\bf k}[G]}(G\langle X\rangle/J).
\label{tria}
\end{equation}
\label{rtri}
\end{proposition}
\begin{proof}
We note first that the algebra ${\mathfrak A}_1$ generated over $H$ by $X$ and defined by the relations
$u_s=0,$ $s\in S$ has the form ${\bf k}[H]\otimes _{{\bf k}[G]}(G\langle X\rangle/J),$
while the algebra ${\mathfrak A}_2$ generated over $H$ by $X^-$ and defined by the relations
$u^-_t=0,$ $t\in T$ has the form $(F\langle X^-\rangle/J^-)\otimes _{{\bf k}[F]}{\bf k}[H].$
Hence it suffices to show that ${\mathfrak A}={\mathfrak A}_2\otimes _{{\bf k}[H]}{\mathfrak A}_1.$ 

Denote by $D_i,$ $D_i^- ,$ $1\leq i\leq n$ the linear maps 
$$D_i:{\bf k}\langle X^-\rangle \rightarrow H\langle X^-\rangle ,\ \ \ 
D_i^-:{\bf k}\langle X\rangle \rightarrow H\langle X\rangle $$
that satisfy the initial conditions 
\begin{equation}
D_i(x_j^-)=D_j^-(x_i)=\delta_i^j(1-g_if_i),
\label{inc}
\end{equation}
and the skew differential Leibniz rules
\begin{equation}
D_i(v^-\cdot w^-)=D_i(v^-)\cdot w^-+p(x_i, v^-)v\cdot D_i(w^-),\ \ \ v^-, w^-\in {\bf k}\langle X^-\rangle ;
\end{equation}
\begin{equation}
D_i^-(u\cdot v)=p(v,x_i^-)D_i^-(u)\cdot v+u\cdot D_i^-(v), \ \ \ u, v\in {\bf k}\langle X\rangle .
\label{sqdf2}
\end{equation}
Lemma 6.2 \cite{Kh4} (under the substitutions 
$k\leftarrow n,$ $n\leftarrow 2n,$  $x_{n+i}\leftarrow x_i^-,$ $G\leftarrow H,$ 
$H\leftarrow {\mathfrak A})$ 
gives the required decomposition provided that there exist homogeneous defining
relations $\{ \varphi _s=0,$ $s\in S\} ,$ and $\{ \psi _t=0,$ $t\in T\} $
 for ${\mathfrak A}_1$ and ${\mathfrak A}_2$ respectively,  such that 
\begin{equation}
D_i(\psi_t)\in H\cdot J^-,\ \  D_i^-(\varphi _s)\in H\cdot J, \ \ \ 1\leq i\leq n,\ s\in S,\ t\in T.
\label{inw}
\end{equation}
Consider the linear maps
\begin{equation}
\tilde{D}_i^-: u\rightarrow \partial_ i^*(u)
-p_{ii}^{-1}p(u,x_i)\partial_i(u)g_if_i, \ \ \ u\in{\bf k}\langle X\rangle ,
\label{sqi1}
\end{equation}
where the partial derivatives are defined in (\ref{calc}) and (\ref{calcdu}). We have 
$\tilde{D}_i^-(x_j)=\delta _i^j(1-g_if_i),$ while relations (\ref{defdif}) and (\ref{difdu}) imply
the differential Leibniz rule 
$\tilde{D}_i^-(u\cdot v)=p(x_i,v)\tilde{D}_i^-(u)\cdot v+u\cdot \tilde{D}_i^-(v).$ Since  
according to (\ref{shar1}) we have $p(x_i,v)=p(v, x_i^-),$ the Leibniz rule and initial values for 
$D_i^-$ coincides with that for $\tilde{D}_i^-.$ Hence $D_i^-=\tilde{D}_i^-.$ In perfect analogy we have
\begin{equation}
D_i(u^-)=\partial_{-i}^*(u^-)p(x_i,u^-)p_{ii}^{-1}-g_if_i \partial_{-i}(u^-),
 \ \ \ u^-\in{\bf k}\langle X^-\rangle ,
\label{sqi2}
\end{equation}
where $\partial _{-i},$ $\partial _{-i}^* $ are left and right partial derivatives on
 ${\bf k}\langle X^-\rangle $ with respect to $x_i^-.$

Now if  $u_s,$ $s\in S$ and $u_t^-,$ $t\in T$ are skew primitive elements
(as this is the case for ${\mathfrak A}=U_q(\frak{g})$) then they are constants 
for $\partial _i,$ $\partial _i^*,$ and $\partial _{-i},$ $\partial _{-i}^*,$
respectively.
Hence (\ref{sqi1}), (\ref{sqi2}) imply $D_i^-(u_s)=D_i(u_t^-)=0$ and \cite[Lemma 6.2]{Kh4} applies.

In the general case by the Heyneman-Radford theorem (see \cite[Corollary 5.4.7]{Mon} or very ``skew 
primitive" version  \cite[Corollary 5.3]{Kh2}) the Hopf ideal $J$ has a nonzero
skew primitive element provided that $J\neq 0.$ Denote by $J_1$ an ideal generated
by all skew primitive elements of $J.$ Clearly $J_1$ is a Hopf ideal. 
Since all homogeneous components of a skew primitive 
element are skew primitive, the Hopf ideal $J_1$ is homogeneous.   Moreover, we have
$D_i^-(u_s)=0,$ $s\in S_1$ where $u_s$ run through the set of all homogeneous skew primitive 
elements of $J.$ Now consider the Hopf ideal $J/J_1$ of the quotient Hopf algebra
$G\langle X\rangle /J_1.$ If $J\neq J_1$ then this ideal also has nonzero skew primitive elements.
Denote by $J_2/J_1$ the ideal generated by all skew primitive elements of $J/J_1,$
where $J_2$ is its preimage with respect to the natural homomorphism
$G\langle X\rangle \rightarrow G\langle X\rangle /J_1.$ Again we have 
$D_i^-(\bar{u}_s)=0,$ $s\in S_2$ in $G\langle X\rangle /J_1,$ where $\bar{u}_s$ 
run through the set of all homogeneous skew primitive elements of $J/J_1.$
In particular the ideal $J_2$ has a set of generators $u_s,$ $s\in S_1\cup S_2$
such that $D_i^-(u_s)\in J_1.$ Continuing the process we shall find a set of generators 
$u_s,$ $s\in S_1\cup S_2\cup S_3\cup \ldots $ for $J$ such that $D_i^-(u_s)\in J,$ all $s.$

In perfect analogy we find a set of generators  $u_t^-,$ 
for $J^-$ such that $D_i(u_t^-)\in J^-,$ all $t.$ Hence \cite[Lemma 6.2]{Kh4} applies. 
\end{proof}

\noindent
{\bf Remark}. In the proof we do not use relations (\ref{KM1}) on the quantification parameters, while relations (\ref{shar1}), (\ref{shar2}) on the characters are essential. Also we are reminded that originally the maps $D_i,$ $D_i^-$ were defined so that $D_i(u^-)=[x_i,u^-],$ $D_i^-(u)=[u,x_i^-]$ in the algebra 
$H\langle X\cup X^-|| [x_i,x_j^-]=\delta _i^j(1-g_if_i)\rangle .$ Hence equalities $D_i^-=\tilde{D}_i^-$
and (\ref{sqi2}) imply a differential representation:
\begin{equation}
[x_i,u^-]=\partial_{-i}^*(u^-)p(x_i,u^-)p_{ii}^{-1}-g_if_i \partial_{-i}(u^-),\ \ \ u^-\in{\bf k}\langle X^-\rangle ,
\label{sqi3}
\end{equation}
\begin{equation}
[u, x_i^-]= \partial_ i^*(u)-p_{ii}^{-1}p(u,x_i)\partial_i(u)g_if_i, \ \ \ u\in{\bf k}\langle X\rangle .
\label{sqi4}
\end{equation}

\section{PBW-generators for coideal subalgebras in 
$U_q^+(\frak{sl}_{n+1}),$ $u_q^+(\frak{sl}_{n+1})$}
Suppose that the quantification parameters $p_{ij}=p(x_i,x_j)=\chi ^{i}(g_{j})$
are connected by the following relations 
\begin{equation}
p_{ii}=q; \ \ p_{i\, i+1}p_{i+1\, i}=q^{-1}; \ \ 
p_{ij}p_{ji}=1,\  |i-j|>1,
\label{a1rel}
\end{equation}
where $q\not= \pm 1.$ By definition 
 $U_q^+(\frak{ sl}_{n+1})$ as a character Hopf algebra is set up by the relations
\begin{equation}
[[x_i,x_{i+1}],x_{i+1}]=[x_i,[x_i,x_{i+1}]]=[x_i,x_j]=0, \ \ |i-j|>1.
\label{rela}
\end{equation}
The structure of this algebra is defined by the following theorem (see, \cite[Theorem $A_n$]{Kh4},
or in other terms \cite[Lemmas pp. 176, 184]{CM}, \cite{CLMT}). 
Recall that $u(k,m)$ $=x_kx_{k+1}\ldots x_m,$ $k\leq m.$
The standard word $u(k,m)$ defines a super-letter
$u[k,m]\stackrel{df}{=}$ $[x_k[x_{k+1}[\ldots[x_{m-1}x_m]\ldots ]]],$
while by definition $u[k,k]=x_k.$ Of course the value of $u[k,m]$ in
$U_q^+(\frak{ sl}_{n+1})$ is independent of the alignment of brackets (Lemma \ref{ind}). 
By $U_q^+(\frak{ sl}_{n+1})^{(2)}$
we denote the ideal $\xi (G\langle X\rangle ^{(2)})$ 
generated by the values of $x_ix_j,$ $1\leq i,j\leq n.$
\begin{theorem}
$1.$ The values of the super-letters $u[k,m],$ $1\leq k\leq m\leq n$ in $U_q^+(\frak{ sl}_{n+1})$
form the set of PBW-generators of $U_q^+(\frak{ sl}_{n+1})$ over ${\bf k}[G].$ All heights are infinite.

$2.$ If $q$ is not a root of unity, then the ideal $U_q^+(\frak{ sl}_{n+1})^{(2)}$ has no nonzero skew-primitive elements.

$3.$ If $u$ is a standard word then either $u=u(k,m)$ or $[u]=0$ in $U_q^+(\frak{ sl}_{n+1}).$
\label{strA}
\end{theorem}
According to the Heyneman-Radford theorem (see \cite{HR}, or \cite[Corollary 5.4.7]{Mon}) every non-zero bi-ideal of a character Hopf algebra always has nonzero skew primitive elements. By this reason the second statement 
implies that Ker$\, \xi ,$ $\  \xi :G\langle X\rangle \rightarrow U_q^+(\frak{ sl}_{n+1}),$ is the biggest
Hopf ideal in $G\langle X\rangle ^{(2)}.$ In particular one can apply Lemma \ref{qsim}
to $U_q^+(\frak{ sl}_{n+1})$ provided $q$ is not a root of 1.

If $q$ is a root of 1 ($q\neq \pm 1$), then by definition
$u_q^+(\frak{ sl}_{n+1})$ is a quotient $U_q^+(\frak{ sl}_{n+1})/{ \Lambda },$ where
$\Lambda $ is the biggest Hopf subideal of $U_q^+(\frak{ sl}_{n+1})^{(2)}.$
Hence we may apply Lemma \ref{qsim} to $u_q^+(\frak{ sl}_{n+1})$ as well.

\smallskip
If {\bf U} is a right coideal subalgebra of $U_q^+(\frak{ sl}_{n+1})$
that contains {\bf k}$\, [G],$ then by Proposition \ref{pro} and Theorem \ref{strA} it has 
PBW-generators of the form  (\ref{vad1}):
\begin{equation}
c_u=u^s+\sum \alpha _iW_i+\ldots \in \hbox{\bf U}, \ \ \ u=u[k,m].
\label{vad25}
\end{equation}
By means of relations (\ref{a1rel}) we have
$p_{uu}=p(x_kx_{k+1}\cdots x_m,$ $ x_kx_{k+1}\cdots x_m)=q.$
Thus, if $q$ is not a root of 1, Lemma \ref{nco1}
shows that in (\ref{vad25})  the exponent $s$ equals 1, 
while  all heights of the $c_u$'s in {\bf U} are infinite.

If $q$ has a finite multiplicative order $t>2,$ then $u[k,m]^t=0$ in $u_q^+(\frak{ sl}_{n+1}),$
see for example \cite[Theorem 3.2]{KA}; that is, by \cite[Lemma 3.3]{KA}
the values of $u[k,m]$ are still the 
PBW-generators of $u_q^+(\frak{ sl}_{n+1}),$ but all of them have the finite height $t.$
By Lemma \ref{nco1} in (\ref{vad25}) we have $s\in \{ 1, t, tl^r\}.$ 
Since $u[k,m]^t=u[k,m]^{tl^r}=0,$ the exponent $s$
in (\ref{vad25})  equals 1, while all heights of the $c_u$'s in {\bf U} equal $t.$

Hence in both cases the PBW-generators of {\bf U} have the following form
\begin{equation}
c_u=u[k,m]+\sum \alpha _iW_i+\sum_j \beta_jV_j \in \hbox{\bf U}.
\label{vad22}
\end{equation}
where $W_i$ are the basis super-words starting with less than $u[k,m]$ super-letters, 
$D(W_i)=$ $D(u[k,m])$ $=x_k+x_{k+1}+\ldots +x_m,$ and $V_j$ are $G$-super-words 
of $D$-degree less than $x_k+x_{k+1}+\ldots +x_m.$

Now, in order to reduce the freedom in construction of the PBW-generators, 
we are going to show that  (\ref{vad22}) with homogeneous $c_u$ implies that 
{\bf U} has an element of the same form that belongs to
a special finite set of elements (\ref{cbr1}).

\begin{proposition} 
If a right coideal subalgebra {\bf U}$\supseteq {\bf k}[G]$ of $U_q^+(\frak{ sl}_{n+1})$
or $u_q^+(\frak{ sl}_{n+1})$ contains a homogeneous element $c$
with the leading term $u[k,m],$ $k\leq m,$ then for a suitable subset $\bf S$
of the interval $[k,m-1]$ the value of the below defined element
$\Psi ^{\hbox{\bf s}}(k,m)$ belongs to {\bf U}.
\label{phi}
\end{proposition}

\begin{definition} \rm
Let {\bf S} be a set of integers from the interval $[1,n].$ We define a
{\it piecewise continuous word related to} {\bf S} as follows
\begin{equation}
u^{\hbox{\bf s}}(k,m) 
\stackrel{df}{=}u(1+s_r,m)u(1+s_{r-1},s_r)\cdots u(1+s_1,s_2)u(k,s_1),
\label{pdw}
\end{equation}
where {\bf S}$\, \cap [k,m-1]=\{ s_1,s_2,\ldots s_r\}, \ $
$k\leq s_1<s_2<\ldots <s_{r-1}<s_r<m.$

If the pair $(k,m)$ \underline{is fixed}, 
we denote ${\bf S}_{\circ }=\, (${\bf S}$\, \cap [k,m-1])\, \cup \, \{ k-1 \} ,$
while ${\bf S}^{\bullet }=\, (${\bf S}$\, \cap [k,m-1])\, \cup \, \{ m \} ;$
respectively we extend $s_0=k-1\in \, {\bf S}_{\circ },$ and 
$s_{r+1}=m\in \, {\bf S}^{\bullet }.$ By $\overline{\bf S}$ we denote a complement of {\bf S}
with respect to $[k,m-1].$

We define a bracketing of a piecewise continuous word by
$$
\Psi ^{\hbox{\bf s}}(k,m)\stackrel{df}{=}
$$
\begin{equation}
 \hbox{\Large [[}\ldots \hbox{\Large [}u[1+s_r,m], u[1+s_{r-1},s_r]\hbox{\Large ]}, \ldots \ ,
u[1+s_1,s_2]\hbox{\Large ]},\, u[k,s_1]\hbox{\Large ]}.
\label{cbr1}
\end{equation}
The leading term of 
$\Psi ^{\hbox{\bf s}}(k,m)$ in the PBW-decomposition 
given in Theorem \ref{strA} is proportional to $u[k,m].$
In particular $\Psi ^{\hbox{\bf s}}(k,m)$ has the form (\ref{vad22}) up to a scalar factor.
\label{eski}
\end{definition}

The elements $y_i=u[1+s_{r-i+1},s_{r-i+2}],$ $1\leq i\leq r+1,$ satisfy
$p(y_i,y_j)p(y_j,y_i)=1,$ $[y_i,y_j]=0$ provided that $|i-j|>1.$
Thus by Lemma \ref{ind} applied to (\ref{cbr1}) the value in $U_q^+(\frak{ sl}_{n+1})$
or $u_q^+(\frak{ sl}_{n+1})$ of the bracketing is independent
of the alignment of the big brackets. In particular we have the following decomposition
\begin{equation}
 \Psi ^{\hbox{\bf s}}(k,m)=\hbox{\Large [} \Psi ^{\hbox{\bf s}}(1+s_i,m),
\Psi ^{\hbox{\bf s}}(k,s_i) \hbox{\Large ]}.
\label{cbrr}
\end{equation}
Of course the word $u(k,m)$ is a piecewise continuous word with empty {\bf S}, 
or more generally, with {\bf S}$\, \cap [k,m-1]=\emptyset .$

Let us choose an arbitrary $s_{r+1},$ $s_r<s_{r+1}<m,$ and denote
$$
u=u[1+s_r,s_{r+1}], \ \ v=u[1+s_{r+1},m].
$$
Then by (\ref{alin}) and decomposition (\ref{cbrr})
we have
$$
\Psi ^{\hbox{\bf s}}(k,m)=\hbox{\Large [} [u,v],
\Psi ^{\hbox{\bf s}}(k,s_r)\hbox{\Large ]},
$$
while
$$
\Psi ^{\hbox{\bf s}\cup \{ s_{r+1}\}}(k,m)=\hbox{\Large [} [v,u],
\Psi ^{\hbox{\bf s}}(k,s_r) \hbox{\Large ]}.
$$
Since $[v,u]=-p_{vu}[u,v]+(1-p_{uv}p_{vu})v\cdot u$ and $p_{uv}p_{vu}=q^{-1},$
using formula (\ref{br2}), we get the following recurrence relation for the complete bracketing:
\begin{equation}
\Psi ^{\hbox{\bf s}\cup \{ s_{r+1}\} }(k,m)=
(1-q^{-1})v\cdot \Psi ^{\hbox{\bf s}}(k,s_{r+1})-p_{vu}\Psi ^{\hbox{\bf s}}(k,m),
\label{cbr3}
\end{equation}
where as above $u=u[1+s_r,s_{r+1}],$ $v=u[1+s_{r+1},m].$

More generally, if $s_{i-1}<t<s_i,$ and we denote $u=u[1+s_{i-1},t],$ $v=u[1+t, s_i,]$
then (\ref{cbr3}) reads
$$
\Psi ^{\hbox{\bf s}\cup \{ t\} }(k,s_i)=
(1-q^{-1})v\cdot \Psi ^{\hbox{\bf s}}(k,t)-p_{vu}\Psi ^{\hbox{\bf s}}(k,s_i).
$$
Let us commute this equality with $\Psi ^{\hbox{\bf s}}(1+s_i,m)$ from the left.
Then decomposition (\ref{cbrr}) and (\ref{br1}) imply a general formula
$$
\Psi ^{\hbox{\bf s}\cup \{ t\} }(k,m)=
(1-q^{-1})[\Psi ^{\hbox{\bf s}}(1+s_i,m),v]\cdot \Psi ^{\hbox{\bf s}}(k,t)-p_{vu}\Psi ^{\hbox{\bf s}}(k,m)
$$
\begin{equation}
=(1-q^{-1})\Psi ^{\hbox{\bf s}}(1+t,m)\cdot \Psi ^{\hbox{\bf s}}(k,t)-p_{vu}\Psi ^{\hbox{\bf s}}(k,m).
\label{cby}
\end{equation}

\smallskip
Denote by $\varphi $ a map from $[1,n]$ to $[1,n]$ given by $\varphi (i)=n-i+1.$
One may replace the main skew primitive generators $x_i,$  $1\leq i\leq n$
with $y_i,$  $1\leq i\leq n,$ where by definition $y_i=x_{\varphi (i)}.$
 Our basic concepts (Definition \ref{eski}) are not invariant
with respect to this replacement. For example, 
$$
\Psi ^{\emptyset }(1,n)=[x_1x_2\cdots x_n]=[y_ny_{n-1}\cdots y_1]=\Psi ^{[1,n-1]}_y(1,n).
$$ 
This fact signifies that the application of already proved formulae to new generators
ought to provide additional information. To get this information we need the following
{\it decoding} lemma. In what follows by $\sim $ we denote the projective equality:
$a\sim b$ if and only if $a=\alpha b,$ where $0\neq \alpha \in {\bf k}.$

\begin{lemma}  
We have 
\begin{equation}
\Psi ^{\hbox{\bf s}}(k,m)\sim \Psi ^{\overline{\varphi(\hbox{\bf s})-1}}_y(\varphi (m),\varphi (k)),
\label{decod}
\end{equation}
where by $\varphi ({\bf S})-1$ we denote $\{ \varphi (s)-1\, |\, s\in {\bf S}\} ,$
while the complement is related to the pair $(\varphi (m),\varphi (k)),$
see Definition $\ref{eski}.$
\label{dec}
\end{lemma}
\begin{proof} We use induction on $m-k.$ If $m=k,$ the equality reduces to $x_k=y_{\varphi (k)}.$
To make the inductive step we consider two cases.

a). $m-1\in {\bf S}.$ In this case $\varphi (m-1)\in \varphi ({\bf S}).$ Since 
$\varphi (m-1)=\varphi (m)+1,$ we get $\varphi (m)\in \varphi ({\bf S})-1.$
In particular we have an equality of sets:
\begin{equation}
[\varphi (m), \varphi (k)-1]\setminus \{ \varphi ({\bf S})-1\}=
[\varphi (m)+1, \varphi (k)-1]\setminus \{ \varphi ({\bf S})-1\}.
\label{kra}
\end{equation}
Using the inductive supposition we have
$$
\Psi ^{\hbox{\bf s}}(k,m)\sim [x_m, \Psi ^{\hbox{\bf s}}(k,m-1)]=
[y_{\varphi (m)}, \Psi ^{\overline{\varphi(\hbox{\bf s})-1}}_y(\varphi (m)+1,\varphi (k))].
$$
If set (\ref{kra}) is empty, the above element equals the element
$u_y[\varphi (m),\varphi (k)]$ that, in the context, coincides with 
$\Psi ^{\overline{\varphi(\hbox{\bf s})-1}}_y(\varphi (m),\varphi (k)).$
Otherwise denote by $t_1$ the minimal element in (\ref{kra}).
Using definition (\ref{cbr1}) we may continue
$$
=\left[y_{\varphi (m)}, [\Psi ^{\overline{\varphi(\hbox{\bf s})-1}}_y(1+t_1,\varphi (k)), 
u_y[\varphi (m)+1,t_1]] \right] .
$$ 
Denote $u=y_{\varphi (m)},$ $v=\Psi ^{\overline{\varphi(\hbox{\bf s})-1}}_y(1+t_1,\varphi (k)),$
and $w=u_y[\varphi (m)+1,t_1].$ Then $[u,v]=0$ since $u,v$ are separated.
By the same reason $p_{uv}p_{vu}=1,$ see (\ref{a1rel}). Hence the
conditional Jacobi identity (\ref{jak4}) implies
$$
[u,[v,w]]\sim [v,[u,w]]=\Psi ^{\overline{\varphi(\hbox{\bf s})-1}}_y(\varphi (m),\varphi (k)).
$$

b). $m-1\notin {\bf S}.$ In this case $s_r<m-1,$ and we have
$$
\Psi ^{\hbox{\bf s}}(k,m)=[u[1+s_r,m], \Psi ^{\hbox{\bf s}}(k,s_r)]
=\left[ [x_{1+s_r},u[2+s_r,m]], \Psi ^{\hbox{\bf s}}(k,s_r) \right] .
$$
Let us denote $u=x_{1+s_r},$ $v=u[2+s_r,m],$ $w=\Psi ^{\hbox{\bf s}}(k,s_r).$
Then $[v,w]=0,$ $p_{vw}p_{wv}=1$ since $v$ and $w$ are separated. Thus the conditional
Jacobi identity (\ref{jak1}) implies $[[u,v],w]\sim [[u,w],v].$ Since 
$[u,w]=\Psi ^{\hbox{\bf s}}(k,1+s_r),$ we may apply the inductive supposition
\begin{equation}
[[u,w],v]=\left[ \Psi ^{\overline{\varphi(\hbox{\bf s})-1}}_y(\varphi (s_r)-1,\varphi (k)), 
[y_{\varphi (s_r)-2}y_{\varphi (s_r)-3}\cdots y_{\varphi (m)}]\right] . 
\label{to}
\end{equation}
We know that {\bf S} has no one of the points $1+s_r,$ $2+s_r,$ $\ldots ,$ $m-1,$ hence
the set $\varphi({\bf S})-1$ has no one of the points 
$\varphi (s_r)-2, \varphi (s_r)-3,$ $\ldots , \varphi (m)+1, \varphi (m),$
and we have the equality of sets
$$
\{ \varphi (m), \varphi (m)+1, \ldots \varphi (s_r)-2 \} \cup
\left( [\varphi (s_r)-1,\varphi (k)-1]\setminus \{ \varphi({\bf S})-1\} \right)
$$
$$
=[\varphi (m), \varphi (k)-1]\setminus \{ \varphi({\bf S})-1\} .
$$
Therefore the right hand side of (\ref{to}) takes up the form of the right hand side of 
(\ref{decod}). The decoding lemma is proved. \end{proof}

\begin{lemma}  
If $t\notin {\bf S},$ $k\leq t<m,$ then
\begin{equation}
 \Psi ^{\hbox{\bf s}}(k,m)\sim \left[ \Psi ^{\hbox{\bf s}}(k,t),
\Psi ^{\hbox{\bf s}}(1+t,m) \right] .
\label{cbry}
\end{equation}
\label{dy}
\end{lemma}
\begin{proof} One should replace variables by means of (\ref{decod}), then apply
(\ref{cbrr}), and replace back the variables by (\ref{decod}). \end{proof}

In particular we have the following formula, convenient for induction.
\begin{equation}
\Psi ^{\hbox{\bf s}}(k,m)\sim \left\{ 
\begin{matrix}[x_m,\Psi ^{\hbox{\bf s}}(k,m-1)], & \hbox{ if } m-1\in {\bf S};\hfill \cr
[\Psi ^{\hbox{\bf s}}(k,m-1),x_m], \hfill & \hbox{ if } m-1\notin {\bf S}.
\end{matrix}
\right.
\label{cin}
\end{equation}

Symmetrically relations (\ref{cbrr}) and (\ref{cbry}) imply
\begin{equation}
\Psi ^{\hbox{\bf s}}(k,m)\sim \left\{ 
\begin{matrix}[x_k,\Psi ^{\hbox{\bf s}}(k+1,m)], & \hbox{ if } k\notin {\bf S};\hfill \cr
[\Psi ^{\hbox{\bf s}}(k+1,m),x_k], \hfill & \hbox{ if } k\in {\bf S}.
\end{matrix}
\right.
\label{cin1}
\end{equation}

\smallskip
Our immediate task is to find a differential subspace generated by $\Psi ^{\hbox{\bf s}}(k,m).$
To attain these ends we display the element $\Psi ^{\hbox{\bf s}}(k,m)$ graphically 
as a sequence of black and white points labeled by the numbers
$k-1,$ $k,$ $k+1, \ldots ,$ $m-1,$ $m,$ where the first point is always white, and
the last one is always black, while an intermediate point labeled by $i$ is black if and only if 
$i\in {\bf S}:$  
\begin{equation}
\stackrel{k-1}{\circ } \ \ \stackrel{k}{\circ } \ \ \stackrel{k+1}{\circ } 
\ \ \stackrel{k+2}{\bullet }\ \ \ \stackrel{k+3}{\circ }\ \cdots
\ \ \stackrel{m-2}{\bullet } \ \ \stackrel{m-1}{\circ }\ \ \stackrel{m}{\bullet }
\label{gra}
\end{equation}
For example the primitive generator $x_i$ is displayed by two dots
$$
\stackrel{i-1}{\circ }\ \ \stackrel{i}{\bullet }
$$
The element $[x_kx_{k+1}\cdots x_{m-1}x_m]$ is pictured as the whitest sequence 
$$
\stackrel{k-1}{\circ } \ \ \stackrel{k}{\circ } \ \ \stackrel{k+1}{\circ } 
\ \ \stackrel{k+2}{\circ }\ \ \ \stackrel{k+3}{\circ }\ \cdots
\ \ \stackrel{m-2}{\circ } \ \ \stackrel{m-1}{\circ }\ \ \stackrel{m}{\bullet }
$$
The most black sequence
$$
\stackrel{k-1}{\circ } \ \ \stackrel{k}{\bullet  } \ \ \stackrel{k+1}{\bullet  } 
\ \ \stackrel{k+2}{\bullet }\ \ \ \stackrel{k+3}{\bullet }\ \cdots
\ \ \stackrel{m-2}{\bullet } \ \ \stackrel{m-1}{\bullet }\ \ \stackrel{m}{\bullet }
$$
corresponds to the ``dual" element $[x_{m}x_{m-1}\cdots x_{k+1}x_k].$

In what follows we denote by $W^{\hbox{\bf s}}(k,m)$ the set of all elements
that are displayed by subsequences of the sequence (\ref{gra})
related to $\Psi ^{\hbox{\bf s}}(k,m):$
\begin{equation}
\{ \Psi ^{\hbox{\bf s}}(a,b), \, |\,  k\leq a\leq b\leq m;\
   b\in \hbox{\bf S} \hbox{ or } b=m;\ a-1\notin \hbox{\bf S} \hbox{ or } a=k \}.
\label{pbse}
\end{equation}

The following theorem shows that the differential subspace generated by 
an element displayed by (\ref{gra}) is spanned by all elements corresponding
to subsequences of (\ref{gra}) and their separated products.
\begin{definition}$\!\!\!$ \rm
If $w$ is a word in $X,$ we define a differential operator $D_w$ by the recurrence
formula $D_{x_{k}u}=\partial_k \circ D_u,$ $D_{\emptyset }=\, $id.
\label{oper}
\end{definition}
\begin{theorem}
The differential subspace generated by $\Psi ^{\hbox{\bf s}}(k,m)$ in $U_q^+(\frak{ sl}_{n+1})$
or $u_q^+(\frak{ sl}_{n+1})$
is spanned by the values of $W^{\hbox{\bf s}}(k,m)$ and by the values of all products of pairwise 
separated $($hence $q$-commuting, see Definition $\ref{sep})$ 
elements from $W^{\hbox{\bf s}}(k,m).$
\label{26}
\end{theorem}
\begin{proof}
We shall need the following properties of the partial derivations.
If $u$ is independent of $x_j$ (or, more generally, if $\partial_j (u)=0$) then
\begin{equation}
\partial_j ([u,v])=p(u,x_j)\left[ u, \partial_j (v)\right] ,\ \ 
\partial_j ([u,x_j])=0.
\label{dif0} 
\end{equation}
If $\partial_j (v)=0,$  then
\begin{equation}
\partial_j ([u,v])=\partial_j (u)\cdot v-p_{uv}p(v,x_j)v\cdot \partial_j (u);
\label{dif00}
\end{equation}
if additionally $p(x_j,v)p(v,x_j)=1$ (and still $\partial_j (v)=0$), then
\begin{equation}
\partial_j ([u,v])=\left [\partial_j (u), v\right] .
\label{dif000}
\end{equation}
These properties are a straightforward consequence of the definition (\ref{sqo})
and skew differential Leibniz rule (\ref{defdif}). Indeed,
$$
\partial_j ([u,v])=\partial_j (uv)-p_{uv}\partial_j (vu)
$$
$$
=\partial_j (u)\cdot v+p(u,x_j)u\cdot \partial_j (v)-
p_{uv}\partial_j (v)\cdot u-p_{uv}p(v,x_j)v\cdot \partial_j (u),
$$
which easely implies (\ref{dif0} -- \ref{dif000}) since $p_{uv}=p(u,x_j)p(u,\partial_j (v))$
provided that $v$ is dependent on $x_j,$ and $p_{uv}=p(\partial_j(u), v)p(x_j,v)$
provided that $u$ is dependent on $x_j.$

We shall prove by induction on $m-k$ the following formula
\begin{equation}
\partial_j (u[k,m])=\left\{ \begin{matrix}0,\hfill & j\neq k;\hfill \cr
(1-q^{-1})u[k+1,m], \hfill & j=k<m; \cr
1, \hfill & \hfill j=k=m. \end{matrix} \right.
\label{der1}
\end{equation}

If $k=m,$ the formula follows from definition (\ref{defdif}). Let $k<m.$
According to recurrence definition (\ref{alin}) we have  $u[k,m]=[x_k,u[k+1,m]].$ 

If $j=k,$ then $\partial_j (u[k+1,m])=0$ since $u[k+1,m]$ is independent of $x_k.$
Hence $(\ref{dif00})$ with $u\leftarrow x_k,$ $v\leftarrow u[k+1,m]$ implies
$$
\partial_j (u[k,m])=v-p(x_k,v)p(v,x_k)\cdot v.
$$
By means of (\ref{a1rel}) we have 
$p(x_k,v)p(v,x_k)=p_{k\, k+1}p_{k\, k+2}\ldots p_{k\, m}\cdot p_{k+1\, k}p_{k+2\, k}\ldots p_{mk}$
$=p_{k\, k+1}p_{k+1\, k}=q^{-1}.$

If $j\neq k,$ relation (\ref{dif0}) with $u\leftarrow x_k,$ $v\leftarrow u[k+1,m]$ implies
$$
\partial_j (u[k,m])=p(u,x_k)[x_k,\partial_j (u[k+1,m])].
$$
The inductive supposition 
yields $\partial_i (u[k+1,m])=\alpha u[k+2,m],$ $\alpha \in {\bf k}.$ Since $x_k$ and 
$u[k+2,m]$ are separated, they $q$-commute (this is true even if $u[k+2,m]$ is empty:
$[x_k,1]=x_k\cdot 1-\chi^{k}({\rm id})1\cdot x_k=0).$ Thus $[x_k,\partial_i (u[k+1,m])]=0,$
which completes the proof of (\ref{der1}).

Since $\Psi ^{\hbox{\bf s}}(k,m)$ is a linear combination of
products of $u[1+s_{i-1},s_i],$ $1\leq i\leq r+1,$ formula (\ref{der1}) implies
\begin{equation}
\partial_j(\Psi ^{\hbox{\bf s}}(k,m))=0,
 \hbox{ if } j-1\notin \, \hbox{\bf S}_{\circ }\stackrel{df}{=}({\bf S}\cap [k,m-1])\cup \{ k-1\} .
\label{der35}
\end{equation}

Let us consider remaining partial derivatives $\partial _j$ when $j-1\in {\bf S}_{\circ },$
that is $j=k$ or $j=1+s_i,$ $s_i\in {\bf S}\, \cap [k,m-1].$

Suppose $j=k.$ Then definition (\ref{cbr1}) implies 
$\Psi ^{\hbox{\bf s}}(k,m)=[u,v],$ where $u=\Psi ^{\hbox{\bf s}}(1+s_1,m),$
$v=u[k,s_1].$ Since $u$ is independent of $x_k,$ relation (\ref{dif0})
yields $\partial_k([u,v])=p(u,x_k)[u,\partial_k(v)].$ Hence by 
(\ref{der1}) we get
\begin{equation}
\partial_k (\Psi ^{\hbox{\bf s}}(k,m))=\left\{ 
\begin{matrix}\lambda \Psi ^{\hbox{\bf s}}(k+1,m), \hfill & \hbox{if }  s_1\neq k<m; \hfill \cr
0, \hfill & \hbox{if } s_1=k<m,\hfill  \end{matrix} \right.
\label{der4}
\end{equation}
provided that {\bf S}$\, \cap [k,m-1]\neq \emptyset ,$ where $\lambda =(1-q^{-1})p(u(1+s_1,m),x_k).$

Let $j=1+s_r.$ Recurrence formula (\ref{cbr3}) reads 
$$
\Psi ^{\hbox{\bf s}}(k,m)=(1-q^{-1})v\cdot \Psi ^{\hbox{\bf s}}(k,s_r)
-p_{uv}\Psi ^{\hbox{\bf s}\setminus \{ s_r\}}(k,m),
$$
where $u=u[1+s_{r-1}, s_r],$ $v=u[1+s_r, m].$ By (\ref{der35}) we have 
$\partial_j (\Psi ^{\hbox{\bf s}\setminus \{ s_r\}}(k,m))=0.$ Since 
$\Psi ^{\hbox{\bf s}}(k,s_r)$ is independent of $x_j,$ 
by means of (\ref{der1}) and skew differential Leibniz rule (\ref{defdif}) we have 
$$
\partial_{1+s_r} (\Psi ^{\hbox{\bf s}}(k,m))=
$$
\begin{equation}
=\left\{ \begin{matrix}(1-q^{-1})^2 u[2+s_r ,m]\cdot 
\Psi ^{\hbox{\bf s}}(k,s_r),\hfill & \hbox{if }1+s_r\neq m; \hfill \cr
(1-q^{-1}) \Psi ^{\hbox{\bf s}}(k,s_r),\hfill & \hbox{if }1+s_r=m.\hfill \end{matrix} \right. 
\label{der8}
\end{equation}

Let $j=1+s_i,$ $1\leq i<r.$ By (\ref{cbrr}) we have $\Psi ^{\hbox{\bf s}}(k,m)=[u,v],$
where $u=\Psi ^{\hbox{\bf s}}(1+s_i,m),$ $v=\Psi ^{\hbox{\bf s}}(k,s_i).$
By (\ref{der4}) the elements $\partial_j (u)$ and $v$ are separated, hence 
$v\cdot \partial_j (u)=p(v, \partial_j(u))\partial_j (u)\cdot v.$
Since $v$ is independent of $x_j,$ one may apply (\ref{dif00}).
We have $\partial_j([u,v])=\partial_j(u)\cdot v(1-p_{uv}p(v,x_j)p(v,\partial_j(u))).$
Due to (\ref{a1rel}) we get $p_{uv}p(v,x_j)p(v,\partial_j(u))$ $=p_{uv}p_{vu}=q^{-1}.$
Thus formula (\ref{der4}) implies
$$
\partial_{1+s_i} (\Psi ^{\hbox{\bf s}}(k,m))=
$$
\begin{equation}
=\left\{ \begin{matrix}\mu \Psi ^{\hbox{\bf s}}(2+s_i,m)\cdot 
\Psi ^{\hbox{\bf s}}(k,s_{i}),\hfill & \hbox{if }s_{i+1}>1+s_i; \hfill \cr
0,\hfill &\hbox{if }s_{i+1}=1+s_i,\hfill \end{matrix} \right.
\label{der7}
\end{equation}
where 
$\mu =(1-q^{-1})^2p(u(1+s_{i+1},m),x_{1+s_i}).$ 

Formulae (\ref{der35} -- \ref{der7}) show that products of pairwise separated
elements from $W^{\hbox{\bf s}}(k,m)$ span a differential subspace, that contains
all first derivatives of $\Psi ^{\hbox{\bf s}}(k,m).$ Hence
by induction  it contains the derivatives of higher order as well. 

To see that any  product of pairwise separated
elements from $W^{\hbox{\bf s}}(k,m)$ is proportional to some derivative
of $\Psi ^{\hbox{\bf s}}(k,m)$ we shall prove the following relation. 
\begin{equation}
\Psi ^{\hbox{\bf s}}(k,m)\cdot  D_w=\alpha \in {\bf k},\, \alpha \neq 0,\
\hbox{ if } \ w=u^{\hbox{\bf s}}(k,m).
\label{der9}
\end{equation}
If {\bf S}$\, \cap [k,m-1]=\emptyset $ then $w=x_kx_{k+1}\ldots x_m$ and relation follows from (\ref{der1}).
Let {\bf S}$\, \cap [k,m-1]\neq \emptyset .$  By definition (\ref{pdw}) we have $w=v\cdot w^{\prime }, $
 where  $v=x_{1+s_r}x_{2+s_r}\ldots x_m,$ $w^{\prime }=u^{\hbox{\bf s}}(k,s_r).$
Hence  
$$
\Psi ^{\hbox{\bf s}}(k,m)\cdot  D_w=\partial_{1+s_r}(\Psi ^{\hbox{\bf s}}(k,m))
\cdot D_{v^{\prime }}D_{w^{\prime }},
$$
where $v^{\prime }=x_{2+s_r}\ldots x_m.$ By (\ref{der8})
the element $\partial_{1+s_r}(\Psi ^{\hbox{\bf s}}(k,m))$ is proportional
to $u[2+s_r]\cdot \Psi ^{\hbox{\bf s}}(k,s_r).$ Since $\Psi ^{\hbox{\bf s}}(k,s_r)$
is independent of $x_j,$ $2+s_r\leq j\leq m,$ skew differential Leibniz rule (\ref{defdif}) implies
$$
(u[2+s_r]\cdot \Psi ^{\hbox{\bf s}}(k,s_r))\cdot D_{v^{\prime }}=
(u[2+s_r]\cdot D_{v^{\prime }})\Psi ^{\hbox{\bf s}}(k,s_r).
$$
By means of the multiple application of (\ref{der1}) we see that 
$\Psi ^{\hbox{\bf s}}(k,m)\cdot D_w$ is proportional to 
$\Psi ^{\hbox{\bf s}}(k,s_r)\cdot D_{w^{\prime }}.$ By induction on $m-k$ 
we get (\ref{der9}).

\smallskip
Now consider a product of separated elements from $W^{\hbox{\bf s}}(k,m),$
\begin{equation}
\Psi ^{\hbox{\bf s}}(a_1,b_1)\cdot  \Psi ^{\hbox{\bf s}}(a_2,b_2)\ldots \Psi ^{\hbox{\bf s}}(a_l,b_l),
\ \ k\leq b_i<a_{i+1}-1<m, 1\leq i\leq l.
\label{spe}
\end{equation}
We shall prove by induction on $l$ that there exists a word $w$ such that 
$\Psi ^{\hbox{\bf s}}(k,m)\cdot D_w$ is proportional to (\ref{spe}).

Assume $l=1,$ $\Psi ^{\hbox{\bf s}}(a,b)\in W^{\hbox{\bf s}}(k,m).$
Let us prove first that $\Psi ^{\hbox{\bf s}}(a,m)$ has the required representation.
If $a=k$ there is nothing to prove. Let $a>k.$ In this case by definition 
$a-1\notin \,${\bf S}, say $s_i<a-1<s_{i+1}$ for some $i,$ $0\leq i\leq r,$
where formally $s_0=k-1,$ $s_{r+1}=m.$ We have $s_{i+1}>1+s_i.$ Hence by (\ref{der7}),
or by (\ref{der4}) provided $i=0,$ or by (\ref{der8}) provided $i=r$,
the element $\partial _{1+s_i}(\Psi ^{\hbox{\bf s}}(k,m))$ is proportional to 
$\Psi ^{\hbox{\bf s}}(2+s_i,m)\cdot \Psi ^{\hbox{\bf s}}(k,s_i),$ where formally 
$\Psi ^{\hbox{\bf s}}(k,s_0)=1.$ Since 
$\Psi ^{\hbox{\bf s}}(2+s_i,m)$ is independent of $x_j,$ $k\leq j\leq s_i,$ formula
(\ref{der9})  shows that $\Psi ^{\hbox{\bf s}}(k,m)\cdot D_u,$ $u=x_{1+s_i}u^{\hbox{\bf s}}(k,s_i)$
is proportional to $\Psi ^{\hbox{\bf s}}(2+s_i,m).$ If $2+s_i=a,$ the required representation 
of $\Psi ^{\hbox{\bf s}}(a,m)$ is found. If $2+s_i<a$ then by means of a multiple use of
(\ref{der4}) we see that $\Psi ^{\hbox{\bf s}}(2+s_i,m)\cdot D_v,$ $v=x_{2+s_i}x_{3+s_i}\ldots x_{a-1}$
is proportional to $\Psi ^{\hbox{\bf s}}(a,m).$

Similarly we may find a word $w$ such that $\Psi ^{\hbox{\bf s}}(a,m)\cdot D_w$
is proportional to $\Psi ^{\hbox{\bf s}}(a,b).$ Indeed, if $b=m,$ we put $w=\emptyset .$
If $b<m,$ then by definition $b\in \,${\bf S}. Denote by $s_i$ an element from {\bf S} such that
$[b,s_i]\subseteq \,${\bf S}, and either $1+s_i\notin \,${\bf S} or $1+s_i=m.$

If $1+s_i\neq m$ then Eq. (\ref{der7}) implies that $\partial _{1+s_i}(\Psi ^{\hbox{\bf s}}(a,m))$
is proportional to $\Psi ^{\hbox{\bf s}}(2+s_i,m)\cdot \Psi ^{\hbox{\bf s}}(a,s_i).$
Since $\Psi ^{\hbox{\bf s}}(2+s_i,m)$ is independent of $x_j,$ $a\leq j\leq s_i,$ formula
(\ref{der9}) implies that $\Psi ^{\hbox{\bf s}}(a,m)\cdot D_y,$ 
$y=x_{1+s_i}u^{\hbox{\bf s}}(2+s-i,m)$ is proportional to $\Psi ^{\hbox{\bf s}}(a,s_i).$

If $1+s_i=m$ then by (\ref{der8}) the element $\partial _{1+s_i}(\Psi ^{\hbox{\bf s}}(a,m))$
itself is proportional to $\Psi ^{\hbox{\bf s}}(a,s_i).$

Finally a multiple use of (\ref{der8}) shows that 
$\Psi ^{\hbox{\bf s}}(a,s_i)\cdot D_z,$ $z=x_{s_i}x_{s_i-1}x_{s_i-2}\ldots x_{b+1}$
is proportional to $\Psi ^{\hbox{\bf s}}(a,b).$ This completes the case $l=1.$

Consider (\ref{spe}) with $l>1.$  By definition $1+b_1<a_2,$ and $b_1\in \,${\bf S}, 
$a_2-1\notin \,$ {\bf S} since obviously $b_1\neq m,$ $a_2\neq k.$ Denote by 
$s_i$ a number such that $[b,s_i]\subseteq \,${\bf S},  $1+s_i\notin \,${\bf S}.
Of course $1+s_i\leq a_2-1.$ Relation (\ref{der7}) implies that
$\partial _{1+s_i}(\Psi ^{\hbox{\bf s}}(k,m))$ is proportional to 
$\Psi ^{\hbox{\bf s}}(2+s_i,m)\cdot \Psi ^{\hbox{\bf s}}(k,s_i).$
We have $k\leq a_1\leq b_1\leq s_i,$ hence by the considered above case
there exists a word $v$ in $Y=\{ x_j\, |\, k\leq j<a_1 \hbox{ or } b_1<j\leq s_i\}$
such that $\Psi ^{\hbox{\bf s}}(k,s_i)\cdot D_v$ is proportional to 
$\Psi ^{\hbox{\bf s}}(a_1,b_1).$ Since $\Psi ^{\hbox{\bf s}}(2+s_i,m)$
is independent of $Y,$ skew differential Leibniz rule (\ref{defdif}) shows that 
$(\Psi ^{\hbox{\bf s}}(2+s_i,m)\cdot \Psi ^{\hbox{\bf s}}(k,s_i))\cdot D_v$
is proportional to $\Psi ^{\hbox{\bf s}}(2+s_i,m)\cdot \Psi ^{\hbox{\bf s}}(a_1,b_1).$
Further, by the inductive supposition there exists a word $u$ in
$Z=\{ x_j\, |\, 2+s_i\leq j\leq m\}$ such that $\Psi ^{\hbox{\bf s}}(2+s_i,m)\cdot D_u$
is proportional to $\Psi ^{\hbox{\bf s}}(a_2,b_2)\ldots \Psi ^{\hbox{\bf s}}(a_l,b_l).$
Now we see that $\Psi ^{\hbox{\bf s}}(k,m)\cdot D_{vu}$ is proportional to (\ref{spe})
because all factors in (\ref{spe}) (skew)commute.
\end{proof}

\smallskip
\noindent
{\it Proof} of Proposition \ref{phi}. 
Since $c$ is homogeneous, in representation (\ref{vad22}) there is no the second sum,
while by definition the $W_i$'s are monotonous products of the basis elements $u[a,b].$
In each $W_i$ exactly one factor depends on $x_k.$ This factor has the form $u[x_k,b],$
$k\leq b<m,$ and, as the maximal super-letter of the super-word $W_i,$ it is located at 
the end of $W_i.$ Hence representation (\ref{vad22}) of $c$ takes the form
\begin{equation}
c=u[k,m]+\sum _{i=k}^{m-1}A_iu[k,i],
\label{ana1}
\end{equation}
where  $A_i,$ $k\leq i<m$  is a linear combination of monotonous 
super-words of degree $x_{i+1}+x_{i+2}+\ldots +x_m.$
Each monotonous super-word of that degree has the form
\begin{equation}
u[1+l_p,m]\cdot u[1+l_{p-1},l_p]\cdot \ \ldots \ \cdot u[1+i,l_1] \stackrel{df}{=}u^{[L]}(1+i, m),
\label{pbr1}
\end{equation}
where $L=\{ l_1<l_2< \ldots <l_p\} $ $\subseteq [1+i, m-1].$
To derivate products of that type we shall prove the following general formula. Assume 
$T=\{ t_1<t_2< \ldots <t_s\}$ is another subset of $[1+i,m-1],$ and let $w=u^{T}(1+i,m),$
see definition (\ref{pdw}). In this case we have 
\begin{equation}
u^{[L]}(1+i,m)\cdot  D_w=\left\{ \begin{matrix}
\alpha , \hfill & \hbox{if } T\subseteq L;\hfill \cr
0 , \hfill & \hbox{otherwise,}\hfill 
\end{matrix}
\right.
\label{xy}
\end{equation}
where $\alpha \in \,${\bf k}, $\alpha \neq 0.$ Indeed, if $T=\emptyset $ then 
$w=x_{1+i}x_{2+i}\ldots x_m,$ while the required equality follows from
(\ref{der1}) and (\ref{defdif}). If $T\neq \emptyset $ then
definition (\ref{pdw}) implies
$w=v\cdot w^{\prime },$ where $v=x_{1+t_s}x_{2+t_s}\ldots x_m,$ $w^{\prime }=u^{T}(1+i,t_s).$ 
If $t_s\notin L$ then (\ref{der1}) implies 
$\partial _{1+t_s}(u^{[L]}(1+i,m))=0,$ hence  $u^{[L]}(1+i,m)\cdot D_w=0,$ which is required.
If $t_s\in L$ say $t_s=l_j,$
then again formulae (\ref{der1}) and (\ref{defdif}) imply $u^{[L]}(1+i,m)\cdot D_v\sim u^{[L]}(1+i,l_j).$
Therefore $u^{[L]}(1+i,m)\cdot D_w\sim u^{[L]}(1+i,t_s)\cdot D_{w^{\prime }},$
which proves (\ref{xy}) by evident induction on $m.$

Suppose that $A_{s_1},$ $A_{s_2},\ldots ,$ $A_{s_r}$ are all nonzero terms in (\ref{ana1}).
Denote {\bf S}$_0=\emptyset ,$ {\bf S}$_t=\{ s_1, s_2, \ldots s_t\} ,$
{\bf S}$_r=\, ${\bf S}. 
The decomposition (\ref{ana1}) shows that
\begin{equation}
\Psi ^{\hbox{\bf s}_t}(k,m)+
\sum_{i=s_{t+1}}^{m-1}A_i\Psi ^{\hbox{\bf s}_t}(k,i)\in \hbox{\bf U},
\label{pit1}
\end{equation}
where $t=0.$ We shall prove  in two steps that inclusion (\ref{pit1})
with a given $t,$ $0\leq t<r$  
implies the same type of inclusion with $t\leftarrow t+1$
and proportional $A_i.$ 
This will imply that (\ref{pit1}) with $t=r,$ which says 
$\Psi ^{\hbox{\bf s} }(k,m)\in \, ${\bf U},
is valid as well.

1.  The element $A_{s_{t+1}}$ is a linear combination of super-words 
(\ref{pbr1}) with $i\leftarrow s_{t+1}.$ Denote by $T$ the maximal with respect to inclusion 
subset of the interval $[1+s_{t+1},m-1],$ such that in the PBW-decomposition of $A_{s_{t+1}}$
the super-word $u^{[T]}(1+s_{t+1}, m)$ appears with a nonzero coefficient.
Let us apply $D_w,$ with $w=u^{T}(1+s_{t+1}, m)$ to (\ref{pit1}). Since $w$ is independent of $x_k,$
$x_{1+s_j},$ $1\leq j\leq t,$ formula (\ref{der35}) with skew differential Leibniz rule (\ref{defdif}) show that
\begin{equation}
B\cdot D_w=\sum _{i=s_{t+1}}^{m-1}
(A_i\cdot D_w)\Psi ^{\hbox{\bf s}_t}(k,i),
\label{s}
\end{equation}
where $B$ is the left hand side of (\ref{pit1}).
The constitution of $w$
contains the constitutions of $A_i,$ $i>s_{t+1}.$ Hence all terms in the sum (\ref{s}),
except one that corresponds 
to $i=s_{t+1},$ are zero. 
Moreover (\ref{xy}) implies that, due to the choice of $T,$
the element $u^{[T]}(1+s_{t+1}, m)\cdot D_w$ is a nonzero scalar,
while $u^{[L]}(1+s_{t+1}, m)\cdot D_w=0$ for any other super-word 
$u^{[L]}(1+s_{t+1}, m)$ that appears in the decomposition of $A_{s_{t+1}}$
with a nonzero coefficient.
Hence  $A_{s_{t+1}}\cdot D_w$ is a nonzero scalar $\mu .$ Finally, we have
\begin{equation}
\Psi ^{\hbox{\bf s}_t}(k,s_{t+1})=\mu ^{-1} B\cdot D_v   \in \hbox{\bf U}.
\label{pit2}
\end{equation}

\smallskip
2. Let us derivate (\ref{pit1}) by $x_{1+s_t}.$ By formulae 
(\ref{der8}) and (\ref{defdif}) we have
$$
\mu u[(2+s_t, m)] \cdot
\Psi ^{\hbox{\bf s}_{t-1}}(k,s_t)
$$
\begin{equation}
+\sum _{i=s_{t+1}}^{m-1}\mu_iA_iu[(2+s_{t},i)]\cdot \Psi ^{\hbox{\bf s}_{t-1}}(k,s_t)\in \, \hbox{\bf U},
\label{pit3}
\end{equation}
where $\mu =(1-q^{-1})^2,$  $\mu_i=(1-q^{-1})^2p(A_i,x_{1+s_t})$ with the only possible 
exception
$\mu _{s_{t+1}}=(1-q^{-1})p(A_i,x_{1+s_t})$ provided $s_{t+1}=1+s_t.$
Here we may apply (\ref{der8}) 
since the number $r$ related to {\bf S}$_t$ 
equals $t.$

Denote by $z$ the piecewise continuous word $u^{\hbox{\bf s}_{t-1}}(k,s_t).$
Let us apply $\cdot D_z$ to (\ref{pit3}). Formula (\ref{der9}) shows that
$\Psi ^{\hbox{\bf s}_{t-1}}(k,s_t)\cdot D_z$ is a nonzero scalar. Hence we get
$$
\mu u[2+s_t,m]+\sum _{i=s_{t+1}}^{m-1}\mu_iA_iu[2+s_t,i]\in \hbox{\bf U}.
$$
Let us apply $\cdot D_w$ with $w=u(2+s_t,s_{t+1})$ to this sum.
Formulae (\ref{defdif}) and (\ref{der4}) imply
\begin{equation}
(1-q^{-1})u[1+s_{t+1},m]+\beta _1A_{s_{t+1}}+
(1-q^{-1})\sum _{i>s_{t+1}}\beta _iA_iu[1+s_{t+1},i]\in \hbox{\bf U}.
\label{pit4}
\end{equation}
where $\beta _1=p(A_{s_{t+1}},x_{1+s_t}w)$ $=p_{vu},$ 
$\beta _i$ $=p(A_i, x_{1+s_t}w)$ $=p_{v_{i}\, u}$
with 
$$
u=u(1+s_t,s_{t+1}), \ v=u(1+s_{t+1},m),\ v_i=u(1+i,m).
$$
Let us multiply the element (\ref{pit4}) from the right by 
$\Psi ^{\hbox{\bf s}_t}(k,s_{t+1}) \in \, ${\bf U}, see (\ref{pit2}),  and subtract the result from
(\ref{pit1}) multiplied by $\beta _1.$ 
We get
$$
\beta _1\Psi ^{\hbox{\bf s}_t}(k,m)+(q^{-1}-1)u[1+s_{t+1},m]\cdot 
\Psi ^{\hbox{\bf s}_t}(k,s_{t+1})
$$
\begin{equation}
+\sum _{i=s_{t+2}}^{m-1}A_i\{\beta _1 \Psi ^{\hbox{\bf s}_t}(k,i)
+(q^{-1}-1)\beta_iu[1+s_{t+1},i]\cdot \Psi ^{\hbox{\bf s}_t}(k,s_{t+1}) \} \in \hbox{\bf  U}.
\label{pit6}
\end{equation}
By the recurrence formula (\ref{cbr3}) the first line of the above formula
equals $-\Psi ^{\hbox{\bf s}_{t+1}}(k,m),$ while the expression in the braces
equals $-\beta_i\Psi^{\hbox{\bf s}_{t+1}}(k,i).$
Thus we get  the required relation
\begin{equation}
\Psi ^{\hbox{\bf s}_{t+1}}(k,m)+\sum _{i=s_{t+2}}^{m-1} 
\beta_iA_i\Psi^{\hbox{\bf s}_{t+1}}(k,i)\in \hbox{\bf U}.
\label{pit7}
\end{equation}
Proposition \ref{phi} is proved.

\begin{corollary} 
If $q$ is not a root of $1,$ then $U_q^+(\frak{ sl}_{n+1})$ has
just a finite number of right coideal subalgebras that include the coradical.
If the multiplicative order of $q$ equals $t>2,$ then 
$u_q^+(\frak{ sl}_{n+1})$ has
just a finite number of homogeneous right coideal subalgebras 
that include the coradical.
\label{fin1}
\end{corollary}
\begin{proof}
This follows from Lemma \ref{odn} and Propositions \ref{pro}, \ref{phi}.
 Indeed, one has $n(n-1)/2$ options for possible value of an {\bf U}-root
(Definition \ref{root}). There exists $2^{n(n-1)/2}$ variants for sets of {\bf U}-roots. For any given 
root $\gamma =x_k+x_{k+1}+\ldots +x_m$ there exists not more than $2^{m-k}<2^n$
options for {\bf S} to define a PBW-generator $\Psi ^{\hbox{\bf s}}(k,m).$ Hence the total number of possible sets of PBW-generators  is less than $n^{(2^n)}\cdot 2^{n(n-1)/2}.$ \end{proof}

\section{Root sequence}

Our next goal is to show that the exact number of (homogeneous) right coideal 
subalgebras  in  $U_q^+(\frak{ sl}_{n+1})$ (in $u_q^+(\frak{ sl}_{n+1})$)
that contain {\bf k}$\, [G]$
equals $(n+1)!.$ In what follows for short we shall denote by $[k:m]$ the element 
$x_k+x_{k+1}+\ldots +x_m\in \Gamma ^+$ considered as an $U_q^+(\frak{ sl}_{n+1})$-root.

\begin{definition} \rm
Let $\gamma _k$ be a simple {\bf U}-root of the form $[k:m]$ with the maximal $m.$
Denote by $\theta _k$ the number $m-k+1,$ the length (weight) of $\gamma _k.$
If there are no simple {\bf U}-roots of the form $[k:m],$ we put $\theta _k=0.$
The sequence $r({\bf U})=(\theta_1, \theta_2, \ldots ,\theta_n)$
satisfies $0\leq \theta_k\leq n-k+1$ and it is uniquely defined by {\bf U}.
We shall call $r({\bf U})$ a {\it root sequence of } {\bf U}, or just an $r$-{\it sequence of} {\bf U}.
By $\tilde{\theta }_k$ we denote $k+\theta_k -1,$ the maximal value of $m$ for the simple
 {\bf U}-roots of the form $[k:m]$ with fixed $k.$
\label{tet}
\end{definition}

\begin{theorem} 
For each sequence 
$\theta=(\theta_1, \theta_2, \ldots ,\theta_n),$ such that $0\leq \theta_k\leq n-k+1,$
$1\leq k\leq n$ there exists one and only one $($homogeneous$)$ right coideal subalgebra
{\bf U}$\, \supseteq G$ of $U_q^+(\frak{ sl}_{n+1})$ $($respectively, of $u_q^+(\frak{ sl}_{n+1}))$
with $r({\bf U})=\theta .$ In what follows we shall denote this subalgebra by {\bf U}$_{\theta }.$
\label{teor}
\end{theorem}
The proof will result from the following lemmas.

\begin{lemma} 
If $[k:m]$ is an {\bf U}-root, then for each $r,$ $k\leq r< m$
either $[k:r]$ or  $[r+1:m]$ is an {\bf U}-root.
\label{su}
\end{lemma}
\begin{proof} By Proposition \ref{phi} we have $\Psi ^{\hbox{\bf s}}(k,m)\in \, ${\bf U}
for a suitable {\bf S}. If $r\in \, ${\bf S}, then Theorem \ref{26} and definition (\ref{pbse}) imply
$\Psi ^{\hbox{\bf s}}(k,r)\in \, ${\bf U}, hence $[k:r]$ is an {\bf U}-root.
If $r\notin \, ${\bf S}, then again Theorem \ref{26} and (\ref{pbse}) with $a=r+1$
imply $\Psi ^{\hbox{\bf s}}(r+1,m)\in \, ${\bf U}, hence $[r+1:m]$ is an {\bf U}-root.
\end{proof}

\begin{lemma} 
If $[k:m]$ is a simple {\bf U}-root, then 
there exists only one subset {\bf S} of the interval $[k,m-1],$
such that $\Psi ^{\hbox{\bf s}}(k,m)\in \, ${\bf U}.
Moreover the set {\bf S} is uniquely defined by the set of all {\bf U}-roots.
\label{su1}
\end{lemma}
\begin{proof} 
Let $\Psi ^{\hbox{\bf s}}(k,m)\in \, ${\bf U}.
By the definition of a simple root for each $r,$ $k\leq r<m$ either $[k:r]$ or $[r+1:m]$
is not an {\bf U}-root. Hence Lemma \ref{su}
and Theorem \ref{26} provide a criterion for $r$ to belong to {\bf S}:
$r\in \,${\bf S} if and only if $[k:r]$ is an {\bf U}-root.
\end{proof}

\begin{lemma} 
A $($homogeneous$)$ right coideal subalgebra {\bf U} is uniquely defined by the set of all its
simple roots.
\label{su2}
\end{lemma}
\begin{proof} Since obviously two subalgebras with the same PBW-basis coincide,
it suffices to find a PBW-basis of {\bf U} that depends only on a set of simple {\bf U}-roots.

We note first that the set of all {\bf U}-roots is uniquely defined by the set of simple {\bf U}-roots.
Indeed, if $[k:m]$ is an {\bf U}-root, then there exists a sequence $k=k_0<k_1<\ldots <k_l=m+1$
such that $[k_i:k_{i+1}-1],$ $0\leq i<l$ are simple {\bf U}-roots. Conversely, 
if there exists a sequence $k=k_0<k_1<\ldots <k_l=m+1$
such that $[k_i:k_{i+1}-1],$ $0\leq i<l$ are simple {\bf U}-roots then 
$f_i=\Psi ^{\hbox{\bf s}_i}(k_i,k_{i+1}-1)\in \, ${\bf U}, $0\leq i<l$ for suitable 
subsets {\bf S}$_i$ of the intervals $[k_i,k_{i+1}-2]$.
 By decomposition (\ref{cbrr}) we have 
\begin{equation}
\Psi ^{\hbox{\bf s}}(k,m)=[[[\ldots [f_{l-1},f_{l-2}],f_{l-3}],\ldots ],f_1], f_0]\in \, \hbox{\bf U},
\label{pet}
\end{equation}
where {\bf S}$\, =\cup_{i=0}^{l-1} S_i\cup \{ k_i-1 \, |\, 0< i<l\}.$
Thus, by definition, $[k:m]$ is an {\bf U}-root. Of course the decomposition of $[k:m]$
in a sum of simple {\bf U}-roots is not unique in general, however we may fix that decomposition for each non-simple {\bf U}-root from the very beginning. 

Now if $[k:m]$ is a simple {\bf U}-root, Lemma \ref{su1} 
shows that the element $\Psi ^{\hbox{\bf s}}(k,m)\in \, ${\bf U}
is uniquely defined by the set of simple {\bf U}-roots. We include this element
in the PBW-basis of {\bf U}. If $[k:m]$ is a non-simple {\bf U}-root with the fixed decomposition
in sum of simple {\bf U}-roots, then we include
in the PBW-basis the above defined element (\ref{pet}).  \end{proof}

\begin{lemma} 
If for $($homogeneous$)$ right coideal subalgebras
{\bf U}, {\bf U}$^{\prime }$ we have $r({\bf U})=r({\bf U}^{\prime }),$
 then {\bf U}$\, =\, ${\bf U}$^{\prime }.$
\label{su3}
\end{lemma}
\begin{proof} 
By Lemma \ref{su2} it suffices to show that the $r$-sequence uniquely defines the set of all simple roots. We use the downward induction on $k,$ the onset of a simple {\bf U}-root.
Suppose $k=n.$ Then the only possible root with the onset $n$ is
$\gamma =[n:n]=x_n,$ in which case
 $\gamma $  is a simple {\bf U}-root if and only if  $\theta _n=1.$

Let $k<n.$ 
By definition there  do not exist simple {\bf U}-roots of the form
$[k:m],$ $m>\tilde{\theta }_k$,  while $[k,\tilde{\theta }_k]$ is a simple {\bf U}-root. 

If $m<\tilde{\theta }_k,$
then $[m+1:\tilde{\theta }_k]$ is an {\bf U}-root if and only if it is a sum of simple {\bf U}-roots
starting with a number greater than $k.$ Hence by induction the $r$-sequence
defines all roots of the form $[m+1:\tilde{\theta }_k],$ $k\leq m<\tilde{\theta }_k.$ 

By Lemma \ref{su}
the weight $[k:m]$ is an {\bf U}-root if and only if $[m+1: \tilde{\theta }_k]$ is not an {\bf U}-root
(recall that $[k:\tilde{\theta }_k]$ is simple). Hence the $r$-sequence also defines  the set of
all   {\bf U}-roots of the form $[k:m],$ $m<\tilde{\theta }_k.$ An {\bf U}-root $[k:m],$
$m<\tilde{\theta }_k$
is simple if and only if there does not exist $r,$ $k\leq r<m$ such that both $[k:r]$ and 
$[r+1:m]$ are {\bf U}-roots. \end{proof}

Our next goal is a construction of a coideal subalgebra with a given root sequence
\begin{equation}
{\bf \theta }=(\theta_1, \theta_2, \ldots ,\theta_n),\hbox{ such that } 0\leq \theta_k\leq n-k+1,\ 
1\leq k\leq n.
\label{secu}
\end{equation}
We shall need the following technical definition. 
\begin{definition} \rm
By downward induction on $k$ we define subsets $R_k,$ $T_k,$ $1\leq k\leq n$
of the interval $[k,n]$ associated to a given sequence (\ref{secu})
as follows. If $\theta _n=0,$ we put $R_n=T_n=\emptyset .$ If 
$\theta _n=1,$ we put $R_n=T_n=\{ n \}.$ Suppose that 
$R_i,$ $T_i,$ $k<i\leq n$ are already defined. If $\theta _k=0,$ then we set 
$R_k=T_k=\emptyset .$ If $\theta _k\neq 0,$ then 
by definition $R_k$ contains 
$\tilde{\theta }_k=k+\theta _k-1$
and all $m$ satisfying the following three properties
\begin{equation}
\begin{matrix}\smallskip
a)\ k\leq m<\tilde{\theta }_k; \hfill \cr \smallskip
b)\ \tilde{\theta }_k \notin T_{m+1}; \hfill \cr
c) \  \forall r (k\leq r<m)\ \ m\in T_{r+1}\Longleftrightarrow \tilde{\theta }_k\in T_{r+1}. \hfill
\end{matrix}
\label{pet1}
\end{equation}
Respectively, 
\begin{equation}
T_k=R_k\cup \bigcup_{s\in R_k\setminus \{ n\} } T_{s+1}.
\label{pet2}
\end{equation}
\label{tski}
\end{definition}

\begin{example} \rm
Assume $n=3,$ ${\bf \theta }=(3,1,0).$

Since $\theta_3=0,$ by definition $R_3=T_3=\emptyset .$ 

Let $k=2.$ We have $\theta_2=1\neq 0,$ hence $\tilde{\theta }_2=2+\theta_2-1=2\in R_2.$
Certainly there are no points $m$ that satisfy $k=2\leq m<\tilde{\theta }_2=2,$
that is $R_2=\{ 2\} .$ Eq. (\ref{pet2}) yields 
$$
T_2=\{ 2\}\cup \bigcup_{s\in \{ 2\}\setminus \{ 3\} }T_{s+1}=\{ 2\}. 
$$

To find $R_1$ we take $k=1$ and consider $\tilde{\theta }_1=1+3-1=3.$ 
Obviously $\theta_1=3\neq 0;$ that is, $\tilde{\theta }_1=3\in R_1.$
There exist two points $m$ that satisfy $k=1\leq m<\tilde{\theta }_1=3,$
they are $m=1,$ $m=2.$ Both of them satisfy condition $b$) since
$\tilde{\theta }_1=3\notin T_{2+1}=\emptyset ,$ $\tilde{\theta }_1=3\notin T_{1+1}=\{ 2\} .$

Let us check condition $c$) for $m=1.$ The interval $1=k\leq r<m=1$ is empty.
Therefore the equivalence $c$) is true (elements from the empty set 
satisfy all conditions, even $r\neq r\, $). Thus $1\in R_1.$

It remains to check condition $c$) for $m=2.$ The interval $k=1\leq r<m=2$
has only one point $r=1.$ For this point we have $T_{r+1}=T_2$
and $m=2\in T_2=\{ 2\} .$ At the same time $\tilde{\theta }_1=3\notin T_{r+1}=T_2=\{ 2\} ;$
that is, condition $c$) for $m=2$ fails. Thus $R_1=\{ 1, 3\} .$

Finally by (\ref{tski}) we find
$$
T_1=\{ 1,3\}\cup \bigcup_{s\in R_1\setminus \{ 3\} }T_{s+1}=\{ 1,3\} \cup T_2=\{ 1,2,3\} .
$$

Thus for ${\bf \theta }=(3,1,0)$ we have $R_3=T_3=\emptyset ,$
$R_2=T_2=\{ 2\} ,$ $R_1=\{ 1,3\} ,$ $T_1=\{ 1,2,3\} .$
\label{rex1}
\end{example}

\begin{example} \rm
Assume $n=3,$ ${\bf \theta }=(2,1,1).$ 

Here $\theta_3\neq 0,$ hence by definition $R_3=T_3=\{ 3\} .$ 

Let $k=2.$ Since $\theta_2=1\neq 0,$ we have $\tilde{\theta }_2=2+1-1=2\in R_2.$
There are no points that satisfy $k=2\leq m<\tilde{\theta }_2=2;$ that is, $R_2=\{ 2\} .$
By Eq. (\ref{tski}) we have
$$
T_2=\{ 2\}\cup \bigcup_{s\in \{ 2\}\setminus \{ 3\} }T_{s+1}=\{ 2\}\cup \{ 3\} =\{ 2,3\}. 
$$

To find $R_1$ we take $k=1$ and consider $\tilde{\theta }_1=1+2-1=2.$ 
Since $\theta_1=3\neq 0,$ we have $\tilde{\theta }_1=2\in R_1.$
There exist only one point $m$ that satisfies $k=1\leq m<\tilde{\theta }_1=2,$
this is $m=1.$ For $m=1$ we have $\tilde{\theta }_1=2\in T_{1+1}=\{ 2,3\},$
hence condition $b$) fails. Thus $R_1=\{2\} .$ 
Finally by (\ref{tski}) we find
$$
T_1=\{ 2\}\cup \bigcup_{s\in \{ 2\} \setminus \{ 3\} }T_{s+1}=\{ 2\} \cup T_3=\{ 2,3\} .
$$

Therefore for ${\bf \theta }=(2,1,1)$ we have $R_3=T_3=\{ 3\} ,$
$R_2=\{ 2\} ,$ $T_2=\{ 2,3\} $ $R_1=\{ 2\} ,$ $T_1=\{ 2,3\} .$
\label{rex2}
\end{example}
\begin{lemma} For each sequence 
${\bf \theta }=(\theta_1, \theta_2, \ldots ,\theta_n),$ such that $0\leq \theta_k\leq n-k+1,$
$1\leq k\leq n$ there exists a homogeneous right coideal subalgebra
{\bf U} with $r(\hbox{\bf U})={\bf \theta }.$
\label{su4}
\end{lemma}
\begin{proof} 
Consider a set $W$ of all elements
$\Psi ^{\hbox{\bf s}}(k,m),$ $1\leq k\leq m\leq n$ with 
\begin{equation}
m\in R_k,\ \hbox{\bf S}\, =T_k.
\label{pet3}
\end{equation}
Denote by $U$ a subalgebra generated by values of $W$ in $U_q^+(\frak{ sl}_{n+1})$
or in $u_q^+(\frak{ sl}_{n+1}).$
 We shall check that {\bf U}$\, \stackrel{df}{=} U\# {\bf k}[G]$
 is a right coideal subalgebra with $r({\bf U})=\theta .$ To attain these ends we shall prove
some properties of the sets $R_k, T_k$ by downward induction on $k.$

\smallskip
\noindent
{\sc Claim} 1. $m\in T_k$ {\it if and only if there exists a sequence 
$k=k_0<k_1<\ldots <k_r=m+1,$  such that $k_{i+1}-1\in R_{k_i},$ $0\leq i<r.$}

\smallskip
\noindent
Assume first $m\in T_k.$
If $m\in R_k,$ we put $k_1=m+1,$ $r=1.$ If $m\notin R_k,$ then by definition
there exists $s\in R_k,$ such that $m\in T_{s+1}.$ We put $k_1=s+1.$
By the inductive supposition applied to $s+1$ there exists a sequence 
$s+1=k_1<k_2<\ldots <k_r=m+1$ such that $k_{i+1}-1\in R_{k_i},$ $1\leq r.$

Conversely. The inductive supposition applied to $k_1>k$ implies $m\in T_{k_1}.$
At the same time $k_1-1\in R_k,$ hence by definition $m\in T_k.$

\smallskip
\noindent
{\sc Claim} 2.  {\it If $s\in T_k,$ $m\in T_{s+1},$ then $m\in T_k.$} 

\smallskip
\noindent
By means of Claim 1 applied to $s$ we find a sequence 
$k_0=k<k_1<k_2<\ldots <k_r=s+1$ such that $k_{i+1}-1\in R_{k_i},$ $1\leq r.$
Again by Claim 1 we have $s\in T_{k_1}.$  
The inductive supposition applied to $k_1$ shows that $m\in T_{k_1}.$
Since $k_1-1\in R_k,$ it remains to apply definition (\ref{pet2}) with $s\leftarrow k_1-1.$

\smallskip
\noindent
{\sc Claim} 3.  {\it If $m\in T_k,$ then for all $s,$ $k\leq s<m$ either $s\in T_k,$ or $m\in T_{s+1}.$} 

\smallskip
\noindent
By Claim 1 there exists a sequence 
$k=k_0<k_1<\ldots <k_r=m+1,$  $k_{i+1}-1\in R_{k_i},$ $0\leq i<r.$
The same claim implies $m\in T_{k_1},$ provided that $r\geq 1.$

Since $k\leq s<m,$ there exists $i,$ $1\leq i\leq r,$ such that 
$k_i\leq s<k_{i+1}.$ If $i\geq 1,$ then the inductive supposition applied 
to $k_1$ implies that either $s\in T_{k_1}$ or $m\in T_{s+1}.$
If $m\in T_{s+1},$ we have got the required condition.
If $s\in T_{k_1},$ then definition (\ref{pet2}) with $s\leftarrow k_1-1$
yields $s\in T_k,$ which required.

Thus it remains to check the case $i=0;$ that is, $k\leq s<k_1,$ $k_1-1\in R_k.$
Claim 2 with $s\leftarrow k_1-1,$ $k\leftarrow s+1$ says that conditions 
$k_1-1\in T_{s+1}$ and $m\in T_{k_1}$ imply $m\in T_{s+1}.$
Hence it is sufficient to show that either $s\in T_k$ or $k_1-1\in T_{s+1}.$
If $s=k_1-1,$ then of course $s=k_1-1\in R_k\subseteq T_k.$
This allows us to replace $m$ with $k_1-1$  and suppose further that $m\in R_k,$ $i=0.$
In this case condition (\ref{pet1}$c$) with $r\leftarrow s$ is
$``m\in T_{s+1}\Longleftrightarrow \tilde{\theta }_k \in T_{s+1}".$
Therefore we have to consider only one case $m=\tilde{\theta }_k.$

Let us suppose that for some $s,$ $k\leq s<\tilde{\theta }_k$ we have 
$s\notin T_k$ and $\tilde{\theta }_k\notin T_{s+1}.$ By induction on $s,$
in addition to the downward induction on $k,$ we shall show
that these conditions imply $s\in R_k,$ which certainly contradicts to 
$s\notin T_k,$ see definition (\ref{pet2}).

Definition (\ref{pet1}) with $m=k$ shows that $k\in R_k$ if and only if 
$\tilde{\theta }_k\notin T_{k+1}.$ Since in our case $\tilde{\theta }_k\notin T_{s+1},$
we have $s\in R_k,$ provided that $s=k.$

Let $s>k.$ Conditions (\ref{pet1}$a$) and (\ref{pet1}$b$) are valid for $m\leftarrow s.$
Suppose that (\ref{pet1}$c$) fails. In this case we may find a number $t,$ $k\leq t<s,$
such that $\neg (s\in T_{t+1}\Longleftrightarrow \tilde{\theta}_k \in T_{t+1}).$

If $s\in T_{t+1}$ but $\tilde{\theta_k }\notin T_{t+1},$ then by the inductive supposition
(induction on s) either $t\in R_k$ or $\tilde{\theta}_k \in T_{t+1};$ that is,
$t\in R_k.$ Definition (\ref{pet2}) implies  $s\in T_k$ --- a contradiction. 

If $\tilde{\theta}_k \in T_{t+1}$ but $s\notin T_{t+1},$ then the inductive supposition
of the downward induction on $k$
with $k\leftarrow t+1$ shows that either $s\in T_{t+1}$ or $\tilde{\theta}_k \in T_{s+1};$
that is, $\tilde{\theta}_k \in T_{s+1}$ --- again a contradiction.

Thus $s$ satisfies all conditions (\ref{pet1}a --- \ref{pet1}$c$), hence $s\in R_k.$

\smallskip
\noindent
{\sc Claim 4}. {\it If $k\leq m<\tilde{\theta }_k,$ then $m\in T_k$ if and only if 
$\tilde{\theta }_k\notin T_{m+1}$}.

\smallskip
\noindent
According to Claim 3 one of the conditions $m\in T_k$ or 
$\tilde{\theta }_k\in T_{m+1}$ always holds. If both of them are valid then
due to Claim 1 we find a sequence 
$k=k_0<k_1<\ldots <k_r=m+1,$  such that $k_{i+1}-1\in R_{k_i},$ $0\leq i<r.$
Due to (\ref{pet1}$b$), we have $m\notin R_k,$ hence $r>1.$ Again by the first
claim  we have $m\in T_{k_1}.$
Since $k_1-1$ belongs to $R_k,$ it satisfies condition (\ref{pet1}$b$), 
$\tilde{\theta }_k\notin T_{k_1}.$ However Claim 2 shows that the conditions 
 $m\in T_{k_1},$ $\tilde{\theta }_k\in T_{m+1} $
 imply $\tilde{\theta }_k\in T_{k_1}.$
A contradiction, that proves the claim.

\smallskip
\noindent
{\sc Claim 5}. {\it The subalgebra $U^{\prime }$ generated by $\Psi ^{\hbox{\bf s}}(k,m),$
$1\leq k\leq m\leq n,$ $m\in T_k,$ {\bf S}$\, =T_k$ is a differential subalgebra.}

\smallskip
\noindent
It suffices to show that all partial derivatives of $\Psi ^{\hbox{\bf s}}(k,m)$ belong to $U.$
By Theorem \ref{26} we have to check that $\Psi ^{T_k}(a,b)\in U$ provided that
$b\in T_k,$ $a-1\notin T_k,$ $k\leq a\leq b\leq m.$ By definition
$\Psi ^{T_a}(a,b)\in U$ since due to the third claim $b\in T_a.$ If 
\begin{equation}
T_k\cap [a, b-1]=T_a\cap [a,b-1],
\label{tor}
\end{equation}
then we have nothing to prove. In general, however, just the inclusion
$T_k\cap [a, b-1]\subseteq T_a\cap [a,b-1]$ holds: if $t\in T_k,$ $a\leq t,$
then Claim 3 with $s\leftarrow a-1$ says $t\in T_a$ (since $a-1\notin T_k$).

We shall prove $\Psi ^{T_k}(a,b)\in U$ by induction on $b-a.$ If $b=a,$
then (\ref{tor}) certainly holds.

Let us choose the minimal $s\in T_a,$ $s\notin T_k.$ Then $T_k\cap [a, s-1]=T_a\cap [a,s-1].$
Hence $\Psi ^{T_k}(a,s)=\Psi ^{T_a}(a,s)\in U.$
By the inductive supposition applied to the interval $[s+1,b]$
we get $\Psi ^{T_k}(s+1,b)\in U.$
By decomposition (\ref{cbrr}) we have 
$$
\Psi ^{{T_k}\cup \{ s \}}(a,b)=[\Psi ^{T_k}(s+1,b), \Psi ^{T_k}(a,s)]\in U.
$$
At the same time (\ref{cby}) implies
$$
\Psi ^{{T_k}\cup \{ s \}}(a,b)-(1-q^{-1})\Psi ^{T_k}(s+1,b)\cdot \Psi ^{T_k}(a,s)=
-p_{vu}\Psi ^{T_k}(a,b).
$$
Therefore $\Psi ^{T_k}(a,b)\in U,$ which is required.

\smallskip
\noindent
{\sc Claim 6}. {\it  {\bf U}$\, =U\#{\bf k}[G]$ is a right coideal subalgebra.}

\smallskip
\noindent
Since $U$ is homogeneous in each variable, we have  
$g^{-1}Ug\subseteq U,$ $g\in G.$
It remains to apply Lemma \ref{qsim}.

\smallskip
\noindent
{\sc Claim 7}. {\it  The set of all \, {\bf U}-roots is $\{ [k:m]\, |\, m\in T_k\}.$ In particular
$\{ \Psi ^{T_k}(k,m)\, |$ $m\in T_k\} $ is a set of PBW-generators of {\bf U} over {\bf k}$[G].$}

\smallskip
\noindent
If $\gamma =[a:b]$ is an {\bf U}-root, then, by definition, in {\bf U} there exists
a homogeneous element (\ref{vad22}) of degree $\gamma .$ Since by definition 
$U$ is generated by  $\{ \Psi ^{T_k}(k,m)\, |\, m\in T_k\} ,$ the degree $\gamma $
is a sum of degrees of the generators: 
$\gamma =[k_1:k_2-1]+[k_2:k_3-1]+\ldots +[k_{r-1}:k_r-1],$ $k_{i+1}-1\in T_{k_i},$ $1\leq i<r.$
The multiple use of Claim 2 yields $k_r-1\in T_{k_1},$ which is required since $a=k_1,$
$b=k_r-1.$

\smallskip
\noindent
{\sc Claim 8}. {\it  The set of all simple {\bf U}-roots is $\{ [k:m]\, |\, m\in R_k\}.$ In particular
$r({\bf U})=\theta ,$ and {\bf U} is generated as an algebra by 
$\Psi ^{T_k}(k,m),\, $  $m\in R_k ,$ $1\leq k\leq n,$ and {\bf k}$[G].$}

\smallskip
\noindent
If $\gamma =[k:m]$ is a simple {\bf U}-root, then due to the above claim $m\in T_k.$
Hence, according to Claim 1, we may find a sequence 
$k=k_0<k_1<\ldots <k_r=m+1,$  such that $k_{i+1}-1\in R_{k_i},$ $0\leq i<r.$
In this case $\gamma =[k:k_1-1]+[k_1:k_2-1]+\ldots +[k_{r-1}:m]$ is a sum of {\bf U}-roots.
Since $\gamma $ is simple this is possible only if $r=1,$ and $m=k_1-1\in R_k.$

Conversely. Let $m\in R_k.$ Then by Claim 7 and definition (\ref{pet2}) $[k:m]$ is an {\bf U}-root.
If it is not simple, then it is a sum of two roots $[k:m]=[k:s]+[s+1:m],$ $s\in T_k,$ $m\in T_{s+1}.$
Then Claim 4 implies $\tilde{\theta }_k\notin T_{s+1},$ hence condition (\ref{pet1}$c$) fails
for $r\leftarrow s.$

In Lemma \ref{su2} we have seen that {\bf U} is uniquely defined by the simple 
{\bf U}-roots. In particular formula (\ref{pet}) provides a representation of PBW-generators
in terms of $\Psi ^{T_k}(k,m),$  $m\in R_k.$
Lemma \ref{su4} and Theorem \ref{teor} are completely proved. 
\end{proof}

\section{Homogeneous right coideal subalgebras}

In this section we consider right coideal subalgebras
in  $U^{+}_q(\mathfrak{sl}_{n+1}),$ (respectively in $u^{+}_q(\mathfrak{sl}_{n+1})$
if $q^{t}=1, t>2$) that do not contain the coradical. First of all we note that for every submonoid 
$\Omega \subseteq G$ the set of all linear combinations {\bf k}$\, [\Omega]$
is a right coideal subalgebra.
Conversely if $U_0\subseteq \,${\bf k}$\, [G]$ is a right coideal subalgebra then
$U_0=\, ${\bf k}$\, [\Omega]$ for $\Omega =U_0\cap G$
since  $a=\sum_i \alpha_i g_i\in U_0$ implies $\Delta (a)=
\sum_i \alpha_i g_i\otimes g_i\in U_0\otimes \, ${\bf k}$\, [G];$ that is, $\alpha_i g_i\in U_0.$
\begin{definition} \rm
For a sequence $\theta=(\theta_1, \theta_2, \ldots ,\theta_n),$ such that $0\leq \theta_k\leq n-k+1,$
$1\leq k\leq n$ we denote by {\bf U}$^1_{\theta }$ a subalgebra with 1 generated by 
$g^{-1}\Psi ^{\hbox{\bf s}}(k,m),$ where $g=g_kg_{k+1}\ldots g_m,$ and $\Psi ^{\hbox{\bf s}}(k,m)$
runs through the set of PBW-generators of {\bf U}$_{\theta },$ see Theorem \ref{teor}
and Claim 7, Section 5.
\label{1te}
\end{definition}
\begin{lemma} The subalgebra
{\bf U}$^1_{\theta }$ is a homogeneous right coideal, and 
{\bf U}$^1_{\theta }\cap G=\{ 1\} .$ 
\label{1su}
\end{lemma}
\begin{proof}
The subalgebra
{\bf U}$^1_{\theta }$ is homogeneous since it is generated by homogeneous elements.
Its zero homogeneous  component equals {\bf k} since among the generators just one,
the unity, has zero degree. 

Denote by $A_{\theta }$ a {\bf k}-subalgebra generated by 
the PBW-generators $\Psi ^{\hbox{\bf s}}(k,m)$ of {\bf U}$_{\theta }.$
The algebra {\bf U}$^1_{\theta }$ is spanned by all elements of the form 
$g_a^{-1}a,$ $a\in A_{\theta }.$ Since {\bf U}$_{\theta }$ is a right coideal, 
for any homogeneous $a\in A_{\theta }$ we have
$\Delta (a)=\sum g(a^{(2)})a^{(1)}\otimes a^{(2)}$ where $a^{(1)}\in A_{\theta },$
$g_a=g(a^{(1)})g(a^{(2)}).$ Therefore 
$\Delta (g_a^{-1}a)=\sum g(a^{(1)})^{-1}a^{(1)}\otimes g_a^{-1}a^{(2)}$ with 
$g(a^{(1)})^{-1}a^{(1)}\in \, ${\bf U}$_{\theta }^1.$
\end{proof}
\begin{theorem} If $U$ is a homogeneous right coideal subalgebra   
of  $U^{+}_q(\mathfrak{sl}_{n+1})\,  ($resp. of $u^{+}_q(\mathfrak{sl}_{n+1}))$
such that $\Omega \stackrel{df}{=} U\cap G$ is a group,
then $U=\, ${\bf U}$_{\theta }^{1}\, ${\bf k }$[\Omega ]$
for a suitable $ \theta .$
\label{orc}
\end{theorem}
\begin{proof}
Let $u=\sum h_ia_i\in U$ be a homogeneous element of degree $\gamma \in \Gamma ^{+}$
with different $h_i\in G,$ and $a_i\in A,$ where by $A$ we denote the {\bf k}-subalgebra
generated by $x_i,$ $1\leq i\leq n.$ Denote by ${\pi }_{\gamma }$ the natural
projection on the homogeneous component of degree $\gamma .$
Respectively $\pi _g,$ $g\in G$ is a projection on the subspace {\bf k}$\, g.$
We have $\Delta (u)\cdot (\pi _{\gamma }\otimes \pi _{h_i})=h_ia_i\otimes h_i.$
Thus $h_ia_i\in U.$

By Theorem \ref{teor} we have {\bf k}$\, [G]U=\,${\bf U}$_{\theta }$ for a suitable $\theta .$
If $u=ha\in U,$ $h\in G,$ $a\in A,$ then $\Delta (u)\cdot (\pi _{hg_a}\otimes \pi _{\gamma })=
hg_a\otimes ha.$ Therefore $hg_a\in U\cap G=\Omega ;$ that is, $u=\omega g_a^{-1}a,$
$\omega \in \Omega .$ Since $\Omega $ is a subgroup we get $g_a^{-1}a\in U.$
It remains to note that all elements $g_a^{-1}a,$ such that $ha\in U$ span the algebra 
{\bf U}$_{\theta }^1.$
\end{proof}
If $U\cap G$ is not a group then $U$ may have a more complicated structure.
\begin{example} \rm 
Let $\Omega $ be a submonoid of $G.$  Denote by $\overline{\Omega }$ an arbitrary family of sets 
$\{ \Omega _{\gamma }, \gamma \in \Gamma ^{+}\} $ that satisfies 
the following conditions
$$
\Omega _0=\Omega, \ \  \ \ 
\Omega _{\gamma }\cdot \Omega _{\gamma ^{\prime }}\subseteq 
\Omega _{\gamma +\gamma ^{\prime }}\subseteq
 \Omega _{\gamma }\cap \Omega _{\gamma ^{\prime }}.
$$
In this case the linear space {\bf U}$_{\theta }^{\overline{\Omega }}$ spanned
by the elements $\omega _{\gamma }a,$ 
$\omega _{\gamma }\in \Omega _{\gamma },$ $a\in \, ${\bf U}$_{\theta }^1,$
deg$\, (a)=\gamma $ is a right coideal subalgebra such that 
{\bf U}$_{\theta }^{\overline{\Omega }}\cap G=\Omega .$ The $\gamma $-homogeneous component
of this algebra equals $\Omega _{\gamma }\, (${\bf U}$_{\theta }^1)_{\gamma }.$
 Hence different $\overline{\Omega }$ define different homogeneous  right coideal subalgebras.
\label{oxc1}
\end{example}
Finally we point out a simplest one-parameter family of inhomogeneous 
right coideal subalgebras that have trivial intersection  with the coradical.
\begin{example} \rm 
Let $a=g_1^{-1}(x_1+\alpha ),$ $\alpha \in \,${\bf k}.
We have 
$$
\Delta (a)=g_1^{-1}x_1\otimes g_1^{-1}+ 1\otimes g_1^{-1}x_1+\alpha g_1^{-1}\otimes g_1^{-1}
=a\otimes g_1^{-1}+1\otimes g_1^{-1}x_1.
$$
Therefore the two-dimensional space spanned by $a$ and $1$ is a right coideal.  
Thus the algebra {\bf k}$\, [a]$ with 1 generated by $a$ is a right coideal subalgebra,
in which case  {\bf k}$\, [a]\cap G=\{ 1\} .$
\label{oxc2}
\end{example}

\section{K\'eb\'e construction and ${\rm ad}_r$-invariant subalgebras}

In this section we characterize ad$_r$-invariant right coideal subalgebras that have trivial intersection 
with the coradical in terms of   K\'eb\'e's construction \cite{Keb, Keb1}. Recall that the right 
adjoint action of a Hopf algebra $H$ on itself is defined by the formula
$$
({\rm ad}_ra)b=\sum \sigma (a^{(1)})ba^{(2)},
$$
where $\sigma $ is the antipode. The map $a\rightarrow \, $ad$_ra$ is a homomorphism 
of algebras ad$_r:H\rightarrow \hbox{End}H.$ In particular a subspace is invariant under 
the action of all operators ad$_rH$ if and only if it is invariant under the actions of ad$_rh_i$
for some set of generators $\{ h_i\}.$ For $H=U_q^+(\mathfrak{sl}_{n+1})$ or for
$H=u_q^+(\mathfrak{sl}_{n+1})$ we have
$$
(\hbox{ad}_rg)b=g^{-1}bg, \ \ g\in G; \ \ \ (\hbox{ad}_rx_i)b=g_i^{-1}(bx_i-x_ib).
$$
The latter equality would be more familiar if we take $b=g_a^{-1}a$ with
$a\in A\stackrel{df}{=}\hbox{\bf k}\langle x_1,\ldots ,x_n\rangle :$
\begin{equation}
(\hbox{ad}_rx_i)(g_a^{-1}a)=g_i^{-1}(g_a^{-1}ax_i-x_ig_a^{-1}a)=-g_i^{-1}g_a^{-1}[x_i,a].
\label{keb1}
\end{equation}
In particular the subalgebra $H^1$ generated by $g_a^{-1}a,$ $a\in A$ 
(in our terms this is {\bf U}$^1_{\theta }$ for $\theta =(1,1,\ldots, 1)$) is ad$_r$-invariant.

The following construction of ad$_r$-invariant right coideal subalgebras appeared in 
\cite{Keb, Keb1},  see also \cite[Section 6]{Let}. Let $\pi $ be a subset of $[1,k].$
Denote by $K(\pi )$ a subalgebra generated by elements of the form
$$
\hbox{ad}_r(x_{i_1}x_{i_2}\ldots x_{i_k})\, g_j^{-1}x_j, \ \ \ j\in \pi,\  \ i_r\notin \pi, 1\leq r\leq k.
$$
The algebra $K(\pi )$ is ad$_r$-invariant right coideal (see, \cite[Lemma 1.2]{Let} up to a left-right symmetry). This is homogeneous, and $K(\pi )\cap G=\{ 1\}$
since due to (\ref{keb1}) the inclusion  $K(\pi )\subseteq H^1$ is valid. Thus by Theorem \ref{orc}  
we have $K(\pi )={\bf U}^1_{\theta }$ for a suitable $\theta .$
 
\begin{theorem}
The following conditions on $U=\,${\bf U}$^1_{\theta }$ are equivalent

i. $U$ is {\rm ad}$_r$-invariant.

ii. The sets $T_k,$ see Definition $\ref{tski},$ have the form $T_k=[j(k),n],$ where 
$$
j(k)\stackrel{df}{=}\, {\rm min}\, \{ j\, |\, k\leq j,\ \ j\in T_j\} .
$$

iii. $U=K(\pi )$ for a suitable $\pi \subseteq [1,n].$ 
\label{Korc}
\end{theorem}
\begin{proof}
{\it  iii}$\, \Rightarrow \,${\it  i}. We have mentioned above.
 
{\it  i}$\, \Rightarrow \,${\it  ii}. By Claim 7 we have $m\in T_k$ if and only if 
$\Psi ^{T_k}(k,m)\in \hbox{\bf U}_{\theta }.$ In particular $j\in T_j$ if and only if 
$x_j\in  \hbox{\bf U}_{\theta },$ or, equivalently, $g_j^{-1}x_j\in U.$
If $j\in T_j$ and $k\leq j,$ then by (\ref{keb1}) we have
$$
{\rm ad}_r(x_{j-1}x_{j-2}\ldots x_k)\ g_j^{-1}x_j=g_{u[k,j]}^{-1}u[k,j]\in U.
$$
Hence, by definition $[k:j]$ is an {\bf U}$_{\theta }$-root; that is, 
according to Claim 7, we get  $j\in T_k.$ Moreover, if $i>j$ then by
(\ref{cin}) we have 
$$
{\rm ad}_r(x_{j+1}x_{j+2}\ldots x_i)\ g_{u[k,j]}^{-1}u[k,j]=g^{-1}[x_i,[x_{i-1},\ldots [x_{j+1}, u[k,j]]\ldots]]
$$
$$
\sim g^{-1}\Psi ^{\{ j+1,j+2,\ldots ,i\}}(k,i),
$$
where $g=g_kg_{k+1}\ldots g_i.$ In particular $[k:i]$ is an {\bf U}$_{\theta }$-root; that is, 
again according to Claim 7, we get  $i\in T_k.$ This proves $[j(k),n]\subseteq T_k.$

If $m$ is the smallest element from $T_k$ then $u[k,m]=\Psi ^{T_k}(k,m)\in \hbox{\bf U}_{\theta },$
hence by multiple use of (\ref{der1}) we get $x_m\in \hbox{\bf U}_{\theta };$ that is, $m\in T_m,$
and $m=j(k).$

{\it  ii}$\, \Rightarrow \,${\it  iii}. Let $\pi =\{ j\, |\, j\in T_j\} .$  For all $k,m,j$ such that 
$m\in T_k,$ $j=j(k)$ we have 
$$
{\rm ad}_r(x_{j-1}x_{j-2}\ldots x_kx_{j+1}x_{j+2}\ldots x_m)\ g_j^{-1}x_j
$$
$$=g^{-1}[x_m,[x_{m-1},\ldots [x_{j+1}, u[k,j]]\ldots]]=g^{-1}\Psi ^{T_k}(k,m),
$$
where $g=g_mg_{m-1}\ldots g_k.$ Since $K(\pi )$ is ad$_r$-invariant, we get
$g^{-1}\Psi ^{T_k}(k,m)\in K(\pi ).$ Now Definition \ref{1te} implies $U\subseteq K(\pi ).$

Since $g_j^{-1}x_j\in U,$ to check $K(\pi )\subseteq U$  it remains to show that
ad$_r(x_i)U\subseteq U$ for $i\notin \pi .$ By (\ref{br1}) and (\ref{keb1}) it suffices to prove that 
$[x_i , \Psi ^{T_k}(k,m)]\in {\bf U}_{\theta }$ for $i\notin \pi ,$ $m\in T_k.$

Let $i=k-1.$ Since $k-1=i\notin \pi ,$ we have $k-1\notin T_{k-1},$ and hence 
$j(k-1)=j(k).$ Eq. (\ref{cbry})
implies $[x_{k-1}, \Psi ^{T_k}(k,m)]\sim \Psi ^{T_{k-1}}(k-1,m)\in {\bf U}_{\theta },$
where $m\in T_{k-1}$ follows from $T_{k-1}=[j(k),n]=T_k.$

If $i<k-1,$ then $[x_i , \Psi ^{T_k}(k,m)]=0$ since $x_i$ and $\Psi ^{T_k}(k,m)$
are separated.

Let $i=m+1.$ Since $m\in T_k=[j(k),n],$ we have $m\geq j(k).$ 
Therefore $m+1>j(k),$ and $m+1\in T_k.$ Now formula (\ref{cin}) yields  
$$[x_{m+1}, \Psi ^{T_k}(k,m)]\sim \Psi ^{T_k}(k,m+1)\in {\bf U}_{\theta }.$$

If $i>m+1,$ then $[x_i , \Psi ^{T_k}(k,m)]=0$ since $x_i$ and $\Psi ^{T_k}(k,m)$
are separated.

We shall show by induction on $m-k$ that in all remaining cases $[x_i , \Psi ^{T_k}(k,m)]=0.$
More precisely, we prove $[x_i,  \Psi ^{\hbox{\bf s}}(k,m)]=0$ provided that 
{\bf S} has the form $[j,n],$ and $k\leq i\leq m,$ $i\neq \hbox{\rm min}\, \{ j,m \}.$ 

If $m-k=1,$ then for $j\geq m$ we have just one option $i=k.$ The required relation 
takes the form $[x_k,[x_k,x_{k+1}]]=0$ which is one of the defining relations (\ref{rela}).
For $j=k$ we also have just one option $i=m=k+1.$ The required relation is 
$[x_m,[x_m,x_{m-1}]]=0.$ This relation is valid in $U_q^{+}(\mathfrak{sl}_{n+1})$
since (\ref{a1rel}) imply $[x_m,[x_m,x_{m-1}]]\sim [[x_{m-1},x_m],x_m],$
see, for example, \cite[Corollary 4.10]{Kh4}.

If $m-k>2$ then either $i<m-1$ or $i>k+1.$ In the former case we have $[x_i,x_m]=0,$
and by the inductive supposition $[x_i , \Psi ^{\hbox{\bf s}}(k,m-1)]=0.$ Hence representation
(\ref{cin}) implies the required equality. In the latter case we have $[x_i,x_k]=0,$
and by the inductive supposition $[x_i , \Psi ^{\hbox{\bf s}}(k+1,m)]=0.$ In this case representation
(\ref{cin1}) implies the required equality.

Finally, suppose that $m-k=2.$ To simplify the notations we put $k=1,$ $m=3.$

If $j\geq 3$ then  $\Psi ^{\hbox{\bf s}}(1,3)=[[x_1,x_2],x_3],$  
and we have two options $i=1,$ $i=2.$ If $i=1,$ we have to show
$[x_1,[[x_1,x_2],x_3]]=0.$ This relation is valid since $x_1$ (skew)commutes both with 
$[x_1,x_2]$ and $x_3$ (but not vise versa: $[x_1,x_2]$ does not (skew)commute with $x_1$
since $[[x_1,x_2],x_1]\neq 0!$)  Let  $i=2.$ We may apply (\ref{bri}) since 
$p_{21}p_{22}p_{23}\cdot p_{12}p_{22}p_{32}=1.$ Thus 
by (\ref{bri}) and (\ref{jak3}) we have  
$$
[x_2,[[x_1,x_2],x_3]]\sim [[[x_1,x_2],x_3],x_2]=[[x_1,[x_2,x_3]],x_2].
$$
The word $x_1x_2x_3x_2$ is standard, and the standard alignment of brackets is precisely
$[[x_1,[x_2,x_3]],x_2].$ Hence by the third statement of Theorem \ref{strA} this is zero in 
$U_q^{+}(\mathfrak{sl}_{n+1}).$

If $j=2,$ then $\Psi ^{\hbox{\bf s}}(1,3)=[x_3, [x_1,x_2]],$ and 
we have two options $i=1,$ $i=3.$ If $i=1$ then
$[x_1,[x_3,[x_1,x_2]]]=0$ since $x_1$ (skew)commutes both with 
$[x_1,x_2]$ and $x_3.$ Let $i=3.$ By (\ref{jak4}) we have $[x_3,[x_1,x_2]]\sim [x_1,[x_3,x_2]]].$
Since $x_3$ (skew)commutes both with $x_1$ and $[x_3,x_2],$ we get 
$[x_3,[x_1,[x_3,x_2]]]=0.$

If $j=1$ then $\Psi ^{\hbox{\bf s}}(1,3)=[[x_3,x_2],x_1],$
 and we have two options $i=2,$ $i=3.$ If $i=3$ then
$[x_3,[[x_3,x_2],x_1]]=0$ since $x_3$ (skew)commutes both with 
$[x_3,x_2]$ and $x_1.$ Let $i=2.$ 
We may use (\ref{bri}) since 
$p_{23}p_{22}p_{21}\cdot p_{32}p_{22}p_{21}=1.$
Thus by (\ref{bri}) and (\ref{jak3}) we have 
$$
[x_2,[[x_3,x_2],x_1]]\sim [[[x_3,x_2],x_1],x_2]=[[x_3,[x_2,x_1]],x_2].
$$
This element in new variables $y_1=x_3,$ $y_2=x_2,$ $y_3=x_3$ takes up the form
$[[y_1,[y_2,y_3]],y_2].$ By the third statement of Theorem \ref{strA}
this is zero in $U_q^{+}(\mathfrak{sl}_{n+1}).$
\end{proof}

\section{Examples}
In this section by means of Theorem \ref{teor} we provide some examples of right coideal subalgebras in $U_q^+(sl_n)$ or $u_q^+(sl_n)$ with their main characteristics:
PBW-generators, the root sequence $r({\bf U}),$ 
the sets $T_i,$ $R_i,$ right coideal subalgebra generators, and maximal Hopf subalgebras.
We start with a  characterization of 
$2^n$ ``trivial" examples --- Hopf subalgebras.

\begin{proposition}  
A right coideal subalgebra {\bf U}$={\bf U}_{\theta }$
is a Hopf subalgebra if and only if 
for every $k,$ $1\leq k\leq n$ either $\theta _k=0$ or $\theta _k=1.$
An algebra {\bf U}$_{\theta },$ with $\theta _i\leq 1$
is generated over {\bf k}$[G]$ by all $x_k$ with $\theta _k=1.$   
\label{hop}
\end{proposition}
\begin{proof} 
If $\theta _k\leq 1,$ $1\leq k\leq n,$ then Definition \ref{pet1} shows that $R_k=\{ k\} $
provided that $\theta _k=1$ and $R_k=\emptyset $ otherwise.
Hence by Claim 8  the algebra {\bf U} is generated over {\bf k}$[G]$ 
by all $x_k$ with $\theta _k=1.$   In particular {\bf U} is a Hopf subalgebra of
$U_q^+(sl_n).$ 

Conversely, let {\bf U} be a Hopf subalgebra. According to Claim 8 the algebra {\bf U} is generated 
over {\bf k}$[G]$ by the elements $a$ of the form $\Psi ^{T_k}(k,m)$  
with $[k:m]$ being  the simple {\bf U}-roots.
We have $\Delta (a)=\sum a^{(1)}\otimes a^{(2)}$ with 
$a^{(1)}, a^{(2)}\in \, ${\bf U}. Since $[k:m]$ $=D(a)$ $=D(a^{(1)})+D(a^{(2)})$ and $[k:m]$ is simple,
we have either $D(a^{(1)})=0,$ or $D(a^{(2)})=0.$ Thus $a$ is a skew primitive element;
that is, $a=x_k$ is the only option for $a$ (see the second statement of Theorem \ref{strA}  for $q^t\neq1,$ and comments after that theorem for $q^t=1).$
In particular all simple roots are of length 1, while Definition \ref{tet} implies 
$\theta _k\leq 1.$  \end{proof}

Now we consider three special cases.

\begin{example} \rm 
Consider  the root sequence with the maximal possible components,
$r({\bf U})=(n,n-1,n-2,\ldots ,2,1).$
 In this case by definition $T_n=R_n=\{ n\}.$ For $k<n$ we have
$\tilde{\theta }_k=k+\theta _k-1=n\in R_k.$ Moreover by downward induction on $k$
one may prove that $T_k=R_k=\{ n\}.$ Indeed, the inductive supposition implies that
condition (\ref{pet1}$b$) fails for all $m,$ $k\leq m<n.$ This yields $R_k=\{ n\}.$
Inductive definition (\ref{pet2}) implies  $T_k=\{ n\}$ as well.
Claim 7 provides a PBW-basis: 
$$
\{  \Psi ^{\{ n\} }(k,n)\, |\, 1\leq k\leq n\}=\{ [x_kx_{k-1}\ldots x_{n-1}x_n]\, |\, 1\leq k\leq n\} .
$$
 Due to formula (\ref{der1}) the subalgebra  {\bf U}
as a right coideal subalgebra
is generated over {\bf k}$[G]$ by a single long skew commutator  
$[x_1x_2\ldots x_n]$ that has the whitest diagram (\ref{gra}). The maximal
Hopf subalgebra of {\bf U} is {\bf k}$[G]\langle x_n\rangle .$
\label{sex1}
\end{example}

\smallskip
\begin{example} \rm
Let $r({\bf U})=(n,0,0,\ldots ,0).$ In this case according to definitions (\ref{pet1}) and (\ref{pet2}) we have $T_i=R_i=\emptyset $ if $1<i\leq n,$ and $T_1=R_1=\{ 1,2,\ldots n-1,n\}.$ By means of Claim 7
we obtain a PBW-basis: 
$$
\{  \Psi ^{T_1}(1,m)\, |\, 1\leq m\leq n\}=\{ [x_mx_{m-1}\ldots x_2x_1]\, |\, 1\leq m\leq n\} .
$$ 
As a right coideal subalgebra {\bf U}
is generated by a single long skew commutator  $[x_nx_{n-1}\ldots x_1]$
that has the most black diagram (\ref{gra}). The maximal
Hopf subalgebra of {\bf U} is {\bf k}$[G]\langle x_1\rangle .$
\label{sex2}
\end{example}

\smallskip
\begin{example} \rm
More generally consider a right coideal subalgebra {\bf U}$^{\hbox{\bf s}}(k,m)$
generated by $\Psi ^{\hbox{\bf s}}(k,m).$ We claim that the root sequence for
this algebra is defined as follows
\begin{equation}
\theta _i=\left\{ \begin{matrix}m-k+1, \hfill &\ \hbox{if } i=k; \hfill \cr
m-i+1,\hfill &\ \hbox{if } i-1\notin {\bf S},\ k<i\leq m; \hfill \cr
0,\hfill &\ \hbox{otherwise;} \hfill \end{matrix} \right. 
\label{roo}
\end{equation}
that is, $\theta _i$ takes the maximal value if $i-1$ is a white point on the diagram
(\ref{gra}), and $\theta_i=0$ if either $i-1$ is a black point or $i-1$ 
is not displayed on the diagram at all.

Let {\bf U}$_{\theta }$ be a right coideal subalgebra  defined by the sequence (\ref{roo}) 
in Lemma \ref{su4}. By downward induction of Lemma \ref{su4} it is easy to see that 
\begin{equation}
R_i=T_i=\left\{ \begin{matrix}({\bf S}\, \cap [i,m])\cup \{ m\} ,\hfill &\ \hbox{if } 
i=k \hbox{ or } i-1\notin {\bf S};\hfill \cr
\emptyset ,\hfill &\ \hbox{otherwise;} \hfill \end{matrix} \right. 
\label{rts}
\end{equation}
that is, $R_i=T_i$ equals {\bf S}$^{\bullet }$ related to the interval 
$[i,m]$ if $i-1$ is a white point on the diagram (\ref{gra}),
and $R_i=T_i$ is empty otherwise.

By Claim 7 the set $W^{\hbox{\bf s}}(k,m)$ defined in (\ref{pbse}) is a set of PBW-generators
for {\bf U}$_{\theta }$ over {\bf k}$[G].$ In particular {\bf U}$_{\theta }$ contains 
$\Psi ^{\hbox{\bf s}}(k,m),$ and according to Theorem \ref{26} it is generated over {\bf k}$[G]$
by $\Psi ^{\hbox{\bf s}}(k,m)$ as a right coideal subalgebra.
\label{sex3}
\end{example}

To describe the maximal Hopf subalgebra of {\bf U}$^{\hbox{\bf s}}(k,m)$ 
we need the following definition. 

\begin{definition}  \rm
A black point $s,$ $s\in \, {\bf S}^{\bullet },$ is said to be a $(k,m)$-{\it entrance into {\bf S}} if 
$s-1$ is a white point; that is, either $s=k$ or $s-1\notin \, ${\bf S}. 
\label{ent}
\end{definition}

\begin{lemma}  
The maximal Hopf subalgebra of {\bf U}$^{\hbox{\bf s}}(k,m)$
is generated by all $x_i,$ where $i$ is a $(k,m)$-entrance into {\bf S}.
\label{ent1}
\end{lemma}
\begin{proof} The element $x_i$ belongs to {\bf U}$_{\theta }$ if and only if $i\in T_i.$
 Hence formula (\ref{rts}) shows that $x_i\in \, ${\bf U}$_{\theta }$ if and only if $i$
 is a  $(k,m)$-entrance into {\bf S}.
\end{proof}

\smallskip
\noindent
\underline{If $n=2$} then by Theorem \ref{teor} we have totally $3!=6$ right coideal subalgebras. Among them  $2^2=4$ are ``trivial" cases (Hopf subalgebras) and two more  right coideal subalgebras
are given in  Example \ref{sex1} and Example \ref{sex2}.

\smallskip
\noindent
\underline{If $n=3$} then we have $4!-2^3=16$ proper (not ``trivial") right coideal subalgebras.
In the tableaux below we provide main characteristics of these 16 right coideal subalgebras. 
We mark off by * the ad$_r$-invariant subalgebras {\bf U}$_{\theta }^1=K(\pi ).$
\begin{center}
\begin{tabular}{*{6}{|c}|}
\hline
$r({\bf U})$ & R & T & PBW-generators & r. c. s. generators \\
\hline 
$\ * $ & $R_3=\{ 3\} \hfill $ & $T_3= \{ 3\} \hfill $ & $x_3,\hfill $ & $\ $ \\
$(3,2,1)$ & $R_2=\{ 3\} \hfill $ & $T_2= \{ 3\} \hfill $ & $[x_2x_3],\hfill $ & $[x_1x_2x_3] \hfill $ \\
$\ $ & $R_1=\{ 3\}\hfill $ & $T_1= \{ 3\}\hfill $ & 
$[x_1x_2x_3]\hfill $ & $\circ \,  \circ \, \circ \ \bullet   \hfill $ \\
\hline 

$\ * $ & $R_3=\emptyset  \hfill $ & $T_3=\emptyset  \hfill $ & $\ $ & $\ $ \\
$(3,2,0)$ & $R_2=\{ 2, 3\} \hfill $ & $T_2= \{ 2, 3\} \hfill $ & $x_2,\  [x_3x_2], \hfill $ & $[x_3[x_1x_2]] \hfill $ \\
$\ $ & $R_1=\{ 2, 3\} \hfill $ & $T_1= \{ 2, 3\} \hfill $ & $[x_1x_2],\  [x_3[x_1x_2]] \hfill $ &
 $\circ \,  \circ  \bullet \ \bullet   \hfill $ \\
\hline

$\ $ & $R_3=\{ 3\} \hfill $ & $T_3=\{ 3\} \hfill $ & $x_3,\hfill  $ & $\  $ \\
$(3,1,1)$ & $R_2=\{ 2\} \hfill $ & $T_2= \{ 2, 3\} \hfill $ & $x_2,\ [x_3x_2], \hfill $ & $ [x_1x_2x_3],\ x_2 \hfill $ \\
$\ $ & $R_1=\{ 3\} \hfill $ & $T_1= \{ 3\} \hfill $ & $[x_1x_2x_3] \hfill $ &
 $\circ \,  \circ \, \circ \ \bullet \, , \ \  \cdot \, \circ \, \bullet \ \cdot \hfill $ \\
\hline

$\ $ & $R_3=\emptyset  \hfill $ & $T_3=\emptyset  \hfill $ & $\ $ & $\ $ \\
$(3,1,0)$ & $R_2=\{ 2\} \hfill $ & $T_2= \{ 2\} \hfill $ & $x_2, \hfill $ & $ [x_3x_2x_1],\ x_2 \hfill $ \\
$\ $ & $R_1=\{ 1, 3\} \hfill $ & $T_1= \{ 1, 2, 3\} \hfill $ & $x_1,\  [x_2x_1],\  [x_3x_2x_1]\hfill $ &
 $\circ \,  \bullet \, \bullet \ \bullet \, , \ \  \cdot \, \circ \, \bullet \ \cdot \hfill $ \\
\hline

$\ $ & $R_3=\{ 3\} \hfill $ & $T_3=\{ 3\} \hfill $ & $x_3,\hfill  $ & $\  $ \\
$(3,0,1)$ & $R_2=\emptyset  \hfill $ & $T_2= \emptyset \hfill $ & $\ \hfill $ & $ [[x_2x_3]x_1] \hfill $ \\
$\ $ & $R_1=\{ 1, 3\} \hfill $ & $T_1= \{ 1, 3\} \hfill $ & $x_1,\  [[x_2x_3]x_1] \hfill $ & 
$\circ \,  \bullet \, \circ \ \bullet \hfill $ \\
\hline

$\ * $ & $R_3=\emptyset  \hfill $ & $T_3=\emptyset  \hfill $ & $\ \hfill $ & $\ $ \\
$(3,0,0)$ & $R_2=\emptyset \hfill $ & $T_2= \emptyset \hfill $ & $\ \hfill $ & $ [x_3x_2x_1]\hfill $ \\
$\ $ & $R_1=\{ 1, 2, 3\} \hfill $ & $T_1= \{ 1, 2, 3\} \hfill $ & $x_1,\  [x_2x_1],\  [x_3x_2x_1] \hfill $ & 
$\circ \,  \bullet \, \bullet \ \bullet  \hfill $ \\
\hline

$\ * $ & $R_3=\{ 3\} \hfill $ & $T_3=\{ 3\} \hfill $ & $x_3,\hfill  $ & $\  $ \\
$(2,2,1)$ & $R_2=\{ 3\} \hfill $ & $T_2= \{ 3\} \hfill $ & $[x_2x_3], \hfill $ & $ [x_3x_2x_1],\  [x_2x_3] \hfill $ \\
$\ $ & $R_1=\{ 1, 2, 3\} \hfill $ & $T_1= \{ 1, 2, 3\} \hfill $ & $x_1,\ [x_2x_1],\ [x_3x_2x_1] \hfill $ & 
$\circ \,  \bullet \, \bullet \ \bullet \, , \ \  \cdot \, \circ \,\circ \ \bullet \hfill $ \\
\hline

$\ $ & $R_3=\emptyset  \hfill $ & $T_3=\emptyset  \hfill $ & $\ $ & $\ $ \\
$(2,2,0)$ & $R_2=\{ 2, 3\} \hfill $ & $T_2= \{ 2, 3\} \hfill $ & $x_2,\ [x_3x_2] \hfill $ & $ [x_1x_2],\ [x_3x_2] \hfill $ \\
$\ $ & $R_1=\{ 2\} \hfill $ & $T_1= \{ 2\} \hfill $ & $[x_1x_2]\hfill $ & 
$\circ \,  \circ \, \bullet \ \cdot \, , \ \  \cdot \, \circ \, \bullet \ \bullet \hfill $ \\
\hline

$\ * $ & $R_3=\{ 3\} \hfill $ & $T_3=\{ 3\} \hfill $ & $x_3,\hfill  $ & $\  $ \\
$(2,1,1)$ & $R_2=\{ 2\} \hfill $ & $T_2= \{ 2, 3\} \hfill $ & $x_2, \ [x_3x_2], \hfill $ & $ [x_1x_2],\ x_3 \hfill $ \\
$\ $ & $R_1=\{ 2\} \hfill $ & $T_1= \{ 2, 3\} \hfill $ & $[x_1x_2],\ [x_3[x_1x_2]] \hfill $ & 
$\circ \,  \circ \, \bullet \ \cdot  \, , \ \  \cdot \, \cdot \, \circ \ \bullet \hfill $ \\
\hline

$\ $ & $R_3=\emptyset  \hfill $ & $T_3=\emptyset  \hfill $ & $\ $ & $\ $ \\
$(2,1,0)$ & $R_2=\{ 2 \} \hfill $ & $T_2= \{ 2 \} \hfill $ & $x_2,\hfill $ & $ [x_1x_2]\hfill $ \\
$\ $ & $R_1=\{ 2\} \hfill $ & $T_1= \{ 2\} \hfill $ & $[x_1x_2]\hfill $ & 
$\circ \, \circ \, \bullet \ \cdot \hfill $ \\
\hline

$\ $ & $R_3=\{ 3\} \hfill $ & $T_3=\{ 3\} \hfill $ & $x_3,\hfill  $ & $\  $ \\
$(2,0,1)$ & $R_2=\emptyset \hfill $ & $T_2= \emptyset \hfill $ & $\ \hfill $ & $ [x_2x_1],\ x_3 \hfill $ \\
$\ $ & $R_1=\{ 1, 2\} \hfill $ & $T_1= \{ 1, 2, 3\} \hfill $ & $x_1,\ [x_2x_1],\ [x_3x_2x_1] \hfill $ & 
$\circ \,  \bullet \, \bullet \ \cdot \, , \ \  \cdot \, \cdot  \circ \, \bullet \hfill $ \\
\hline

$\ $ & $R_3=\emptyset  \hfill $ & $T_3=\emptyset  \hfill $ & $\ $ & $\ $ \\
$(2,0,0)$ & $R_2=\emptyset \hfill $ & $T_2= \emptyset \hfill $ & $\ \hfill $ & $ [x_2x_1] \hfill $ \\
$\ $ & $R_1=\{ 1, 2\} \hfill $ & $T_1= \{ 1, 2\} \hfill $ & $x_1,\  [x_2x_1]\hfill $ & 
$\circ \,  \bullet \, \bullet \ \cdot \hfill $ \\
\hline

$\ $ & $R_3=\{ 3\} \hfill $ & $T_3= \{ 3\} \hfill $ & $x_3,\hfill $ & $\ $ \\
$(1,2,1)$ & $R_2=\{ 3\} \hfill $ & $T_2= \{ 3\} \hfill $ & $[x_2x_3],\hfill $ & $[x_2x_3], \  x_1\hfill $ \\
$\ $ & $R_1=\{ 1\}\hfill $ & $T_1= \{ 1, 3\}\hfill $ & $x_1,\  [[x_2x_3]x_1]\hfill $ & 
$\cdot \, \circ  \circ \, \bullet \, , \  \circ  \bullet \, \cdot \, \cdot \hfill $ \\
\hline

\end{tabular}
\end{center}

\begin{center}
\begin{tabular}{*{6}{|c}|}
\hline
$r({\bf U})$ & R & T & PBW-generators & r. c. s. generators \\
\hline

$\ * $ & $R_3=\emptyset  \hfill $ & $T_3=\emptyset  \hfill $ & $\ $ & $\ $ \\
$(1,2,0)$ & $R_2=\{ 2, 3\} \hfill $ & $T_2= \{ 2, 3\}\hfill $ & $x_2,\ [x_3x_2] \hfill $ & $ [x_3x_2],\  x_1\hfill $ \\
$\ $ & $R_1=\{ 1\} \hfill $ & $T_1= \{ 1, 2, 3\} \hfill $ & $x_1,\  [x_2x_1], [x_3x_2x_1]\hfill $ & 
$\cdot \, \circ  \bullet \, \bullet \, , \  \circ  \bullet \, \cdot \, \cdot \hfill $ \\
\hline

$\ $ & $R_3=\{ 3\} \hfill $ & $T_3= \{ 3\} \hfill $ & $x_3,\hfill $ & $\ $ \\
$(0,2,1)$ & $R_2=\{ 3\} \hfill $ & $T_2= \{ 3\} \hfill $ & $[x_2x_3]\hfill $ & $[x_2x_3]\hfill $ \\
$\ $ & $R_1=\emptyset \hfill $ & $T_1= \emptyset \hfill $ & $\ \hfill $ & 
$\cdot \, \circ  \circ \, \bullet  \hfill $ \\
\hline

$\ $ & $R_3=\emptyset  \hfill $ & $T_3=\emptyset  \hfill $ & $\ $ & $\ $ \\
$(0,2,0)$ & $R_2=\{ 2, 3\} \hfill $ & $T_2= \{ 2, 3\}\hfill $ & $x_2,\ [x_3x_2] \hfill $ & $ [x_3x_2]\hfill $ \\
$\ $ & $R_1=\emptyset  \hfill $ & $T_1= \emptyset  \hfill $ & $\ \hfill $ & 
$\cdot \, \circ  \bullet \, \bullet  \hfill $ \\
\hline
\end{tabular}
\end{center}

\section{Triangular decomposition in $U_q(\frak{ sl}_{n+1})$}

At this point it is naturally to conjecture that any ($\Gamma $-homogeneous) right coideal subalgebra 
of $U_q(\frak{ sl}_{n+1})$ (of $u_q(\frak{ sl}_{n+1})$) that contains {\bf k}$\, [H]$ has the triangular decomposition and for any two right coideal subalgebras 
{U}$^-\subseteq U_q^- (\frak{ sl}_{n+1}),$ {U}$^+\subseteq U_q^+ (\frak{ sl}_{n+1})$ 
(respectively, {U}$^-\subseteq u_q^- (\frak{ sl}_{n+1}),$ {U}$^+\subseteq u_q^+ (\frak{ sl}_{n+1})$ ),
the tensor product 
\begin{equation}
 \hbox{\bf U}\, =\, \hbox{\bf U}^-\otimes _{{\bf k}[F]} {\bf k}[H]
\otimes _{{\bf k}[G]}\hbox{\bf U}^+
\label{tru}
\end{equation}
is a right coideal subalgebra. In this hypothesis just one statement fails, the tensor product indeed 
is a right coideal but not always a subalgebra.

\begin{lemma} If $q$ is not a root of $1$ then every 
right coideal subalgebra {\bf U} of $U_q(\frak{ sl}_{n+1}), \, $
 ${\bf U}\supseteq {\bf k}[H]$ has a decomposition
$(\ref{tru}),$ where {\bf U}$^+\supseteq {\bf k}[G]$ and 
{\bf U}$^-\supseteq {\bf k}[F]$ are right coideal subalgebras of
$U_q^+(\frak{ sl}_{n+1})$ and $U_q^-(\frak{ sl}_{n+1})$ respectively.
If $q$ has finite multiplicative order $t>2,$ 
then this is the case for $\Gamma $-homogeneous
right coideal subalgebras of $u_q(\frak{ sl}_{n+1}).$
\label{raz2}
\end{lemma}
\begin{proof} 
Due to the triangular decompositions (\ref{tr}), (\ref{tr1})
the values of super-letters $[x_kx_{k+1}\ldots x_m],$ $[x_k^-x_{k+1}^-\ldots x_m^-]$
form a set of PBW-generators over {\bf k}$[H]$ for $U_q(\frak{ sl}_{n+1}).$ 

Let us fix the following order on the skew-primitive generators
\begin{equation} 
x_1>x_2>\ldots >x_n>x_1^->x_2^->\ldots >x_n^-.
\label{orr}
\end{equation}
By  Lemma \ref{nco1} and Proposition \ref{pro} (see the arguments above Eq. (\ref{vad22}))
the subalgebra {\bf U} has PBW-generators of the form
\begin{equation}
c=[u]+\sum \alpha _iW_i+\sum_j\beta_jV_j \in \hbox{\bf U},
\label{vad10}
\end{equation}
where $W_i$ are the basis super-words starting with less than $[u]$ super-letters, 
$D(W_i)=D(u),$ and $V_j$ are $G$-super-words of $D$-degree less than $D(u),$ while
the leading term $[u]$ equals either $[x_kx_{k+1}\ldots x_m]$ or 
$[x_k^-x_{k+1}^-\ldots x_m^-].$ 
Certainly the leading terms here are defined by the degree function into the
additive monoid $\Gamma ^+\oplus \Gamma ^-$ generated by $X\cup X^-$
(but not into the group $\Gamma $!).
In particular all $W_i$ in (\ref{vad10}) have the same constitution in 
$X\cup X^-$ as the leading term $[u]$ does. Thus all $W_i$'s and 
the leading term $[u]$
belong to the same component of the triangular decomposition. Hence
it remains to show that there are no  terms $V_j.$

If $q$ is not a root of 1 then
by Corollary \ref{odn11} the algebra {\bf U} is $\Gamma $-homogeneous.
Hence (in both cases) the PBW-generators may be chosen to be $\Gamma $-homogeneous as well.  
In this case all terms $V_j$ have the same $\Gamma $-degree and smaller 
$\Gamma ^+\oplus \Gamma ^-$-degree. However this is impossible.

Indeed, if the leading term is $[x_k^-x_{k+1}^-\ldots x_m^-]$ 
then the $\Gamma ^+\oplus \Gamma ^-$-degree
of $V_j$ should be less than $x_k^-+x_{k+1}^-+\ldots +x_m^-.$ Hence due to
definitions (\ref{orr}) and (\ref{ord}) we have $V_j\in U_q^-(\frak{ sl}_{n+1}),$
(respectively, $V_j\in u_q^-(\frak{ sl}_{n+1})$),
and the  $\Gamma $-degree of $V_j$ coincides with the $\Gamma ^+\oplus \Gamma ^-$-degree.
A contradiction.

Suppose that the leading term is $[x_k^+x_{k+1}^+\ldots x_m^+].$
 Let $d=\sum s_ix_i+\sum r_ix_i^-$
be the $\Gamma ^+\oplus \Gamma ^-$-degree of $V_j.$ Since 
\begin{equation} 
d<x_k+x_{k+1}+\ldots +x_m,
\label{ine}
\end{equation}
we have $s_i=0$ provided that $i<k.$ Since the $\Gamma $-degree of $V_j$ equals
$x_k+x_{k+1}+\ldots +x_m,$ we have $s_i-r_i=1$ provided that $k\leq i\leq m,$
and $r_i=s_i$ otherwise. Inequality (\ref{ine}) implies
$s_k\leq1.$ Together with $s_k-r_k=1$ this yields $s_k=1,$ $r_k=0.$
Hence, again inequality (\ref{ine}) implies
$s_{k+1}\leq1,$ and again together with $s_{k+1}-r_{k+1}=1$ 
this yields $s_{k+1}=1,$ $r_{k+1}=0.$ In this way we get $s_i=1,$ $r_i=0,$
$k\leq i\leq m,$ that contradicts (\ref{ine}).

Now we see that each PBW-generator (\ref{vad10})
belongs to either $U_q^+(\frak{ sl}_{n+1})$ or $U_q^-(\frak{ sl}_{n+1})$ 
(respectively, to either $u_q^+(\frak{ sl}_{n+1})$ or $u_q^-(\frak{ sl}_{n+1})$). 
Therefore {\bf U} has the decomposition (\ref{tru}). 
\end{proof}  

\noindent
{\bf Remark}. Certainly, if {\bf U} does not contain {\bf k}$\, [H]$ then decomposition
(\ref{tru}) is not valid. A {\bf k}-subalgebra  generated by $g_1f_1$
has no triangular decomposition of any type provided that $f_1\notin G,$ $g_1\notin F.$
In the single parameter case (when $G=F=H$) the subalgebra with 1 generated by 
$g^{-1}_1(x_1^-+x_1)$ provides a similar example. Moreover 
in the single parameter case many of the (left) coideal subalgebras 
studied by M. Noumi  and G. Letzter do not admit a triangular decomposition.

\smallskip
\begin{corollary} 
If $q$ is not a root of 1, then $U_q(\frak{ sl}_{n+1})$ has
just a finite number of right coideal subalgebras containing the coradical.
If $q$ has finite multiplicative order $t>2,$ then
$u_q(\frak{ sl}_{n+1})$ has
just a finite number of $\Gamma $-homogeneous 
right coideal subalgebras containing the coradical.
\label{fin2}
\end{corollary}
\begin{proof}
This follows from the above lemma and Corollary \ref{fin1} applied to $U_q^{\pm }(\frak{ sl}_{n+1}),$
$u_q^{\pm }(\frak{ sl}_{n+1}).$
 \end{proof}

Our next goal is to understand when tensor product (\ref{tru}) is a subalgebra
and then to find a way to calculate the total number of ($\Gamma $-homogeneous)
right coideal subalgebras.
\begin{lemma} The tensor product $(\ref{tru})$
is a right coideal subalgebra if and only if 
\begin{equation}
 [{\bf U}^+,{\bf U}^-]\subseteq \, {\bf U}^-\otimes _{{\bf k}[F]} {\bf k}[H]
\otimes _{{\bf k}[G]} {\bf U}^+.
\label{rud}
\end{equation}
\label{raz1}
\end{lemma}
\begin{proof} Of course if {\bf U} is a subalgebra then (\ref{rud}) holds.
Conversely, it is clear that {\bf U} is a right coideal. Relation (\ref{rud})
implies $u^+\cdot v^-$ $=[u^+, v^-]$ $+p(u^+,v^-) v^-\cdot u^+$ $\in \, ${\bf U}, 
where $u^+\in {\bf U}^+,$ $v^-\in {\bf U}^-.$
Hence $(u^-\cdot u^+)(v^-\cdot v^+)$ $=u^-(u^+\cdot v^-)v^+$
$\in \,${\bf U}, with arbitrary $v^+\in {\bf U}^+,$ $u^-\in {\bf U}^-.$

Since ${\bf U}={\bf U}^-\cdot H\cdot {\bf U}^+,$
it remains to check that ${\bf U}^-\cdot H$ $=H\cdot {\bf U}^-,$ and
${\bf U}^+\cdot H$ $=H\cdot {\bf U}^+.$ 
Since ${\bf U}^+$ contains $G,$ it is homogeneous with respect to the grading 
(\ref{grad}). If $u\in ({\bf U}^+)^{\chi },$ $f\in F,$ then $uf=\chi (f)fu.$
Hence  ${\bf U}^+\cdot F=F\cdot {\bf U}^+.$ Similarly 
${\bf U}^-\cdot G=G\cdot {\bf U}^-.$ \end{proof}

\section{Consistency condition}
In this section we are going to find sufficient condition for consistency relation
 (\ref{rud}) to be valid. 
In what follows we denote by $\Psi_-^{\hbox{\bf s}}(i,j)$ a polynomial that appears from
$\Psi ^{\hbox{\bf s}}(i,j)$ given in (\ref{cbr1}) under the substitutions $x_t\leftarrow x_t^-,$
 $1\leq t\leq n$ with skew commutators defined by (\ref{sqo}) in $U_q^-(\frak{ sl}_{n+1}).$ By
pr$W^{\hbox{\bf s}}(k,m)$ (respectively, pr$W_-^{\hbox{\bf s}}(i,j)$)
we denote a subspace spanned  by proper derivatives 
of $\Psi ^{\hbox{\bf s}}(k,m)$ (respectively, of $\Psi_-^{\hbox{\bf s}}(i,j)$), see Theorem \ref{26}. 
Consider two elements 
$\Psi ^{\hbox{\bf s}}(k,m)$ and $\Psi _-^{T}(i,j).$ Let us display them graphically
as defined in (\ref{gra}):
\begin{equation}
\begin{matrix}
\stackrel{k-1}{\circ } \ & \cdots \ & \stackrel{i-1}{\bullet } 
\ & \stackrel{i}{\bullet }\ \ & \stackrel{i+1}{\circ }\ & \cdots &
\ & \stackrel{m}{\bullet } \ & \ & \stackrel{j}{\cdot } \cr
\ \ & \  \ & \circ  
\ & \circ \ \ & \bullet \ & \cdots &
\ & \bullet  \ & \cdots  \ & \bullet
\end{matrix}
\label{gra1}
\end{equation}
We shall prove that 
\begin{equation}
 [\Psi ^{\hbox{\bf s}}(k,m),\Psi _-^{T}(i,j)]\in 
{\rm pr}\, W_-^{T}(i,j)\cdot {\rm pr}\, W^{\hbox{\bf s}}(k,m)
\label{rdi}
\end{equation}
if one of the following two options fulfills:

a) Representation (\ref{gra1}) has no fragments of the form \label{use}
\begin{equation}
\begin{matrix}
\stackrel{t}{\circ } \ & \cdots & \stackrel{l}{\bullet } \cr
\circ  
\ & \cdots  & \bullet 
\end{matrix}
\label{gra2}
\end{equation}

b) Representation (\ref{gra1}) has the form \label{use1}
\begin{equation}
\begin{matrix}
\stackrel{k-1}{\circ } \ & \cdots & \circ & \cdots & \bullet & \cdots & \stackrel{m}{\bullet } \cr
\circ  
\ & \cdots  &  \bullet & \cdots & \circ & \cdots &  \bullet 
\end{matrix}
\label{gra3}
\end{equation}
(in particular $i=k,$ $j=m$),
where no one intermediate column has points of the same color.

Suppose that diagram (\ref{gra1}) satisfies condition a). In this case all black-black
columns are located before all white-white columns. Let us choose the closest
black-black and white-white pair of columns. Then (\ref{gra1}) takes up the form
\begin{equation}
\underbrace{
\begin{matrix}\ & \  &\circ & \bullet  
& \bullet & \cdots 
& \stackrel{t}{\bullet } \cr 
\cdots & \circ &
\bullet & \bullet  
& \circ & \cdots 
& \bullet 
\end{matrix}}_{\hbox{mainly black}}\,
\underbrace{
\begin{matrix} \stackrel{t+1}{\circ } & \bullet & \cdots \cr
\bullet & \circ & \cdots 
\end{matrix}}_{\hbox{equality}}\, 
\underbrace{
\begin{matrix} \stackrel{l}{\circ } & \circ & \bullet  
& \circ & \bullet  & \cdots \cr
\circ & \circ & \circ 
& \bullet & \
\end{matrix}}_{\hbox{mainly white}}
\label{gra4}
\end{equation}
Here in the ``mainly black" zone there are no white-white columns; in the ``equality" zone
we have just black-white, and white-black columns; while in the ``mainly white"
zone there are no black-black columns. Of course the ``mainly black" zone may be empty.
In this case we may omit the ``equality" zone as well, since all the diagram has no 
black-black columns at all. In the same way the ``mainly white" zone may be empty too.

Recall that in Definition \ref{eski} for a fixed pair $(k,m)$ 
we define ${\bf S}_{\circ }=\, (${\bf S}$\, \cap [k,m-1])\, \cup \, \{ k-1 \} ,$
while ${\bf S}^{\bullet }=\, (${\bf S}$\, \cap [k,m-1])\, \cup \, \{ m \} ;$
respectively $s_0=k-1\in \, {\bf S}_{\circ },$ and 
$s_{r+1}=m\in \, {\bf S}^{\bullet }.$ 

All black-black columns are labeled by numbers from ${\bf S}^{\bullet }\cap T^{\bullet},$
where the bullets correspond to the pairs $(k,m)$ and $(i,j),$ respectively. 
Similarly all white-white columns are labeled by numbers from
 $(\overline{\bf S})_{\circ  }\cap (\overline{T})_{\circ },$
where $\overline{\bf S}, \overline{T}$ are the complements of ${\bf S},T$ with respect to $[k, m-1],$
$[i,j-1],$ respectively.
Thus condition a) is equivalent to the inequality 
\begin{equation}
\sup \{ {\bf S}^{\bullet }\cap T^{\bullet}\} <\inf \{ (\overline{\bf S})_{\circ  }\cap (\overline{T})_{\circ }\} .
\label{ineq}
\end{equation}
We are reminded that the supremum and infimum of the empty set equal $-\infty $
and $\infty ,$ respectively.

Condition b), in turn, means that $i=k,$ $j=m,$ $T=\overline{\bf S}.$

\smallskip
We go ahead with a number of useful notes. If $u$ is a word in $X,$ then by $u^-$ we denote
a word in $X^-$ that appears from $u$ under the substitution $x_i\leftarrow x_i^-.$
We have $p(v,w^-)=\chi ^v(f_w)=p(w,v),$ while $p(w^-,v)=(\chi ^w)^{-1}(g_v)=p(w,v)^{-1}.$
Thus  $p(v,w^-)p(w^-,v)=1.$ Therefore the Jacobi and antisymmetry identities take up their
original ``colored" form (see, (\ref{jak1})):
\begin{equation} 
[[u,v],w^-]=[u,[v,w^-]]+p_{wv}[[u,w^-],v];
\label{uno}
\end{equation}
\begin{equation} 
[u,w^-]=-p_{wu}[w^-,u].
\label{dos}
\end{equation}
In the same way
\begin{equation} 
[[u^-,v^-],w]=[u^-,[v^-,w]]+p_{vw}^{-1}[[u^-,w],v^-].
\label{tres}
\end{equation}
Using antisymmetry and (\ref{tres}) we have also
\begin{equation} 
[u,[v^-,w^-]]=[[u,v^-],w^-]+p_{vu}[v^-[u,w^-]].
\label{cua}
\end{equation}
In these relations $u,v,w$ are words in $X.$ To simplify further
calculations we may extend the brackets to the set of all 
$H$-words: We put $\chi ^{hu}=\chi ^{u},$ $g_{hu}=hg_u,$
$h\in H,$ and define the skew-brackets by the same formula (\ref{sqo}).
In this case we have
\begin{equation} 
[u, hv]=\chi ^u(h)\, h[u,v], \ \ \ h\in H;
\label{cuq1}
\end{equation}
\begin{equation} 
[hu,v]=h[u,v]+p_{uv}(1-\chi^v(h))\, h\, v\cdot u,  \ \ \ h\in H.
\label{cuq2}
\end{equation}
\begin{equation} 
[hu,v]=\chi^v(h)\, h[u,v]+(1-\chi^v(h))\, h\, u\cdot v,  \ \ \ h\in H.
\label{cuq21}
\end{equation}
To calculate the coefficients it is convienient to have in mind the following 
consequences of (\ref{a1rel}):
\begin{equation} 
\chi ^k(g_{k-1}f_{k-1})=\chi ^{k-1}(g_{k}f_{k})=q^{-1}, \ \ \chi^k(g_if_i)=1,\  \hbox{ if } |i-k|>1.
\label{cuq22}
\end{equation}
Of course all basic formulae (\ref{jak1}), (\ref{jak3}), (\ref{br1}) and their consequences 
remain valid. However we must stress that
once we apply relations (\ref{rela3}), or other ``inhomogeneous in $H$" relations
(for example the third option of (\ref{derm1}), see below), 
we have to fix the curvature of the brackets as soon as
the inhomogeneous substitution applies to the right factor in the brackets:
\begin{equation}
[u,[x_i,x_i^-]]=u(1-g_if_i)-\chi ^u(g_if_i)(1-g_if_i)u=(1-\chi^u(g_if_i))u,
\label{cuq3}
\end{equation}
but not $[u,[x_i,x_i^-]]=[u,1-g_if_i]=[u,1]-[u,g_if_i]=0.$ At the same time
\begin{equation}
[[x_i,x_i^-],u]=(1-g_if_i)u-u(1-g_if_i)=(\chi^u(g_if_i)-1)\, g_if_i\cdot u,
\label{cuq4}
\end{equation}
and $[[x_i,x_i^-],u]=[1-g_if_i, u]$ $=[1,u]-[g_if_i,u]$ is valid since the inhomogeneous
substitution has been applied to the left factor in the brackets. In what follows we shall 
denote for short $h_i=g_if_i,$ and $\bar{h}_{ki}=h_kh_{k+1}\cdots h_{i-1},$
where $k<i.$

\smallskip
Now we consider relation (\ref{rdi}) when $i=j.$

\begin{proposition}  If $k<m,$ then 
$[\Psi ^{\hbox{\bf s}}(k,m),x_i^-]\in {\rm pr}W^{\hbox{\bf s}}(k,m)$ if and only if
$i$ is not a $(k,m)$-entrance into {\bf S} $($see Definition $\ref{ent}).$
\label{sin}
\end{proposition}
\begin{proof} 
To prove the proposition we shall need the following four formulae (\ref{derm1}--\ref{derm4}).
In that formulae the notation $a\sim b$ means $a=\alpha b,$ $0\neq \alpha\, $ $\in \,${\bf k}. 
\begin{equation}
[u[k,m], x_i^-]\sim \left\{ \begin{matrix}0,\hfill & \hbox{if }k<i<m;\hfill \cr
h_k \cdot u[k+1,m], \hfill & \hbox{if } i=k<m; \hfill \cr 
u[k,m-1], \hfill & \hbox{if } k<i=m.\hfill \end{matrix} \right.
\label{derm1}
\end{equation}

For $i=k$ we have
$[u[k,m],x_k^-]=[[x_k,u[k+1,m]],x_k^-].$
If we denote $u=x_k,$ $v=u[k+1,m],$ $w^-=x_k^-,$ then 
formulae (\ref{uno}) and (\ref{cuq4}) imply $[[u,v],w^-]$ $=\mu _{k,m}[[u,w^-],v]$ 
$=\mu _{k,m}\, (q^{-1}-1)\,h_k\cdot v$ 
since $p(v,x_k)p(x_k,v)=q^{-1}$ due to (\ref{a1rel}).

For $k<i=m$ we put $u=u[k,m-1],$ $v=x_m,$ $w^-=x_m^-.$ Then 
$[u[k,m],x_k^-]=[[u,v],w^-],$ while (\ref{uno}) and (\ref{cuq3}) imply
$$
[[u,v],w^-]=[u,[v,w^-]]=(1-q^{-1})u
$$
since according to (\ref{a1rel}) we have $p(u,x_m)p(x_m,u)$ $=q^{-1}.$

Finally, for $k<i<m$ we put $u=u[k,i-1],$ $v=u[i,m],$ $w^-=x_i^-.$ Then 
$[u[k,m],x_k^-]$ $=[[u,v],w^-],$ while (\ref{uno}) and considered above case $i=k$ 
with (\ref{cuq1}) imply 
$$
[[u,v],w^-]=[u,[v,w^-]]\sim [u, h_iu[i+1,m]]
$$
$$
=h_i\chi^u(h_i)\cdot [u,u[i+1,m]]=0,
$$
which proves (\ref{derm1}).
\begin{equation}
[\Psi ^{\hbox{\bf s}}(k,m), x_k^-]\sim \left\{ \begin{matrix}\Psi ^{\hbox{\bf s}}(k+1,m),
\hfill &\hbox{ if } s_1=k;\hfill \cr
h_k\cdot \Psi ^{\hbox{\bf s}}(k+1,m), \hfill &\hbox{ if } s_1>k. \hfill \end{matrix} \right.
\label{derm2}
\end{equation}
Let us put $u=\Psi ^{\hbox{\bf s}}(1+s_1,m),$ $v=u[k,s_1],$ $w^-=x_k^-.$
The definition (\ref{cbr1}) shows that $\Psi ^{\hbox{\bf s}}(k,m)=[u,v],$
while (\ref{uno}) implies $[[u,v],w^-]=[u,[v,w^-]].$

 If $s_1=k,$ then by (\ref{cuq3}) we have $[u,[v,w^-]]$
$=u-\chi ^u(h_k)u=(1-q^{-1})u.$

If $s_1>k,$ then (\ref{derm1}) yields
$[v,w^-]\sim h_k\cdot u[k+1,s_1].$  Therefore $[u,[v,w^-]]\sim h_k[u, u[k+1,s_1]],$
see (\ref{cuq1}).
It remains to note that $[u, u[k+1,s_1]]=\Psi ^{\hbox{\bf s}}(k+1,m)$ due to (\ref{cbr1}).
 Thus formula (\ref{derm2}) is proved.
\begin{equation}
[\Psi ^{\hbox{\bf s}}(k,m), x_m^-]\sim 
\left\{ \begin{matrix}h_m\cdot \Psi ^{\hbox{\bf s}}(k,m-1),\hfill & \hbox{if }k\leq s_r=m-1;\hfill \cr
\Psi ^{\hbox{\bf s}}(k,m-1), \hfill & \hbox{if } k\leq s_r<m-1. \hfill \end{matrix} \right.
\label{derm3}
\end{equation}
Let us put $u=u[1+s_r,m],$ $v=\Psi ^{\hbox{\bf s}}(k,s_r),$ $w^-=x_m^-.$
Decomposition (\ref{cbrr}) with $i=r$ shows that $\Psi ^{\hbox{\bf s}}(k,m)=[u,v],$
while (\ref{uno}) implies $[[u,v],w^-]\sim [[u,w^-],v].$

If $s_r=m-1,$ then by (\ref{cuq4}) we get $[[u,w^-],v]\sim h_m \cdot v$
since $p(v,x_m)p(x_m,v)=q^{-1}\neq 1.$
 
If $s_r\neq m-1,$ then by (\ref{derm1}) we have 
$[u,w^-]\sim u[1+s_r,m-1].$ Therefore by decomposition (\ref{cbrr}), we get
$[[u,w^-],v]$ $\sim \Psi ^{\hbox{\bf s}}(k,m-1),$ that proves (\ref{derm3}).
\begin{equation}
[\Psi ^{\hbox{\bf s}}(k,m), x_i^-]\sim
\left\{ \begin{matrix}h_i\cdot \Psi ^{\hbox{\bf s}}(i+1,m)\cdot  \Psi ^{\hbox{\bf s}}(k,i-1), &\!\!\!
i-1\in \, {\bf S},\ i\notin \, {\bf S};\hfill \cr
\Psi ^{\hbox{\bf s}}(i+1,m)\cdot 
\Psi ^{\hbox{\bf s}}(k,i-1), \hfill &\!\!\!
i-1\notin \, {\bf S},\ i\in \, {\bf S}; \hfill \cr
0,\hfill &\!\!\! \hbox{otherwise,}\hfill \end{matrix}
\right.
\label{derm4}
\end{equation}
where of course $k<i<m.$

Suppose that $i-1\in \, ${\bf S}. 
Let us put $u=\Psi ^{\hbox{\bf s}}(i,m),$ $v=\Psi ^{\hbox{\bf s}}(k,i-1),$ $w^-=x_i^-.$
Decomposition (\ref{cbrr}) shows that $\Psi ^{\hbox{\bf s}}(k,m)=[u,v],$
while (\ref{uno}) implies $[[u,v],w^-]\sim [[u,w^-],v].$

If $i\notin \, ${\bf S}, then we apply the second option of (\ref{derm2}) and (\ref{cuq3}):
$$
[[u,w^-],v]=[h_i\cdot \Psi ^{\hbox{\bf s}}(i+1,m), v]\sim h_i\, \Psi ^{\hbox{\bf s}}(i+1,m)\cdot v
$$
This proves the first option of (\ref{derm4}).

If $i\in \, ${\bf S} (and still $i-1\in \, ${\bf S}) then we may use the first option of (\ref{derm2}):
$$
[[u,w^-],v] \sim [\Psi ^{\hbox{\bf s}}(i+1,m), v]=0.
$$

Suppose that $i-1\notin \, ${\bf S}.  If also $i\notin \,${\bf S}, then (\ref{derm4})
follows from the first case of (\ref{derm1}) and definition (\ref{cbr1}).

If $i\in \, ${\bf S}, we put $u=\Psi ^{\hbox{\bf s}}(1+i,m),$ $v=\Psi ^{\hbox{\bf s}}(k,i),$ $w^-=x_i^-.$
Decomposition (\ref{cbrr}) shows that $\Psi ^{\hbox{\bf s}}(k,m)=[u,v],$
while (\ref{uno}) implies $[[u,v],w^-]=[u,[v,w^-]].$ 
By (\ref{derm3}) we have
$[v,w^-]\sim \Psi ^{\hbox{\bf s}}(k,i-1).$ We may use (\ref{derm3}) 
taking into account the curvature of the skew commutator:
$$
[u,[v,w^-]]\sim u\cdot \Psi ^{\hbox{\bf s}}(k,i-1)-\chi ^u(g_vf_w)\, \Psi ^{\hbox{\bf s}}(k,i-1)\cdot u. 
$$   
Since $\Psi ^{\hbox{\bf s}}(k,i-1)$ and $u$ are separated elements, we have 
 $$
u\cdot \Psi ^{\hbox{\bf s}}(k,i-1)=\alpha \, \Psi ^{\hbox{\bf s}}(k,i-1)\cdot u
$$
with the coefficient $\alpha =\chi ^u(g_kg_{k+1}\cdots g_{i-1}).$
In this case $\alpha ^{-1} \chi ^u(g_vf_w)=\chi ^u(h_i)$ 
$=p(u,x_i)p(x_i,u)$ $=q^{-1}\neq 1,$ which completes the proof of (\ref{derm4}).

\smallskip

\smallskip
If $i$ is a $(k,m)$-entrance into {\bf S}, then according to definition (\ref{pbse}) we have
$\Psi ^{\hbox{\bf s}}(k,i)\in {\rm pr}W^{\hbox{\bf s}}(k,m),$ while 
$\Psi ^{\hbox{\bf s}}(k,i-1)\notin {\rm pr}W^{\hbox{\bf s}}(k,m).$
Hence (\ref{derm2}--\ref{derm4}) imply 
$[\Psi ^{\hbox{\bf s}}(k,m),x_i^-]\notin {\rm pr}W^{\hbox{\bf s}}(k,m).$

If we compare formulae (\ref{derm2}--\ref{derm4}) with (\ref{der4}--\ref{der7}),
we see that $[\Psi ^{\hbox{\bf s}}(k,m), x_i^-]\notin {\rm pr}W^{\hbox{\bf s}}(k,m)$
only in the following three cases:
a) if $i=k$ and $k\in \, ${\bf S}; 
b) if $i=m$ and $i-1\notin \, ${\bf S};
c) if $i-1\notin \,${\bf S} and $i\in \,${\bf S}, \ $k<i<m.$
In all these cases $i$ by definition is a $(k,m)$-entrance into {\bf S}. \end{proof}

\smallskip
Now by means of the following lemmas we are going to show 
that (\ref{rdi}) fulfills provided that diagram (\ref{gra1})
takes up the form (\ref{gra4}). 

\begin{lemma}  
If ${\bf S}^{\bullet }\cap T^{\bullet }={\bf S}_{\circ }\cap T_{\circ }=\emptyset $ then
\begin{equation}
 [\Psi ^{\hbox{\bf s}}(k,m),\Psi _-^{T}(i,j)]=0.
\label{rdi1}
\end{equation}
Here ${\bf S}^{\bullet },$  ${\bf S}_{\circ }$ and $T^{\bullet },$ $T_{\circ }$ correspond to the pair $(k,m)$ and $(i,j)$ respectively, see Definition ${\ref{eski}}.$
\label{zer}
\end{lemma}
\begin{proof}
Since relations (\ref{a1rel}) are invariant
under the substitutions $p_{ij}\leftarrow p_{ji}^{-1},$ $q\leftarrow q^{-1},$
formulae (\ref{derm1}--\ref{derm4}) remain valid under the substitutions
$x_i\leftarrow x_i^-,$ $x_i^-\leftarrow x_i.$ 

We show firstly that
\begin{equation}
[u[k,m], u[i,j]^-]=0,\ \ \hbox{ if }k\neq i ,\ \ j\neq m.
\label{dm1}
\end{equation}
Indeed, if $k<i\leq j<m,$ or $i<k\leq m<j,$ this follows from the first case of (\ref{derm1}) and its dual.

If $k<i=m<j,$ we have $[u[k,m], u[m,j]^-]$ $=[u[k,m],[x_m^-, u[m+1,j]^-]]$
$=[[u[k,m],x_m^-], u[m+1,j]^-]$ $\sim [u[k,m-1],u[m+1,j]^-]=0.$

If $k<i<m<j,$ the  induction on $m-i$ provides
$$[u[k,m], u[i,j]^-]=[u[k,m],[x_i^-, u[i+1,j]^-]]=0.$$

The remaining case,  $i<k<j<m,$ due to (\ref{dos}), is dual 
to one considered above, $k<i<m<j.$ Thus, (\ref{dm1}) is proved.

If $S^{\bullet }\cap T^{\bullet }=S_{\circ }\cap T_{\circ }=\emptyset ,$
then (\ref{dm1}) implies 
$$[u[1+s_a,s_{a+1}],u[1+t_b,t_{b+1}]^-]=0,$$
where $s_a,$ $t_b$ are given in definition (\ref{cbr1}). Hence 
(\ref{uno}-\ref{cua}) imply (\ref{rdi1}). \end{proof}

\begin{lemma}  
If ${\bf S}^{\bullet }\cap T^{\bullet }=\emptyset ,$ while  
${\bf S}_{\circ }\cap T_{\circ }\neq \emptyset $ then
$$
 [\Psi ^{\hbox{\bf s}}(k,m),\Psi _-^{T}(i,j)]
$$
\begin{equation}
=\Psi ^T_-(i,\nu -1)\left(
\sum _{a=\nu +1}^{\mu +1}\alpha _{a}\,    
\bar{h}_{\nu \, a}\, \Psi ^T_-(a,j)\cdot \Psi ^{\hbox{\bf s}}(a,m)\right) \Psi ^{\hbox{\bf s}}(k,\nu -1)
\label{rdi2}
\end{equation}
where $\mu =\min \{ m,j\} ,$ $\nu =\max \{ k,i \} ,$
while $\alpha _{a}=0$ provided that  $a-1\in \, {\bf S}\, \cup T$
with the only exception,  $\alpha _{\mu +1}\neq 0.$
Here by definition
$h_i=g_if_i,$ $\bar{h}_{\nu \,  a}$ $=h_{\nu }h_{\nu +1}\cdots h_{a-1},$ 
and  $\Psi ^{\hbox{\bf s}}(k,k-1)=$ $\Psi ^T_-(i,i-1)=1$
 $($in particular always either the first or the last factor in $(\ref{rdi2})$ is trivial$\, ).$
\label{zer1}
\end{lemma}
\begin{proof} We start with a particular case:
\begin{equation}
[u[k,m], u[k,j]^-]=\sum _{a=k+1}^{\mu +1}\alpha _a \,  
\bar{h}_{ka} \, u[a,j]^-\cdot u[a,m],\ j\neq m,
\label{dm3}
\end{equation}
where $\mu ={\rm min}\{ m, j\} ,$  
$\alpha_a\in \,${\bf k}, $\alpha_a\neq 0,$
and $u[j+1,j]^-=u[m+1,m]=1.$ 

If $k=j,$ the formula follows from (\ref{derm1}).
If $k=m,$  the formula follows from  (\ref{dos}) and dual (\ref{derm1}). 

In the general case we use induction on $j-k.$ Denote $u=u[k,m],$ $v^-=x_k^-,$ $w^-=u[k+1,j]^-.$
The left hand side of (\ref{dm3}) equals $[u, [v^-,w^-]].$ According to (\ref{dm1})
we have $[u,w^-]=0.$ Hence $[u, [v^-,w^-]]=[[u, v^-],w^-].$ Using (\ref{derm1})
and (\ref{cuq2}) we have  
$$
[[u, v^-],w^-]\sim [h_k\cdot u[k+1,m],w^-]
$$
$$
\sim h_k\cdot [u[k+1,m],w^-]]+\alpha_{k+1} \, h_k\, w^-\cdot u[k+1,m],
$$
where due to formula (\ref{cuq2}) and relations (\ref{cuq22})
 the coefficient $\alpha_{k+1} $ equals $p(u(k+1),m),w^-)(1-q)\neq 0.$
To prove (\ref{dm3}) it remains 
to apply the inductive supposition to the first summand.

Our next step is to prove (\ref{rdi2}) for {\bf S}$\, \neq \emptyset ,$ $T=\emptyset $ and $i=k:$
\begin{equation}
[\Psi ^{\hbox{\bf s}}(k,m), u[k,j]^-]=\sum _{a=k+1}^{\mu +1}\alpha _a   \,
\bar{h}_{ka} \, u[a,j]^-\cdot \Psi ^{\hbox{\bf s}}(a,m),\ \ \hbox{ if }j\notin {\bf S}^{\bullet },
\label{dm31}
\end{equation}
where $\mu ={\rm min }\{ m, j \},$ while   $\alpha _a=0$ if  
$a-1\in \, ${\bf S}$\, \cap \, [k,m-1],$ and  formally 
$\Psi ^{\hbox{\bf s}}(m+1,m)=[j+1,j]^-=1.$

Suppose firstly that $j>m.$ In this case $\mu =m.$ We use induction on $m-k.$ If $m$ $=k,$
the formula follows from dual (\ref{derm1}). If $m>k,$ we put
$u=\Psi ^{\hbox{\bf s}}(1+s_1,m),$
$v=u[k,s_1],$
$w^-=u[k,j]^-.$
According to (\ref{cbr1}) we have $[u,v]$ $=\Psi ^{\hbox{\bf s}}(k,m),$
while $[u,w^-]=0$ due to (\ref{rdi1}) with $T=\emptyset ,$ $T_{\circ }=\{ a-1\} .$
By Jacobi identity (\ref{uno}) and (\ref{dm3})
we get
\begin{equation}
[\Psi ^{\hbox{\bf s}}(k,m),w^-]=[u,[v,w^-]]
=\hbox{\huge [}u, \sum _{a=k+1}^{1+s_1}\alpha _a\, \bar{h}_{ka}\, 
u[a,j]^-\cdot u[a,s_1]\hbox{\huge ]}.
\label{rrr}
\end{equation}
Relation (\ref{rdi1}) with $T=\emptyset ,$ $T_{\circ }=\{ k-1\} $ implies
$[u,u[a,j]^-]$ $=0$ unless $a=1+s_1.$ Using ad-identities (\ref{br1}), (\ref{cuq1})
we may continue
$$
=\alpha \, \bar{h}_{k\, 1+s_1}\, [u,u[1+s_1,j]^-]
+\sum _{a =k+1}^{s_1}\alpha_a \, \bar{h}_{ka}\, u[a,j]^-\cdot \Psi ^{\hbox{\bf s}}(a,m),
$$
where $\alpha \neq 0.$
Since $1+s_1>k,$ we may apply the inductive supposition to the first summand.

Let $j<m ,$ say $s_l<j<s_{l+1},$ $0\leq l\leq r+1.$
Suppose that
$l=r;$ that is, $s_r<j<m.$ 
Denote 
$u=u[1+s_r,m],$
$v=\Psi ^{\hbox{\bf s}}(k,s_r),$
$w^-=u[k,j]^-.$
By definition $[u,v]$ $=\Psi ^{\hbox{\bf s}}(k,m),$
while $[u,w^-]=0$ due to (\ref{dm1}). By Jacobi identity (\ref{uno}) and 
the considered above case
$j>m$ with $m\leftarrow s_l$ we get
$$
[\Psi ^{\hbox{\bf s}}(k,m),w^-]=[u,[v,w^-]]=
\hbox{\large [}u,\sum _{a=k+1}^{1+s_l}
\alpha_a \, \bar{h}_{ka}\, u[a,j]^-\cdot \Psi ^{\hbox{\bf s}}(a,s_l) \hbox{\large ]},
$$
where $\alpha _a=0$ if  $a-1\in \,  \hbox{\bf S}\, \cap \, [k,s_l-1].$
 Relation (\ref{dm1}) imply
$[u,u[a,j]^-]=0$ unless $a=1+s_l.$
Hence by ad-identities (\ref{br1}), (\ref{cuq1}) we may continue
$$
=\alpha \, \bar{h}_{k\, 1+s_l}\, [u,u[1+s_l,j]]+
\sum _{a=k+1}^{s_l}
\alpha_a \, \bar{h}_{ka}\, u[a,j]^-\cdot \Psi ^{\hbox{\bf s}}(a,m),
$$
where $\alpha \neq 0,$  $\alpha _a=0$ if $a-1\in \,  \hbox{\bf S}\, \cap \, [k,s_l-1].$
It remains to apply (\ref{dm3}) to the first summand.

If, finally, $l\leq  r,$ we put
$u=\Psi ^{\hbox{\bf s}}(1+s_{l+1},m),$
$v=\Psi ^{\hbox{\bf s}}(k,s_{l+1}),$
$w^-=u[k,j]^-.$
Then still $[u,v]$ $=\Psi ^{\hbox{\bf s}}(k,m),$ see decomposition (\ref{cbrr}),
and $[u,w^-]=0.$ By the considered above case with $m\leftarrow s_{l+1}$ we have
$$
[v,w^-]=
\sum _{a=k+1}^{j+1}
\alpha_a \, \bar{h}_{ka}\, u[a,j]^-\cdot \Psi ^{\hbox{\bf s}}(a,s_{l+1}),
$$
where $\alpha _a=0$ if $a-1\in \, {\bf S}\, \cap [k,s_{l+1}-1].$
In this case $[u,u[a,j]^-]=0$ since $j<s_{l+1}.$
Thus by ad-identities (\ref{br1}), (\ref{cuq1}), and decomposition (\ref{cbrr}) we get 
$$
[u,[v,w^-]]=\sum _{a=k+1}^{j+1}
\alpha_a\, \bar{h}_{ka}\, u[a,j]^-\cdot \Psi ^{\hbox{\bf s}}(a,m),
$$
where $\alpha _a=0$ if $a-1\in \, {\bf S}\, \cap [k,s_{l+1}-1].$
Of course
$$({\bf S}\, \cap [k,m-1])\cap [k,j]
=({\bf S}\, \cap [k,s_{l+1}-1])\cap [k,j],
$$
hence in the above sum $\alpha _a=0$  if $a-1\in \, {\bf S}\, \cap [k,m-1].$
Since Jacobi identity (\ref{uno}) implies $[\Psi ^{\hbox{\bf s}}(k,m),w^-]=[[u,v],w^-],$
formula (\ref{dm31}) is proved.

\smallskip
Now we are ready to consider the general case. Since $S^{\bullet }\cap T^{\bullet }=\emptyset ,$
the intersection $S_{\circ }\cap T_{\circ  }$ equals either $\{ k-1\}$ or $\{ i-1\} .$
More precisely, this intersection has the only point $\{ \nu -1\}.$

Suppose firstly that $i=k.$ In this case we shall prove 
(\ref{rdi2}) by induction on $j-k.$   If $j=k,$ one may apply the second option of (\ref{derm2}).
Let $j>i=k.$ If $T=\emptyset ,$ the formula is already proved, see (\ref{dm31}). 
Suppose that $T\neq \emptyset .$
Let us denote $u=\Psi ^{\hbox{\bf s}}(k,m),$ $v^-=\Psi ^{T}_-(1+t_1,j),$
$w^-=u[k,t_1].$ Then according to definition (\ref{cbr1}) the left hand side of (\ref{rdi2})
equals $[u,[v^-,w^-]].$  By (\ref{rdi1}) we have $[u,v^-]=0.$ Hence Jacobi identity
(\ref{cua}) implies $[u,[v^-,w^-]]\sim [v^-,[u,w^-]].$ Using (\ref{dm31}) we get
\begin{equation}
[u, w^-]=\sum _{a=k+1}^{\mu +1}\alpha _a  \,
\bar{h}_{ka} \, u[a,j]^-\cdot \Psi ^{\hbox{\bf s}}(a,m),
\label{dmr}
\end{equation}
where $\mu ={\rm min }\{ m, t_1\},$ while   $\alpha _a=0$ if 
$a-1\in \, ${\bf S}$\, \cap \, [k,m-1],$ and  formally 
$\Psi ^{\hbox{\bf s}}(m+1,m)=[j+1,j]^-=1.$ 
Of course, $m\neq t_1$ since $S^{\bullet }\cap T^{\bullet }=\emptyset .$

If $m>t_1,$ then we have
$$
[v^-,[u, w^-]=\alpha _{1+t_1}  \,
[v^-, \bar{h}_{k\, 1+t_1} \, \Psi ^{\hbox{\bf s}}(1+t_1,m)]+
\sum _{a=k+1}^{t_1}\alpha _a  \,
[ v^-,\bar{h}_{ka} \, u[a,t_1]^-\cdot \Psi ^{\hbox{\bf s}}(a,m)].
$$
By (\ref{dos}) and  (\ref{rdi1}) we have $[v^-, \Psi ^{\hbox{\bf s}}(a,m)]=0$ if $k<a\leq t_1.$
Hence by (\ref{cuq1}) and (\ref{br1}) we get
$
[ v^-,\bar{h}_{ka} \, u[a,t_1]^-\cdot \Psi ^{\hbox{\bf s}}(a,m)]
\sim \bar{h}_{ka}\, \Psi ^T_-(a,j)\cdot \Psi ^{\hbox{\bf s}}(a,m).
$
It remains to apply (\ref{cuq1}) and then the inductive supposition to the first summand.

If $m<t_1,$ then 
$$
[v^-,[u, w^-]=
\sum _{a=k+1}^{m+1}\alpha _a  \,
[ v^-,\bar{h}_{ka} \, u[a,t_1]^-\cdot \Psi ^{\hbox{\bf s}}(a,m)].
$$
Since $m<t_1,$ we have $[v^-, \Psi ^{\hbox{\bf s}}(a,m)]=0.$
Hence, again by (\ref{cuq1}) and (\ref{br1}), we get
$$
[ v^-,\bar{h}_{ka} \, u[a,t_1]^-\cdot \Psi ^{\hbox{\bf s}}(a,m)]
\sim \bar{h}_{ka}\, \Psi ^T_-(a,j)\cdot \Psi ^{\hbox{\bf s}}(a,m),
$$
which completes the case $i=k.$

Suppose that $i>k.$ In this case $k-1\in T\cap [i,j-1],$ say $k=1+t_l.$
Let us put $u=\Psi ^{\hbox{\bf s}}(k,m),$ $v^-=\Psi ^T_-(k,j),$ $w^-=\Psi ^T_-(i,k-1).$
Decomposition (\ref{cbrr}) implies $\Psi ^T_-(i,j)=[v^-,w^-].$ Since $[u,w^-]=0,$
we have $[u,[v^-,w^-]]=[[u,v^-],w^-].$ To find $[u,v^-]$ we may use already considered case: 
$$
[u,v^-]=\sum _{a=k+1}^{\mu +1}\alpha _a\bar{h}_{ka}\Psi ^T_-(a,j)\cdot \Psi ^{\hbox{\bf s}}(a,m)
$$
with $\alpha _a=0$ if $a-1\in {\bf S}\cup T,$ $a\neq \mu +1.$ Certainly
$[\Psi ^{\hbox{\bf s}}(a,m), w^-]$ $=[\Psi ^T_-(a,j), w^-]=0$ since $a>k.$
By means of (\ref{cuq22}) we have
$$
\chi ^{w^-}(\bar{h}_{ka})=\chi ^{w}_-(h_kk_{k+1}\cdots h_{a-1})=\chi ^{k-1}_-(h_k)=q\neq 1.
$$
Now formula (\ref{cuq2}) shows that $[[u,v^-],w^-]\sim w^-\cdot [u,v^-],$ which is required.

In perfect analogy, if $k<i$ then $i-1\in {\bf S}\cap [k,m-1],$ $i=1+s_l.$ 
Let us put $u=\Psi ^{\hbox{\bf s}}(i,m),$ $v=\Psi ^{\hbox{\bf s}}(k,i-1),$ $w^-=\Psi ^T_-(i,j).$
Decomposition (\ref{cbrr}) implies $\Psi ^{\hbox{\bf s}}(k,m),=[u,v].$ 
Since $[v,w^-]=0,$ by Jacobi identity (\ref{uno}) we have $[[u,v],w^-]\sim [[u,w^-],v].$ 
To find $[u,w^-]$ we use already considered case: 
$$
[u,w^-]=\sum _{a=i+1}^{\mu +1}\alpha _a\bar{h}_{ia}\Psi ^T_-(a,j)\cdot \Psi ^{\hbox{\bf s}}(a,m)
$$
with $\alpha _a=0$ if $a-1\in {\bf S}\cup T,$ $a\neq \mu +1.$ 
Of course,
$[\Psi ^{\hbox{\bf s}}(a,m), v]$ $=[\Psi ^T_-(a,j), v]=0$ since $a>i.$
By means of (\ref{cuq22}) we have
$$
\chi ^{v}(\bar{h}_{ia})=\chi ^{v}(h_ih_{i+1}\cdots h_{a-1})=\chi ^{i-1}(h_i)=q^{-1}\neq 1.
$$
Now formula (\ref{cuq21}) shows that $[[u,w^-],v]\sim [u,w^-]\cdot v,$ which is required.
\end{proof}

We have mentioned  above that our main concepts (Definition \ref{eski}) are not invariant
with respect to the replacement of $x_i,$ $x_i^-,$ $1\leq i\leq n$
by $y_i,$ $y_i^-,$ $1\leq i\leq n,$ where by definition $y_i=x_{\varphi (i)},$
$y_i^-=x_{\varphi (i)}^-,$ $\varphi (i)=n-i+1.$
Hence the application of already proved lemmas to generators $y_i,$ $y_i^{-}$
provides an additional information. In this way we are going to prove the following two statements.

\begin{lemma} 
If $(\overline{\bf S})_{\circ }\cap (\overline{T})_{\circ }
=(\overline{\bf S})^{\bullet }\cap (\overline{T})^{\bullet }=\emptyset $ then
\begin{equation}
 [\Psi ^{\hbox{\bf s}}(k,m),\Psi _-^{T}(i,j)]=0.
\label{rdi1b}
\end{equation}
\label{zerb}
\end{lemma}

\begin{lemma}  
If 
$(\overline{\bf S})_{\circ }\cap (\overline{T})_{\circ }= \emptyset ,$ while
$(\overline{\bf S})^{\bullet }\cap (\overline{T})^{\bullet }\neq \emptyset $
 then
$$
 [\Psi ^{\hbox{\bf s}}(k,m),\Psi _-^{T}(i,j)]
$$
\begin{equation}
=\Psi ^T_-(\mu +1,j)\left(
\sum _{b=\nu -1}^{\mu -1}\alpha _{b}\,    
\bar{h}_{b+1 \, \mu +1}\, \Psi ^T_-(i,b)\cdot \Psi ^{\hbox{\bf s}}(k,b)\right) \Psi ^{\hbox{\bf s}}(\mu +1,m)
\label{rdi2b}
\end{equation}
where $\mu =\min \{ m,j\} ,$ $\nu =\max \{ k,i \} ,$
while $\alpha _{b}=0$ provided that  $b\notin {\bf S} \cap T$
with the only exception,  $\alpha _{\nu -1}\neq 0.$
\label{zerb1}
\end{lemma}
\begin{proof}
To derive these statements from Lemma \ref{zer} and Lemma  \ref{zer1}
we use  the  {\it decoding} lemma, Lemma \ref{dec}.
Let us apply (\ref{decod}) to the left hand side of (\ref{rdi2b}). In this case we have
$$
({\overline{\varphi(\hbox{\bf S})-1}})_{\circ }=\varphi \{ (\overline{\bf S})^{\bullet }+1\},\ 
({\overline{\varphi(\hbox{\bf S})-1}})^{\bullet }=\varphi \{ (\overline{\bf S})_{\circ }+1\},
$$
where in the left hand sides the operators (Definition \ref{eski}) correspond
to $(\varphi (m), \varphi (k)), $ while in the right hand sides to $(k,m).$ In particular 
$(\overline{\bf S})^{\bullet }\cap (\overline{T})^{\bullet }=\emptyset $ is equivalent to
$({\overline{\varphi (\hbox{\bf S})-1}})_{\circ }\cap ({\overline{\varphi (T)-1}})_{\circ } =\emptyset ,$
while
$(\overline{\bf S})_{\circ }\cap (\overline{T})_{\circ }=\emptyset $ is equivalent to 
$({\overline{\varphi (\hbox{\bf S})-1}})^{\bullet }\cap ({\overline{\varphi (T)-1}})^{\bullet } =\emptyset .$
Hence we may  
use relations (\ref{rdi1}), (\ref{rdi2}). Relation (\ref{rdi1}) proves Lemma \ref{zerb}.
In the case of Lemma \ref{zerb1} we 
again apply decoding formula (\ref{decod}) in order to get (\ref{rdi2b}). \end{proof}

\begin{proposition} 
Condition $a),$ p.$\pageref{use}$ implies $(\ref{rdi}).$
\label{coa}
\end{proposition}
\begin{proof} 
We have seen that condition a) is equivalent to inequality (\ref{ineq}).

If ${\bf S}^{\bullet }\cap {T}^{\bullet }=\emptyset $ 
then we may use Lemma \ref{zer} and Lemma \ref{zer1}. Let us show that all factors in
(\ref{rdi2}) in terms with nonzero $\alpha _a$ belong to either pr$\, W^{T}_-(i,j)$
or  pr$\, W^{\hbox{\bf s}}(k,m).$ 

If $a=\mu +1,$ say $a-1=m<j,$ then $a-1=m$
is a white point on the diagram of $\Psi ^{T}_-(i,j).$ Hence by Theorem \ref{26}
we have $\Psi ^{T}_-(m+1,j)\in {\rm pr}\,W^{T}_-(i,j).$ In the same way $a-1=j<m$
implies $\Psi ^{\hbox{\bf s}}(j+1,m) \in {\rm pr}\, W^{\hbox{\bf s}}(k,m)$
(we note that $j\neq m$ since ${\bf S}^{\bullet }\cap {T}^{\bullet }=\emptyset $ yet).

If $a-1\notin {\bf S}\cup T,$ $a\leq \min \{ m,j\}$ then $a-1$ is a white-white point,
hence again Theorem \ref{26} implies $\Psi ^{T}_-(a,j)\in {\rm pr}\, W^{T}_-(i,j)$
and $\Psi ^{\hbox{\bf s}}(a,m) \in {\rm pr}\, W^{\hbox{\bf s}}(k,m).$

For the first and the last factors (if $k\neq i,$ otherwise they do not exist) we have 
${\bf S}_{\circ }\cap T_{\circ }=\{ \nu -1\} $
(since ${\bf S}^{\bullet }\cap {T}^{\bullet }=\emptyset $).
Hence, if $i<k=\nu $ then we have $\nu -1\in T\cap [i,j-1],$
and by Theorem \ref{26} the first factor belongs to pr$\, W^{T}_-(i,j).$
Similarly, if $k<i=\nu $ then  the last factor belongs to pr$\, W^{\hbox{\bf s}}(k,m)$
since $\nu -1\in {\bf S}\cap [k,m-1].$

In perfect analogy if $(\overline{\bf S})_{\circ }\cap (\overline{T})_{\circ }=\emptyset ,$
Lemma \ref{zerb}, Lemma \ref{zerb1}, and Theorem \ref{26} imply (\ref{rdi}).

Suppose that both ${\bf S}^{\bullet }\cap T^{\bullet }\neq \emptyset ,$
and  $(\overline{\bf S})_{\circ }\cap (\overline{T})_{\circ }\neq \emptyset .$ Then the diagram (\ref{gra1})
takes up the form (\ref{gra4}) with nonempty ``mainly black" and ``mainly white" zones.

Let us put $u=\Psi ^{\hbox{\bf s}}(t+1,m),$ $v=\Psi ^{\hbox{\bf s}}(k,t),$
$w^-=\Psi ^T_-(t+1,j),$ $z^-=\Psi ^T(i,t),$
where $t$ is the label of the last ``black-black" column 
(that is, $t=\max \{ {\bf S}^{\bullet }\cap T^{\bullet }\}$).
We have $[u,z^-]=0,$ $[v,w^-]=0.$ Using decomposition (\ref{cbrr})
and Jacobi identities (\ref{uno}) and (\ref{cua})  we may write
$$
[\Psi ^{\hbox{\bf s}}(k,m),\Psi ^T_-(t+1,j)]=[[u,v],[w^-,z^-]]
$$
$$
=[u,[v,[w^-,z^-]]]+p_{wz,\, v}[[u,[w^-,z^-]],v]
$$
$$
=\alpha [u,[w^-,[v,z^-]]]+\beta [[[u,w^-],z^-],v].
$$
To find $[v,z^-]$ and $[u,w^-]$ we may use (\ref{rdi2b}) and (\ref{rdi2}) respectively:
\begin{equation}
[v,z^-]=\sum _{b=\nu -1}^{t-1}\alpha _{b}\,    
\bar{h}_{b+1 \, \mu +1}\, \Psi ^T_-(i,b)\cdot \Psi ^{\hbox{\bf s}}(k,b);
\label{alre}
\end{equation}
$$
[u,w^-]=
\sum _{a=t+2}^{\mu +1}\alpha _{a}\,    
\bar{h}_{t+1 \, a}\, \Psi ^T_-(a,j)\cdot \Psi ^{\hbox{\bf s}}(a,m).
$$
By means of (\ref{cuq1}) and ad-identity (\ref{br1}) we have $[w^-,[v,z^-]]=0$
since $w$ is separated both from $\Psi ^T(i,b),$ and $\Psi ^{\hbox{\bf s}}(k,b).$

At the same time (\ref{cuq2}) and ad-identity (\ref{br1f}) imply
$$
[[u,w^-],z^-]=\sum _{a=t+2}^{\mu +1}\alpha _{a}\, \beta_a \,   
\bar{h}_{t+1 \, a}\, z^-\cdot \Psi ^T_-(a,j)\cdot \Psi ^{\hbox{\bf s}}(a,m)
$$
since due to (\ref{cuq22}) we have $\chi^z_-(\bar{h}_{t+1 \, a})=q\neq 1.$
Again by  (\ref{cuq2}) and (\ref{br1f}), taking into account that 
$v$ is separated both from $\Psi ^{\hbox{\bf s}}(a,m),$ and $\Psi ^T(a,j),$
we find 
$$
[\bar{h}_{t+1 \, a}\, z^-\cdot \Psi ^T_-(a,j)\cdot \Psi ^{\hbox{\bf s}}(a,m), v]
=\bar{h}_{t+1 \, a}\, z^-\cdot  \Psi ^T_-(a,j)\cdot v\cdot \Psi ^{\hbox{\bf s}}(a,m) 
$$
$$
+\varepsilon \bar{h}_{t+1 \, a}\, z^-\cdot \Psi ^T_-(a,j)\cdot [v,z^-]\cdot \Psi ^{\hbox{\bf s}}(a,m).
$$
By Theorem (\ref{26}) we have $v\in {\rm pr}\, W^{\hbox{\bf s}}(k,m),$ and $z^-\in {\rm pr}\, W^{T}_-(i,j),$
while $[v,z^-]$ already appears in (\ref{alre}). Since all factors in (\ref{alre}) with nonzero 
$\alpha _b$ belong to either pr$\, W^{\hbox{\bf s}}(k,m),$ or pr$\, W^{T}_-(i,j),$ we get
$[[[u,w^-],z^-],v]\in {\rm pr}\, W^{T}_-(i,j)\cdot {\rm pr}\, W^{\hbox{\bf s}}(k,m).$ 
The proposition is completely proved.\end{proof}

\begin{proposition}  We have
\begin{equation}
\left[ \Psi ^{\hbox{\bf s}}(k,m),\Psi_-^{\overline{\hbox{\bf s}}}(k,m)\right] \sim 1-\bar{h}_{k\, m+1}.
\label{sh}
\end{equation}
In particular condition $b)$ p.$\pageref{use1}$ implies $(\ref{rdi}).$
\label{cob}
\end{proposition}
\begin{proof} We fix $k$ and use induction on $m.$ If $m=k$ the statement is clear.
Let us put $u=\Psi ^{\hbox{\bf s}}(k,m-1),$ $v=x_m,$ $w^-=x_m^-,$ 
$z^-=\Psi_-^{\overline{\hbox{\bf s}}}(k,m-1)$ (in fact $v=w$). 
In this case $[v,z^-]=0,$ $[u,w^-]=0.$

Suppose firstly that $m-1\notin {\bf S}.$ In this case $m-1\in \overline{\bf S},$ 
hence  (\ref{cin}) implies 
$$
\Psi ^{\hbox{\bf s}}(k,m)\sim [\Psi ^{\hbox{\bf s}}(k,m-1),x_m]=[u,v].
$$ 
$$
\Psi_-^{\overline{\hbox{\bf s}}}(k,m)=[x_m^-,\Psi_-^{\overline{\hbox{\bf s}}}(k,m-1)]=[w^-,z^-].
$$
By means of Jacobi identity (\ref{cua}) we have 
$$
[[u,v],[w^-,z^-]]=[u,[v,[w^-,z^-]]]+p_{wz,\, v}[[u,[w^-,z^-]],v]
$$
\begin{equation}
=[u,[[v,w^-],z^-]]]+p_{wz,\, v}p_{w,\, u}[[w^-,[u,z^-]],v].
\label{sh1}
\end{equation}
Relations (\ref{a1rel}) imply  $p_{wz,\, v}p_{w,\, u}=p_{wv}\cdot p_{zv}p_{wu}=q\cdot q^{-1}=1.$
Using (\ref{cuq4}) and then (\ref{cuq1}) we get
$$
[u,[[v,w^-],z^-]]]=(\chi^z(h_m)-1)\chi^u(h_m)h_m[u,z^-]=(1-q^{-1})\varepsilon (1-\bar{h}_{k\,m}),
$$
where the inductive supposition yields $[u,z^-]=\varepsilon (1-\bar{h}_{k\,m}),$ $\varepsilon \neq 0.$
Using again the inductive supposition and taking
 into account the comment to (\ref{cuq3}), we have
$[w^-,[u,z^-]]=\varepsilon (1-\chi^m(\bar{h}_{k\, m}))x^-_m.$ Thus
$$
[[w^-,[u,z^-]],v]=\varepsilon (1-\chi^m_-(\bar{h}_{k\, m}))[x^-_m,x_m]=
\varepsilon (1-q)(-q^{-1})(1-h_m).
$$
Now (\ref{sh1}) implies 
$$
\left[ \Psi ^{\hbox{\bf s}}(k,m),\Psi_-^{\overline{\hbox{\bf s}}}(k,m)\right]\sim
[[u,v],[w^-,z^-]]=\varepsilon (1-q^{-1})(1-\bar{h}_{k\, m+1}).
$$

If $m-1\in {\bf S},$ then, of course, $m-1\notin \overline{\bf S},$ 
hence  (\ref{cin}) implies 
$$
\Psi ^{\hbox{\bf s}}(k,m)=[x_m,\Psi ^{\hbox{\bf s}}(k,m-1)]=[v,u].
$$
$$
\Psi_-^{\overline{\hbox{\bf s}}}(k,m)\sim [\Psi_-^{\overline{\hbox{\bf s}}}(k,m-1),x_m^-]=[z^-,w^-].
$$
By means of Jacobi identity (\ref{cua}) we have 
$$
[[v,u],[z^-,w^-]]=[v,[u,[z^-,w^-]]]+p_{wz,\, u}[[v,[z^-,w^-]],u]
$$
\begin{equation}
=[v,[[u,z^-],w^-]]]+p_{wz,\, u}p_{z,\, v}[[z^-,[v,w^-]],u].
\label{sh2}
\end{equation}
Relations (\ref{a1rel}) imply  $p_{wz,\, u}p_{z,\, v}=p_{wu}p_{zv}\cdot p_{zu}=q^{-1}\cdot q=1.$
Using the inductive supposition and (\ref{cuq1}), we have
$$
[v,[[u,z^-],w^-]]]=\varepsilon (\chi^m_-(\bar{h}_{k\, m})-1)[x_m,\bar{h}_{k\, m}x_m^-]
=\varepsilon (q-1)q^{-1}\bar{h}_{k\, m}(1-h_m).
$$
Using (\ref{cuq3}), (\ref{dos}) and then the inductive supposition we get
$$
[[z^-,[v,w^-]],u]=(1-\chi^z_-(h_m))[z^-,u]=\varepsilon (1-q)(-p_{zu}^{-1}) (1-\bar{h}_{k\, m})
$$
$$
=-\varepsilon (1-q)q^{-1}(1-\bar{h}_{k\, m})=\varepsilon (1-q^{-1})(1-\bar{h}_{k\, m}).
$$
Now (\ref{sh2}) implies 
$
[[v,u],[z^-,w^-]]=\varepsilon (1-q^{-1})(1-\bar{h}_{k\, m+1}).
$
The proposition is proved.\end{proof}

\section{Right coideal subalgebras in $U_q(\frak{ sl}_{n+1}), u_q(\frak{ sl}_{n+1})$}
\begin{theorem}  Let $U_{\theta }^+,$ $U_{\theta ^{\prime}}^-$ be right coideal subalgebras 
of positive and negative quantum Borel subalgebras
defined over the related coradical by $r$-sequences 
$\theta =(\theta _1, \theta _2, \ldots , \theta _n),$ and 
$\theta ^{\prime }=(\theta _1^{\prime }, \theta _2^{\prime }, \ldots , \theta _n^{\prime }),$
respectively. The tensor product 
\begin{equation}
 \hbox{\bf U}\, =\, U_{\theta ^{\prime }}^-\otimes _{{\bf k}[F]} {\bf k}[H]
\otimes _{{\bf k}[G]} U_{\theta }^+
\label{tru1}
\end{equation}
is a right coideal subalgebra  if and only if for each pair $(k,i),$
$1\leq k,i\leq n$ one of the following two conditions is fulfilled. In that conditions the sets
$T_k,$ $T_i^{\prime }$ are defined according to Definition $\ref{tski}$
by $\theta ,$ $\theta ^{\prime },$ 
respectively,  while $\tilde{\theta }_k=k+\theta _k-1,$ 
$\tilde{\theta }_k^{\prime }=k+\theta _k^{\prime }-1.$
\begin{eqnarray} 
& &
\begin{matrix}
1.\ \hfill & \sup \left\{ a\, |\, k\leq a\leq  \tilde{\theta }_k,\, i\leq a\leq  \tilde{\theta }_i^{\prime },
\, a\in T_k, \,  a\in T_i^{\prime } \right\} \hfill \cr 
\ \ & < \inf \left\{ b\, |\, k-1\leq b< \tilde{\theta }_k,\, i-1\leq b< \tilde{\theta }_i^{\prime },
\, b\notin T_k, \,  b\notin T_i^{\prime } \right\} \hfill ;
\end{matrix}
\label{yo1}  \\
& &
\begin{matrix} 
2. \hfill &\ i=k,\, \theta _k^{\prime }=\theta _k, \hbox{ and}\hfill \cr 
\ \hfill & \left\{ a\, |\, k\leq a< \tilde{\theta }_k,\, a\in T_k, \, a\in T_k^{\prime } \right\} =\emptyset ,\hfill \cr 
\ \hfill & \left\{ a\, |\, k\leq a< \tilde{\theta }_k,\,  a\notin T_k, \,  a\notin T_k^{\prime } \right\} =\emptyset .
\end{matrix}
\label{yo2} 
\end{eqnarray}

If $q$ is not a root of $1,$ 
then every right coideal subalgebra {\bf U}$\supseteq {\bf k}[H]$ of $U_q(\frak{ sl}_{n+1})$
has form $(\ref{tru1})$ with $\theta ,$ $\theta ^{\prime }$ satisfying the above property.
If $q$ has finite multiplicative order $t>2,$ then this is the case for $\Gamma $-homogeneous
right coideal subalgebras of $u_q(\frak{ sl}_{n+1}).$
\label{osn5}
\end{theorem}
\begin{proof} By Lemma \ref{raz1} we have to show that 
$[U_{\theta }^+,U_{\theta ^{\prime }}^-]\subseteq {\bf U}.$
The first condition in the theorem means that $\Psi ^{T_k}(k,\tilde{\theta }_k)$ and
$\Psi_-^{T_i^{\prime }}(k,\tilde{\theta }_i^{\prime })$ satisfy condition a) p.\pageref{use},
while the second one is equivalent to condition b) p.\pageref{use1} for these elements.
By definition $[k:\tilde{\theta }_k]$ is a simple $U_{\theta }^+$-root, such that any other 
simple  $U_{\theta }^+$-root of the form $[k:m]$ satisfies $m<\tilde{\theta }_k.$
Similarly each simple $U_{\theta ^{\prime }}^-$-root of the form $[i:j]$
satisfies $j\leq \tilde{\theta }_i^{\prime }.$ Condition a)  certainly remain valid
for subdiagrams, while proper subdiagrams of (\ref{gra3}) satisfy condition a).
Therefore, due to Proposition \ref{coa} and Proposition \ref{cob}, for each pair of a simple 
$U_{\theta }^+$-root, $[k:m],$ and a simple $U_{\theta ^{\prime }}^-$-root, $[i:j],$ we have 
\begin{equation}
 [\Psi ^{T_k}(k,m),\Psi _-^{T_i^{\prime }}(i,j)]\in {\bf U}.
\label{rdii}
\end{equation}
By Claim 8 the algebras $U^+_{\theta },$ and $U^-_{\theta ^{\prime }}$
 are generated by $\Psi ^{T_k}(k,m),$
and $\Psi _-^{T_i^{\prime }}(i,j),$ respectively, where $[k:m]$ and $[i:j]$ run through 
the sets of simple roots. 
To show that $[U^+_{\theta },U^-_{\theta ^{\prime }}]\subseteq {\bf U},$
it remains to apply ad-identities (\ref{br1f}), (\ref{br1}) and evident induction
on degree (we remark that in (\ref{rdi}) the degree of factors diminishes).

Conversely, suppose that $[U^+_{\theta },U^-_{\theta ^{\prime }}]\subseteq {\bf U}.$
Let us choose any pair $(k,i),$ and denote 
$$
t=\sup \left\{ a\, |\, k\leq a\leq  \tilde{\theta }_k,\, i\leq a\leq  \tilde{\theta }_i^{\prime },
\, a\in T_k, \,  a\in T_i^{\prime } \right\},
$$
$$ 
l= \inf \left\{ b\, |\, k-1\leq b< \tilde{\theta }_k,\, i-1\leq b< \tilde{\theta }_i^{\prime },
\, b\notin T_k, \,  b\notin T_i^{\prime } \right\}.
$$
 If one of these sets is empty then 
condition (\ref{yo1}) is valid. Suppose that $t<l.$ The point $t$ is a white-white point 
on the diagram (\ref{gra1}) of the elements $\Psi ^{T_k}(k,\tilde{\theta }_k),$ 
$\Psi _-^{T_i^{\prime }}(i,\tilde{\theta }_i^{\prime }),$ while $l$ is a black-black one.
In particular, due to Theorem \ref{26}, we have 
$z\stackrel{df}{=}\Psi ^{T_k}(1+t,l)\in U^+_{\theta },$
$z^{\prime }\stackrel{df}{=}\Psi_-^{T_k^{\prime }}(1+t,l)\in U^-_{\theta ^{\prime }}.$
Using both decompositions (\ref{cbrr}) and (\ref{cbry}) we get
\begin{equation}
 \Psi ^{T_k}(k,\tilde{\theta }_k)\sim [\Psi ^{T_k}(1+l,\tilde{\theta }_k),
[\Psi ^{T_k}(k,t),z]].
\label{ooo}
\end{equation}
\begin{equation}
\Psi_-^{T_i^{\prime }}(i,\tilde{\theta }_i^{\prime })\sim 
[\Psi_-^{T_i^{\prime }}(1+l,\tilde{\theta }_i^{\prime }),
[\Psi_-^{T_i^{\prime }}(i,t),z^{\prime }]].
\label{oo1}
\end{equation}
Since between $t$ and $l$ there are no columns marked by one color,
we may apply Proposition \ref{cob}, $[z,z^{\prime }]=\varepsilon (1-h),$
where $h=g_zf_z\in H.$ This relation and (\ref{ooo}), (\ref{oo1}) allow us to consider $z$ and 
$z^-=\varepsilon ^{-1}z^{\prime }$ as a pair of new variables and turn to
a new set of variables
$$\{ x_1,\ldots ,x_t, z, x_{l+1},\ldots , x_n, x_1^-,\ldots ,x_t^-, z^-, x_{l+1}^-,\ldots , x_n^-\}.$$
In particular we may apply formula (\ref{derm4}):
$$
\Psi ^{T_k}(l+1,\tilde{\theta }_k)\cdot \Psi ^{T_k}(k,t)\sim 
[\Psi ^{T_k}(k,\tilde{\theta }_k),z^-]\in U^+_{\theta }\cap {\bf U}=U^+_{\theta }.
$$
This implies that both $[l+1:\tilde{\theta }_k]$ and $[k:t]$ are $U^+_{\theta }$-roots.
Hence the simple $U^+_{\theta }$-root $[k:\tilde{\theta }_k]$ is a sum of three 
$U^+_{\theta }$-roots, $[k:t]+[1+t:l]+[1+l:\tilde{\theta }_k].$ This is a contradiction,
unless $t=k-1,$ $\tilde{\theta }_k=l.$ In perfect analogy we have $t=i-1,$
$\tilde{\theta }^{\prime }_k=l;$ that is, condition (\ref{yo2}) is valid.

The last statement follows from Lemma \ref{raz2}. \end{proof}

\

{\bf Acknowledgements}

The autors were supported by PAPIIT IN 108306-3, UNAM, M\'exico and
partially by the macro-project ``Tecnolog\'\i as para la Universidad de la Informaci\'on y la Computaci\'on, UNAM, M\'exico''. We thank Dr. Marco
V. Jos\'e for the facility given to use the supercomputer of the
 ``Instituto de Investigaciones Biom\'edicas" of the UNAM and Juan R. 
 Bobadilla for valuable
 technical assistance. We also thank the referee for the careful
 reading of the paper and constructive suggestions.


\begin{thebibliography}{90}

\bibitem{AG} N. Andruskiewitsch, M. Gra\~na, Braided Hopf algebras over non abelian finite group,
Bol. Acad. Nac. Cienc. Cordoba, 63(1999), 46--78.

\bibitem{AS} N. Andruskiewitsch, H.-J. Schneider, Pointed Hopf algebras, in: S. Montgomery,
H.-J. Schneider (Eds.) New Directions in 
Hopf Algebras, MSRI Publications, 43(2002), 1--68.

\bibitem{BGH} N. Bergeron, Y. Gao, N. Hu, Drinfel'd doubles and Lusztig's symmetries of two-parameter quantum groups, Journal of Algebra, 301(2006), 378--405.

\bibitem{CM} W. Chin, M. Musson, Multiparameter quantum enveloping algebras,
Journal of Pure and Applied Algebra, 107(1996), 171--191.

\bibitem{CLMT} W. Chin, R. Larson, M. Musson, J. Towber,
The first two terms of the coradical filtration of multiparameter ${\mathfrak gl}_N,$
Comm. Algebra, 24,N2(1996), 3845--3883.

\bibitem{CV} M. Costantini, M. Varagnolo, Quantum double and multiparameter quantum groups,
Comm. Algebra, 22, N15(1994), 6305--6321. 

\bibitem{Dij} M.S. Dijkhuizen, Some remarks on the construction of quantum symmetric spaces,
In: Representations of Lie groups, Lie Algebras and Their Quantum Analogies,
Acta Appl. Math., 44, N1-2(1996), 59--80. 

\bibitem{DN} M.S. Dijkhuizen and M. Noumi, A family of quantum projective spaces
and related $q$-hypergeometric orthogonal polynomials,
Trans. Amer. Math. Soc., 350, N8(1998), 3269--3296.

\bibitem{FMP} V.O. Ferreira, L.S.I. Murakami, and A. Paques, A Hopf-Galois correspondence for
free algebras, J. Algebra, 276(2004), 407--416.

\bibitem{FG} D. Flores and J.A. Green, Quantum symmetric algebras, Algebras and Representation Theory, 4, N1(2001), 55--76. 

\bibitem{FG1} D. Flores de Chela, J.A. Green, Quantum symmetric algebras II,
Journal of Algebra, 269(2003), 610--631.

\bibitem{Flo} D. Flores de Chela, Quantum symmetric algebras as braided Hopf algebras, 
in:  Algebraic Structures
and Their Representations, J.A. de la Pe\~{n}a, E. Vallejo, N Atakishieyev (Eds),
Contemporary Mathematics, 376, AMS, Providence, 2005, 261--271. 

\bibitem{GMO} M. Grefraths, G. McGuire, M. O'sullivan, On Plotkin-optimal codes 
over finite Frobenius rings, Journal of Algebra and Its Applications, 5, N6(2006), 799--815.

\bibitem{HR} R.G. Heyneman, D. E. Radford, Reflexivity and coalgebras of finite type, J. of Algebra,
29(1974), 215-246.

\bibitem{Jac} N. Jacobson, Lie Algebras, Intercsience Publishers, New York --- London, 1962.

\bibitem{JL} A. Joseph and G. Letzter, Separation of variables for quantized enveloping algebras,
American Journal of Math. 116(1994), 127--177.

\bibitem{Kac} V. Kac, Infinite dimencional Lie algebras, Cambridge University Press, 1990.

\bibitem{Keb} M.S. K\'eb\'e, ${O}$-alg\`ebres quantiques, C.R. Acad. Sci. Paris Ser. I Math. 322(1996) ,1--4. 

\bibitem{Keb1} M.S. K\'eb\'e, Sur la classification des ${O}$-alg\`ebres quantiques, J. Algebra, 212, N2(1999), 626--659.

\bibitem{Khar} V.K. Kharchenko, An algebra of skew primitive elements, Algebra i Logica,
v. 37, N2(1998), 181--224. English translation: Algebra and Logic, 37, N2(1998), 101--127 (arXiv:math.QA/0006077).

\bibitem{Kh3}  V.K. Kharchenko, A quantum analog of the  Poincar\`e-Birkhoff-Witt theorem, 
Algebra i Logika, 38, N4(1999), 476--507. 
English translation: Algebra and Logic, 38, N4(1999), 259--276 (arXiv:math.QA/0005101).

\bibitem{Kh2} V.K. Kharchenko, Skew primitive elements in Hopf algebras and related identities,
Journal of Algebra, 238(2001), 534--559.

\bibitem{Kh4} V.K. Kharchenko, A combinatorial approach to the quantifications of Lie algebras,
Pacific Journal of Mathematics, 203, N1(2002), 191--233.

\bibitem{Kh1} V.K. Kharchenko, Constants of coordinate differential calculi
defined by Yang--Baxter operators, Journal of Algebra, 267, N1(2003), 96--129. 

\bibitem{Kh5} V.K. Kharchenko, Quantum Lie algebras and related problems,
Proceedings of the Third International Algebra Conference, June 16-July 1, 2002,
pp. 67--114, Kluwer Academic Publishers, 2003. 

\bibitem{KA} V.K. Kharchenko, A. A. Alvarez, On the combinatorial rank of Hopf algebras,
Contemporary Mathematics, v. 376(2005), 299--308.

\bibitem{KhT} V.K. Kharchenko, PBW-bases of coideal subalgebras and a freeness theorem,
Transactions of the American Mathematical Society, in press. 

\bibitem{KhQ} V.K. Kharchenko, On right coideal subalgebras, preprint, arXiv:math.QA/0702504.

\bibitem{Koo} T. Koornwinder, Askey-Wilson polynomials as zonal spherical functions
on the $SU(2)$ quantum group, SIAM Journal Math. Anal. 24, N3(1993), 795--813.

\bibitem{Let1} G. Letzter, Subalgebras which appear in quantum Iwasawa decompositions,
Canadian Journal of Mathematics, 49, N6(1997), 1206--1223.

\bibitem{Let2} G. Letzter, Symmetric pairs for quantized enveloping algebras,
Journal of Algebra, 220, N2(1999), 729--767.

\bibitem{Let3} G. Letzter, Harish-Chandra modules for quantum symmetric pairs,
Representation Theory, An electronic Journal of AMS, 4 (1999), 64--96.

\bibitem{Let} G. Letzter, Coideal subalgebras and quantum symmetric pairs,
in: S. Montgomery, H.-J. Schneider (Eds.) New Directions in 
Hopf Algebras, MSRI Publications, 43(2002), 117--165.

\bibitem{Lot} M. Lothaire, Algebraic Combinatorics on Words, Cambridge Univ. Press, 2002.

\bibitem{Luz1} G. Lusztig, Finite-dimensional Hopf algebras arising from quantized
enveloping algebras, J. Amer. Math. Soc. 3(1)(1990) 257--296. 

\bibitem{Luz2} G. Lusztig, Introduction to Quantum Groups, in: Progress in Mathematics 110,
Birkh\"auser, Boston, 1993.

\bibitem{Mil1} A. Milinski, Actions of pointed Hopf algebras on prime algebras, Comm. Algebra, 23(1995), 313-333. 

\bibitem{Mil2} A. Milinski, Operationen punktierter Hopfalgebren auf primen Algebren (1995), Ph. D
thesis (M\"{u}nchen).

\bibitem{MS} A. Milinski, H.-J. Schneider, Pointed indecomposable Hopf algebras
over Coxeter groups, Contemp. Math. 267(2000), 215--236.  

\bibitem{Mul} A. E. M\"{u}ller, Some topics on Frobenious-Lustig kernels, Journal of Algebra,
206(1998), 642--681.


\bibitem{Mon} S. Montgomery, Hopf Algebras and Their Actions on Rings, 
CBMS, 82, AMS, Providence,1993. 

\bibitem{Nic}  W. Nichols, Bialgebras of type one, Comm. Algebra, 6(15)(1978), 1521--1552.

\bibitem{Nou} M. Noumi, Macdonald's symmetric polynomials as zonal spherical
functions on some quantum homogeneous spaces, Advances in
Mathematics, 123, N1(1996), 16--77.

\bibitem{NS} M. Noumi and T. Sugitani, Quantum symmetric spaces and related
$q$-orthogonal polynomials, In: Group Theoretical Methods in Physics
(ICGTMP) (Toyonaka, Japan, 1994), Word Sci. Publishing, River Edge, N.J. (1995), 28--40.

\bibitem{Rad} D. E. Radford, Hopf algebras with projection, Journal of Algebra, 92(1985), 322-347.

\bibitem{Res} N. Reshetikhin, Multiparameter quantum groups and twisted quasitriangular
Hopf algebras, Lett.Math. Phys. 20(1990), 331--335.

\bibitem{RTF} N. Reshetikhin, N. Takhtajan, L. Faddeev, Quantization of Lie groups and Lie algebras,
Leningrad Math. J., 1(1990), 193--225. 

\bibitem{Ros} M. Rosso, Quantum groups and quantum Shuffles,
Invent. Math., 113(2)(1998) 399--416.

\bibitem{Scha} P. Schauenburg, A characterization of the Borel-like subalgebras
of quantum enveloping algebras, Comm. in Algebra 24(1996), 2811--2823.

\bibitem{pSh1} A.I. Shirshov, On free Lie rings, Mat. Sb. 45, 87(2)(1958), 113-122. 

\bibitem{pSh2} A.I. Shirshov, Some algorithmic problems for Lie algebras, 
Sibirskii Math. J., 3(2)(1962), 292--296.

\bibitem{SV} S. Sinelshchikov, L. Vaksman, Harish-Chandra embedding and $q$-analogues
of bounded symmetric domains, Lecture Notes in Phys., 509(1998), 312--316. 

\bibitem{Skr} S. Skryabin, Projectivity and freeness over comodule algebras, 
Trans. Amer. Math. Soc., 359, N6(2007), 2597--2623.

\bibitem{Tak1} M. Takeuchi, Survey of braided Hopf algebras,
in: New Trends in Hopf Algebra Theory, Contemp. Math., vol. 267, AMS, Providence 
RI, 2000, pp. 301-324.

\bibitem{Tow} J. Towber, Multiparameter quantum forms of the enveloping
algebra $U_{gl_n}$ related to the Faddeev-Reshetikhin-Takhtajan $U(R)$
construction, J. Knot Theory Ramifications 4(5)(1995) 263--317.

\bibitem{Wor} S.L. Woronowicz, 
Differential calculus on compact matrix pseudogroups (quantum groups),
Comm. Math. Phys., 122(1989), 125--170.

\bibitem{Y} T.Yanai, Galois correspondence theorem for Hopf algebra actions, in:
Algebraic structures and their representations, 393-411, Contemporary Mathematics, vol. 376, AMS, Providence, RI, 2005.
\end{thebibliography}
\end{document}